\documentclass{article}[12pt]
\usepackage{amsfonts,amssymb}
\usepackage{psfig}

\usepackage{a4wide}
\usepackage{graphicx}
\usepackage{color}
\usepackage{pictex}
\usepackage{amsthm}
\usepackage{amsmath}
\usepackage{mathrsfs}
\usepackage{epsfig}

\newtheorem{lemma}{Lemma}[section]
\newtheorem{theorem}{Theorem}[section]
\newtheorem{proposition}{Proposition}[section]
\newtheorem{remark}{Remark}[section]
\newtheorem{definition}{Definition}[section]

\newtheorem{corollary}{Corollary}[section]

\begin{document}

\title{{The Birman-Wenzl-Murakami algebra, Hecke algebra and representations of $U_{q}(osp(1|2n))$}}

\author{Sacha C. Blumen \footnote{School of Mathematics and Statistics, University of Sydney, 2006, NSW, Australia. \newline e-mail: sachab@maths.usyd.edu.au}}

\date{2nd July, 2006}

\maketitle
\begin{abstract}
A representation of the Birman-Wenzl-Murakami algebra
$BW_{t}(-q^{2n},q)$ exists in the centraliser algebra
$\mathrm{End}_{U_{q}(osp(1|2n))}\big(V^{\otimes t}\big)$,
where $V$ is the fundamental $(2n+1)$-dimensional irreducible $U_{q}(osp(1|2n))$-module.
This representation is defined using 
permuted $R$-matrices acting on $V^{\otimes t}$.
A complete set of projections onto and intertwiners between 
irreducible $U_{q}(osp(1|2n))$-summands of $V^{\otimes t}$
exists via this representation, 
proving that $\mathrm{End}_{U_{q}(osp(1|2n))}\big(V^{\otimes t}\big)$ is 
generated by the set of permuting $R$-matrices acting on $V^{\otimes t}$.

We also show that a representation of the the Iwahori-Hecke algebra $H_{t}(-q)$ of type $A_{t-1}$
exists in the centraliser algebra 
$\mathrm{End}_{U_{q}(osp(1|2))}\left[\big(V^{+}_{1/2}\big)^{\otimes t}\right]$, 
where $V^{+}_{1/2}$ is a two-dimensional irreducible representation of $U_{q}(osp(1|2))$.
\end{abstract}

\section{Introduction}

The Birman-Wenzl-Murakami algebra, introduced independently 
by Murakami \cite{murakami} and Birman and Wenzl \cite{bw}
in connection with the study of the Kauffman knot polynomial, 
has been the subject of much research, including its relation with the centraliser algebras of
repeated tensor products of the fundamental irreducible
representations of the quantum algebras
$U_{q}(sp(2n))$ and $U_{q}(so(2n+1))$ \cite{cp,tw}.
Evidence has pointed towards a connection between
the Birman-Wenzl-Murakami algebra and representations of
$U_{q}(osp(1|2n))$ from a variety of sources, including the relation of
the Brauer algebra, which the Birman-Wenzl-Murakami algebra can be seen as a deformation of, with tensor products of representations of $U(osp(1|2n))$ \cite{benkart},
and the isomorphism between tensorial irreducible representations of
$U_{-q}(so(2n+1))$ and $U_{q}(osp(1|2n))$ \cite{zsuper}.

Let $\mathfrak{g}=osp(1|2n)$ and $V$ be the $(2n+1)$-dimensional irreducible
$U_{q}(\mathfrak{g})$-module.
In this paper we show that a subalgebra of the centraliser algebra 
$\mathrm{End}_{U_{q}(\mathfrak{g})}(V^{\otimes t})$
yields a representation of the Birman-Wenzl-Murakami algebra $BW_{t}(-q^{2n},q)$
and that this subalgebra, which is generated by
a representation of the Braid group $B_{t}$ on $t$ strings, 
is in fact equal to $\mathrm{End}_{U_{q}(\mathfrak{g})}(V^{\otimes t})$.

We also show a connection between the Hecke algebra $H_{t}(-q)$ of type
$A_{t-1}$ and repeated tensor powers of
representations of $U_{q}(osp(1|2))$ that has not, 
to our knowledge, appeared in the literature.

It may be thought that the connection between $BW_{t}(-q^{2n},q)$ and 
$\mathrm{End}_{U_{q}(\mathfrak{g})}(V^{\otimes t})$ is a trivial
consequence of the isomorphism between integrable tensorial representations of 
$U_{q}(osp(1|2n))$ and $U_{-q}(so(2n+1))$ shown by Zhang \cite{zsuper}. 
However, this isomorphism is moreso indicative of such a connection
and we define our representation of $BW_{t}(-q^{2n},q)$ without reference to
representations of any other quantum algebra.

The structure of this paper is as follows.
In Section \ref{section:preliminaries} we introduce the notation we use throughout this paper and necessary algebraic concepts, eg $\mathbb{Z}_{2}$-graded vector spaces.
In Section \ref{sec:algebraschapter2} we introduce  
the quantum superalgebra $U_{q}(\mathfrak{g})$.
In Section \ref{subsec:100} we discuss the finite-dimensional irreducible representations of $U_{q}(osp(1|2n))$, including the fundamental $(2n+1)$-dimensional irreducible module $V$.
We recall the isomorphism between finite-dimensional integrable tensorial representations of $U_{q}(osp(1|2n))$ and $U_{-q}(so(2n+1))$ and that
$V^{\otimes t}$ is completely reducible for each $t$.

In Section \ref{subsec:RRRR} we introduce $R$-matrices for representations of
$U_{q}(\mathfrak{g})$.
In Section \ref{sec:'central'element} we investigate the properties of two useful elements of a completion
$\overline{U}^{+}_{q}(\mathfrak{g})$ of $U_{q}(\mathfrak{g})$.
In Section \ref{sec:spectam(a)} we determine the spectral decomposition of the element $\check{\cal{R}}_{V,V}$,
which is a representation of a Braid group generator in 
$\mathrm{End}_{U_{q}(\mathfrak{g})}(V \otimes V)$.
In Section \ref{eq:theXfactor(a)} we show that there is
a representation of  $BW_{t}(-q^{2n},q)$ in the complex algebra ${\cal{C}}_{t}$
generated by the $\check{\cal{R}}_{V,V}$-matrices 
acting on the $i^{th}$ and $(i+1)^{st}$ tensor powers of
$V^{\otimes t}$ for $i=1, \ldots, t-1$.
In Section \ref{eq:caseydonovan(a)} we recall the ideas of Bratteli diagrams and path algebras.
In Section \ref{subsec:projectontome} we construct 
projections from $V^{\otimes t}$
onto all the irreducible $U_{q}(\mathfrak{g})$-submodules of $V^{\otimes t}$
using elements of $\overline{U}^{+}_{q}(\mathfrak{g})$ from Section \ref{sec:'central'element}.
In Section \ref{subsec:rhapsodyinred} we construct a complete set of matrix units in
$\mathrm{End}_{U_{q}(\mathfrak{g})}(V^{\otimes t})$ 
from  matrix units in $BW_{t}(-q^{2n},q)$
and prove that the algebra ${\cal{C}}_{t}$ is in fact equal to 
$\mathrm{End}_{U_{q}(\mathfrak{g})}(V^{\otimes t})$.
In Section \ref{sec:repoftheHeckealgebra} we show that there exists a representation of
$H_{t}(-q)$ in $\mathrm{End}_{U_{q}(osp(1|2))}\left[\big(V^{+}_{1/2}\big)^{\otimes t}\right]$.

The work in this paper elucidates part of the author's Ph.D thesis
\cite{blumen05} together with some further results. 
In writing this paper, the author became aware of the work of Lehrer and Zhang on `strongly multiplicity free modules' of quantum algebras \cite{lehrerzhang}.
Their work may extend in a natural way to cover representations of quantum superalgebras, including part of the work in this paper, but we have not investigated this.

\section{Preliminaires}
\label{section:preliminaries}

Throughout, we fix $\mathbb{Z}$ to be the integers,
$\mathbb{Z}_{+}$ to be the non-negative integers and 
$\mathbb{N} = \mathbb{Z} \backslash \{0\}$.
Let $q \in \mathbb{C}$ be not equal to $0$ or $1$, then for each $k \in \mathbb{Z}_{+}$ we define
$$[k]_{q} = \frac{1-q^{k}}{1-q}, \hspace{5mm}
[k]_{q}! = [k]_{q} [k-1]_{q} \cdots [1]_{q}, \hspace{5mm}
[0]_{q}! = 1.$$

A vector space $W$ is said to be $\mathbb{Z}_{2}$-graded if $W$ 
can be expressed as a direct sum of vector subspaces: $W = W_{0} \oplus W_{1}$.
We say that an element $w \in W$ is homogeneous if $w \in W_{0} \cup W_{1}$.
Each homogeneous element $w \in W$ has a grading: 
$w$ has an even grading (resp. an odd grading) if $w \in W_{0}$ 
(resp. $w \in W_{1}$).
We denote the grading of the homogeneous element 
$w \in W$ by $[w] = i$, where $w \in W_{i}$.

Let $V$ and $W$ be $\mathbb{Z}_{2}$-graded vector spaces.
The graded permutation operator $P: V \otimes W \rightarrow W \otimes V$ is defined for
homogeneous $v \in V$ and $w \in W$ by 
$$P(v \otimes w) = (-1)^{[v][w]} w \otimes v,$$ and is extended to inhomogeneous elements of $V$ and $W$ by linearity.

Let $V = \bigoplus_{i \in \mathbb{Z}_{2}} V_{i}$ and $W = \bigoplus_{j \in \mathbb{Z}_{2}}W_{j}$ be 
$\mathbb{Z}_{2}$-graded vector spaces. 
A $\mathbb{C}$-linear map
$\psi: V \rightarrow W$ is said to be homogeneous of degree $k \in \mathbb{Z}_{2}$ if
$\psi(V_{i}) \subseteq W_{i+k}$ for each $i \in \mathbb{Z}_{2}$.
The vector space $\mathrm{Hom}_{\mathbb{C}}(V,W)$ of $\mathbb{C}$-linear maps from $V$ to $W$ admits a $\mathbb{Z}_{2}$-gradation: 
$\mathrm{Hom}_{\mathbb{C}}(V,W) = \bigoplus_{i \in \mathbb{Z}_{2}} \mathrm{Hom}_{\mathbb{C}}(V,W)_{i}$
where $\mathrm{Hom}_{\mathbb{C}}(V,W)_{i} = 
\{ \psi \in \mathrm{Hom}_{\mathbb{C}}(V,W) | \ \mbox{$\psi$ is homogeneous of degree $i$} \}$.
We fix $\mathrm{End}_{\mathbb{C}}(V) = \mathrm{Hom}_{\mathbb{C}}(V,V)$.

The dual space to $V$ is $V^{*} =  \mathrm{Hom}_{\mathbb{C}}(V,\mathbb{C})$
where we regard
$\mathbb{C}$ as a $\mathbb{Z}_{2}$-graded vector space with zero odd subspace. 
$V^{*}$ is $\mathbb{Z}_{2}$-graded.

One can define $\mathbb{Z}_{2}$-graded algebras, $\mathbb{Z}_{2}$-graded co-algebras, and
$\mathbb{Z}_{2}$-graded Hopf algebras in a standard way (see \cite{gzb} for further information).

Let $A$ be a $\mathbb{Z}_{2}$-graded algebra. 
A homomorphism of $\mathbb{Z}_{2}$-graded $A$-modules is by definition a homomorphism of
the $A$-modules as well as of the underlying $\mathbb{Z}_{2}$-graded vector spaces.
Each such homomorphism is $A$-linear and homogeneous of degree $0$.
Let $V = V_{0} \oplus V_{1}$ be a $\mathbb{Z}_{2}$-graded $A$-module.
If $V' = V'_{0} \oplus V'_{1}$ with $V'_{0} = V_{1}$ and $V'_{1} = V_{0}$ then we regard
$V$ and $V'$ as not being isomorphic.
 
Throughout this paper we fix $\delta_{ij} = \left\{ \begin{array}{ll}
1, & \ \mbox{if } i=j, \\
0, & \ \mbox{otherwise}.
\end{array} \right.$

\begin{section}{The quantum superalgebra $U_{q}(osp(1|2n))$}
\label{sec:algebraschapter2}
%\markright{\text{The quantum superalgebra $U_{q}(osp(1|2n))$}}

In this section we introduce the quantum superalgebra $U_{q}(osp(1|2n))$ \cite{kt,y,z1}.
Let us begin by describing the root system of $osp(1|2n)$.
Let $H^{*}$ be a vector space over $\mathbb{C}$ with a basis
$\left\{ \epsilon_{i} | \ 1 \leq i \leq n \right\}$ and let 
\begin{equation}
\label{eq:starstar(1)}
( \cdot, \cdot): H^{*} \times H^{*} \rightarrow \mathbb{C},
\end{equation}
be a non-degenerate bilinear form defined by
$(\epsilon_{i}, \epsilon_{j})=\delta_{ij}$.  
The set of simple roots of $osp(1|2n)$ is $\left\{ \alpha_{i} | \ 1 \leq i \leq n \right\}$ where
\begin{eqnarray*}
\alpha_{i} & = & \left\{
\begin{array}{ll}
\epsilon_{i}-\epsilon_{i+1}, &  i=1, \ldots, n-1,   \\
\epsilon_{n},                &  i=n,
\end{array} \right.
\end{eqnarray*}
which forms another basis of $H^{*}$.

The set of the positive roots of $osp(1|2n)$ is
$$\Phi^{+}=\left\{
\epsilon_{i} \pm \epsilon_{j}, \epsilon_{k}, 2\epsilon_{k}| \ \ 1 \leq i < j \leq n, \ \ 1 \leq k \leq n \right\},$$
and we further define
the subsets of positive even roots (resp. positive odd roots) $\Phi^{+}_{0}$ (resp. $\Phi^{+}_{1}$) of $\Phi^{+}$ and
also a subset $\overline{\Phi}_{0}^{+} \subset \Phi^{+}_{0}$ by
$$\begin{array}{ll}
\Phi^{+}_{0} = \{\epsilon_{i} \pm \epsilon_{j}, 2\epsilon_{k}| \ \ 1 \leq i < j \leq n, \ \ 1 \leq k \leq n \}, &
\hspace{10mm} \Phi^{+}_{1} = \{ \epsilon_{k} | \ \ 1 \leq k \leq n\}, \\
\overline{\Phi}_{0}^{+} = \{\alpha \in \Phi^{+}_{0} | \ \ \alpha/2 \notin \Phi^{+}_{1}\}.  &
\end{array}$$
The set of negative roots of $osp(1|2n)$ is
$\Phi^{-} = -\Phi^{+}$, and $\Phi = \Phi^{+} \cup \Phi^{-}$ is the set of all roots of $osp(1|2n)$.

We denote by $2\rho \in H^{*}$ the graded sum of the positive roots of $osp(1|2n)$:
$$2 \rho = \sum_{\alpha \in \Phi^{+}_{0}} \alpha - \sum_{\beta \in \Phi^{+}_{1}} \beta.$$
Explicitly, $2\rho=\sum_{i=1}^{n}(2n-2i+1) \epsilon_{i}$. 
The element $2\rho$ satisfies
$\left(2 \rho,\alpha_{i}\right)=(\alpha_{i},\alpha_{i})$ for each $1 \leq i \leq n$.

The Cartan matrix of $osp(1|2n)$ is $A=\left(a_{ij}\right)_{i,j=1}^{n}$ where
$a_{ij}=2\left(\alpha_{i}, \alpha_{j}\right) / \left(\alpha_{i}, \alpha_{i}\right)$:
$$A=\left( \begin{array}{rrrrrr}
2 & -1 & 0 & \cdots & 0 & 0 \\
-1 & 2 & -1 & \cdots & 0 & 0 \\
0 & -1 & 2 & \cdots & 0 & 0 \\
\vdots & \vdots & \vdots & \ddots & \vdots & \vdots \\
0 & 0 & 0 & \cdots & 2 & -1 \\
0 & 0 & 0 & \cdots & -2 & 2
\end{array} \right).$$

The Lie superalgebra $\mathfrak{g}=osp(1|2n)$ over $\mathbb{C}$ 
can be defined in terms of a Serre presentation
with generators $\{E_{i}, F_{i}, H_{i} | \ 1 \leq i \leq n \}$
subject to the relations
\begin{eqnarray}
\left[E_{i},F_{j}\right]=\delta_{ij}H_{i}, & \hspace{10mm} & \left[H_{i},H_{j}\right]=0, \hspace{20mm}  \forall \ i,j, \nonumber \\
\left[H_{i},E_{j}\right]=(\alpha_{i}, \alpha_{j}) E_{j}, 
& \hspace{10mm} & \left[H_{i},F_{j}\right]=-(\alpha_{i}, \alpha_{j})F_{j}, \hspace{5mm} \forall \ i,j, \nonumber \\
\label{eq:heylittlepiggy(1)}
(\mbox{ad} \ E_{i})^{1-a_{ij}}E_{j}=0, & \hspace{10mm} &  (\mbox{ad} \ F_{i})^{1-a_{ij}}F_{j}=0, 
\hspace{10mm}   \forall \ i \neq j, 
\end{eqnarray}
where $[ \cdot, \cdot ]$ denotes the $\mathbb{Z}_{2}$-graded Lie bracket and the adjoint operation is 
$(\mbox{ad} \ a)b= [a,b]$.
The $\mathbb{Z}_{2}$-grading of the generators is
$$[E_{i}] = [F_{i}] = [H_{j}] = 0, \hspace{8mm} [E_{n}] = [F_{n}] = 1, 
\hspace{8mm}  1 \leq i \leq n-1, \hspace{5mm} 1 \leq j \leq n. $$
\noindent

The universal enveloping algebra $U(\mathfrak{g})$ of $\mathfrak{g}$
 is a unital associative 
$\mathbb{Z}_{2}$-graded algebra,
which may be considered as being generated by 
$\{E_{i}, F_{i}, H_{i} | \ 1 \leq i \leq n \}$ subject to relations that are formally the same as
(\ref{eq:heylittlepiggy(1)})
but with the bracket $[ \cdot, \cdot]$ interpreted as a
$\mathbb{Z}_{2}$-graded commutator
$[ \cdot, \cdot ]: U(\mathfrak{g}) \times U(\mathfrak{g}) \rightarrow U(\mathfrak{g})$ defined by
\begin{equation}
\label{eq:starstar(2)}
[X, Y] = XY -(-1)^{[X][Y]} YX.
\end{equation}
If each of two elements $X,Y \in U(\mathfrak{g})$ has a grading, 
then the grading of $XY \in U(\mathfrak{g})$ is defined by
$$[XY]= \big( [X] + [Y] \big) \bmod{2}.$$
The graded commutator $[X, Y]$ of any two homogeneous elements of $U(\mathfrak{g})$ is defined by
(\ref{eq:starstar(2)}) and is extended to inhomogeneous elements of $U(\mathfrak{g})$ by linearity.
Note that $U(\mathfrak{g}) \otimes U(\mathfrak{g})$ has a natural 
$\mathbb{Z}_{2}$-graded associative algebra structure, with the grading defined for
homogeneous $X, Y \in U(\mathfrak{g})$ by
$$[X \otimes Y] = \big( [X] + [Y] \big) \bmod{2}.$$

The quantum superalgebra $U_{q}(\mathfrak{g})$ is some ``$q$-deformation'' of $U(\mathfrak{g})$.
We describe its Jimbo version here.
\begin{definition}
The quantum superalgebra $U_{q}(\mathfrak{g})$ over $\mathbb{C}$ 
is an associative $\mathbb{Z}_{2}$-graded unital
algebra generated by
$\{e_{i},f_{i},K_{i},K_{i}^{-1}| \ 1 \leq i \leq n\}$ subject to the relations
$$\left[e_{i},f_{j}\right]= \delta_{ij}\frac{K_{i}-K_{i}^{-1}}{q-q^{-1}},$$
$$K_{i}e_{j}K_{i}^{-1}= q^{\left(\alpha_{i}, \alpha_{j}\right)}e_{j}, \hspace{10mm}
K_{i}f_{j}K_{i}^{-1}=q^{-\left(\alpha_{i}, \alpha_{j}\right)}f_{j},$$
$$\big[K_{i}^{\pm 1},K_{j}^{\pm 1}\big]=\big[K_{i}^{\pm 1},K_{j}^{\mp 1}\big]=0, 
\hspace{10mm} K_{i}^{\pm 1}K_{i}^{\mp 1}=1,$$
\begin{equation}
\label{equat:two}
(\mathrm{ad}_{q}e_{i})^{1-a_{ij}}e_{j}=0, \hspace{10mm} (\mathrm{ad}_{q}f_{i})^{1-a_{ij}}f_{j}=0, \hspace{10mm} \forall i,j,
\end{equation}
where $0 \neq q \in \mathbb{C}$ and $q^{2} \neq 1$, and the adjoint actions are defined by
$$(\mathrm{ad}_{q}e_{i})X  = e_{i}X -(-1)^{[e_{i}][X]} K_{i}      X K_{i}^{-1} e_{i},$$
$$(\mathrm{ad}_{q}f_{i})X =  f_{i}X -(-1)^{[f_{i}][X]} K_{i}^{-1} X K_{i}      f_{i},$$
for all homogeneous $X \in U_{q}(\mathfrak{g})$.  
The grading of each generator is even except for $e_{n}$ and $f_{n}$, which are odd.
In (\ref{equat:two}), the bracket $[ \cdot, \cdot]$ is as defined in (\ref{eq:starstar(2)}).
\end{definition}
For each $\beta = \sum_{i=1}^{n} m_{i} \alpha_{i}$ where $m_{i} \in \mathbb{Z}$, we define
$K_{\beta} = \prod_{i=1}^{n} \left( K_{i} \right)^{m_{i}}$.

As is well known,
there exists a $\mathbb{Z}_{2}$-graded Hopf algebra structure on $U_{q}(\mathfrak{g})$ with 
the co-multiplication $\Delta$, the co-unit $\epsilon$, and the antipode $S$ defined on each generator by
$$
\Delta(e_{i}) =  e_{i} \otimes K_{i} + 1 \otimes e_{i},   \hspace{5mm}
\Delta(f_{i})  = f_{i} \otimes 1 + K_{i}^{-1} \otimes f_{i}, 
\hspace{5mm} \Delta\left(K_{i}^{\pm 1}\right) =  K_{i}^{\pm 1} \otimes K_{i}^{\pm 1}.$$
$$\epsilon(e_{i}) = \epsilon(f_{i}) =  0, \hspace{5mm} \epsilon\left(K_{i}^{\pm 1}\right)=\epsilon(1) = 1.$$
$$ S(e_{i}) = -e_{i}K^{-1}_{i},      \hspace{5mm}
S(f_{i}) = -K_{i}f_{i},           \hspace{5mm}
S(K_{i}^{\pm 1}) = K_{i}^{\mp 1}.$$

There are a number of quite useful subalgebras of $U_{q}(\mathfrak{g})$.
We define $U_{q}(\mathfrak{b}_{+})$ (resp. $U_{q}(\mathfrak{b}_{-})$)
 to be the subalgebra of $U_{q}(\mathfrak{g})$ 
generated by $\left\{ e_{i}, K^{\pm 1}_{i} | \ 1 \leq i \leq n \right\}$  
(resp. $\left\{ f_{i}, K^{\pm 1}_{i} | \ 1 \leq i \leq n \right\}$).

We shall sometimes refer to the quantum superalgebra in the sense of Drinfel'd which we denote by $U_{h}(\mathfrak{g})$.
This is a $\mathbb{Z}_{2}$-graded Hopf algebra over the ring $\mathbb{C}[[h]]$ in an indeterminate $h$, completed with respect to the $h$-adic topology \cite{kt}.
Fix $q = e^{h}$.
The algebra $U_{h}(\mathfrak{g})$ 
is generated by $\{E_{i}, F_{i}, H_{i} | \ 1 \leq i \leq n \}$ with relations
\begin{eqnarray*}
\left[ E_{i}, F_{j} \right] = \delta_{ij} \frac{e^{h H_{i}} - e^{-h H_{i}}}{e^{h}-e^{-h}}, & \hspace{5mm} & [H_{i},H_{j}]=0, \\
\left[H_{i}, E_{j} \right] = (\alpha_{i}, \alpha_{j}) E_{j}, & \hspace{5mm} & [H_{i}, F_{j}] = -(\alpha_{i},\alpha_{j}) F_{j},  \\
(\mathrm{ad}_{q} E_{i})^{1-a_{ij}}E_{j}=0, & \hspace{5mm} & 
(\mathrm{ad}_{q} F_{i})^{1-a_{ij}}F_{j}=0, \hspace{5mm} \forall i,j,
\end{eqnarray*}
where the adjoint functions are defined by
$$(\mathrm{ad}_{q}E_{i})X  = E_{i}X-(-1)^{[E_{i}][X]} e^{h H_{i}}  X e^{-h H_{i}} E_{i},$$
$$(\mathrm{ad}_{q}F_{i})X  = F_{i}X-(-1)^{[F_{i}][X]} e^{-h H_{i}} X e^{h H_{i}}  F_{i}.$$

A $\mathbb{Z}_{2}$-graded quasitriangular Hopf algebra $A$ is a 
$\mathbb{Z}_{2}$-graded Hopf algebra that admits a universal $R$-matrix, which is
an invertible even element $R \in A \otimes A$ satisfying the relations
\begin{eqnarray}
R \Delta(x) & = & \Delta'(x) R, \hspace{5mm} \forall x \in A \label{eq:definingrelationsR1} \\
(\Delta \otimes \mathrm{id}) R & = & R_{13} R_{23} \label{eq:definingrelationsR2} \\
(\mathrm{id} \otimes \Delta) R & = & R_{13} R_{12} \label{eq:definingrelationsR3} 
\end{eqnarray}
where upon writing $R = \sum_{s} \alpha_{s} \otimes \beta_{s}$, we write
$R_{12} = R \otimes \mathrm{id}$, 
$R_{13} = \sum_{s} \alpha_{s} \otimes \mathrm{id} \otimes \beta_{s}$ and
$R_{23} = \mathrm{id} \otimes R$, and also write 
$\Delta'(x) = P (\Delta(x))$ where $P$ is the graded permutation operator.
It is important to note that $U_{h}(\mathfrak{g})$ admits a universal $R$-matrix \cite{kt} and is a $\mathbb{Z}_{2}$-graded ribbon Hopf algebra \cite{zg}.

Given a finite-dimensional $U_{q}(\mathfrak{g})$-module $W$, we define 
the quantum supertrace of $f \in \mathrm{End}_{\mathbb{C}}(W)$ to be
$$\mathrm{str}_{q}(f) = \mathrm{str}\big(\pi_{W}(K_{2\rho}) \cdot f \big).$$

\end{section}

\begin{section}{Finite dimensional irreducible $U_{q}(osp(1|2n))$-modules}
\label{subsec:100}
%\markright{\text{Finite dimensional irreducible $U_{q}(osp(1|2n))$-modules}}

Given a $\mathbb{Z}_{2}$-graded algebra $A$, we will denote by
$V_{\lambda}$ an $A$-module labelled by $\lambda \in I$ for some index set $I$, and
we denote by $\pi_{\lambda}$ the representation of $A$ afforded by $V_{\lambda}$.

In this section we assume that $q$ is non-zero and not a root of unity;
in this case the representation theory of $U_{q}(\mathfrak{g})$
is completely understood \cite{ab,zsuper,zou}.

We say that an element $\lambda \in H^{*}$ is {\emph{integral}} if
$$l_{i} = \frac{2(\lambda,\alpha_{i})}{(\alpha_{i},\alpha_{i})} \in \mathbb{Z}, 
\ \ \forall i < n, \hspace{10mm} l_{n} = \frac{(\lambda,\alpha_{n})}{(\alpha_{n},\alpha_{n})} \in 
\mathbb{Z},$$
where $(\cdot, \cdot): H^{*} \times H^{*} \rightarrow \mathbb{C}$ is the
bilinear form from (\ref{eq:starstar(1)}) and we let
${\cal{P}}$ be the set of all integral elements of $H^{*}$.  
Furthermore, we say that an element $\lambda \in {\cal{P}}$ is {\emph{integral dominant}} if 
$l_{i} \in \mathbb{Z}_{+}$ for all $i$ and
denote by ${\cal{P}}^{+}$ the set of all integral dominant elements of $H^{*}$.

Zou investigated the representation theory of the quantum superalgebra 
$U_{q}(\mathfrak{g})$ over the quotient field 
$\mathbb{C}(v)$ for an indeterminate $v$ \cite{zou}, which is related to $q$ 
via $q=v^{2}$.  
Zou's results can be adapted to our setting where we take $q$ to be generic.  
Let $\sqrt{q}$ be any square root of the complex number $q$.  
Call a $U_{q}(\mathfrak{g})$-module $V$ {\emph{integrable}} if $V$ 
is a direct sum of its weight spaces and
if $e_{i}$ and  $f_{i}$ act as locally nilpotent endomorphisms of $V$ for each 
$i=1, \ldots, n$.

Let $\overline{V}(\omega)$ be a highest weight $U_{q}(\mathfrak{g})$-module 
with highest weight vector $v$ satisfying 
$K_{i} v = \omega_{i} v$, $\omega_{i} \in \mathbb{C}$, for each $i=1, \ldots, n$, then 
$\overline{V}(\omega)$ has a unique maximal proper
$U_{q}(\mathfrak{g})$-submodule $\overline{M}(\omega)$, and the quotient
$$V(\omega) = \overline{V}(\omega) / \overline{M}(\omega)$$
is an irreducible $U_{q}(\mathfrak{g})$-module with highest weight
$\omega = (\omega_{1}, \omega_{2}, \ldots, \omega_{n})$.
Theorem 3.1 of \cite{zou} can be stated in our setting as follows.
\begin{theorem}
The irreducible highest weight $U_{q}(\mathfrak{g})$-module $V(\omega)$, 
with $\omega = (\omega_{1}, \omega_{2}, \ldots, \omega_{n})$, is integrable if and only if
$$\omega_{i} = \zeta_{i} q^{m_{i}}, \hspace{5mm} 1 \leq i \leq n-1,$$
where $m_{i} \in \mathbb{Z}_{+}$, $\zeta_{i}^{2}=1$, and
$$ \omega_{n} = \left\{ \begin{array}{ll}
\pm q^{m}, & \mbox{ if } m \in \mathbb{Z}_{+}, \\
\pm\sqrt{-1} \ q^{m}, & \mbox{ if } m \in \mathbb{Z}_{+} + \frac{1}{2}.
\end{array} \right.$$
\end{theorem}
\noindent
Note that every finite dimensional integrable $U_{q}(\mathfrak{g})$-module is semisimple
\cite[Sec. 5]{zou}.

If $\omega = (q^{m_{1}}, q^{m_{2}}, \ldots, q^{m_{n}})$ with $m_{i} \in \mathbb{Z}_{+}$
for each $i$, there exists an irreducible $U(\mathfrak{g})$-module 
$V(\omega)_{\mathfrak{g}}$
 with highest weight $\lambda \in {\cal{P}}^{+}$ satisfying $(\lambda, \alpha_{i}) = m_{i}$ 
 for each $i$, and
 $V(\omega)$ and $V(\omega)_{\mathfrak{g}}$ have the same weight space decompositions \cite{zsuper}.
 In this paper, we are more interested in these irreducible $U_{q}(\mathfrak{g})$-modules, and
 for each $\lambda \in {\cal{P}}^{+}$ 
 we let $V_{\lambda}$ denote $V(\omega)$ and call $\lambda$ 
 the highest weight of $V_{\lambda}$.

For convenience we fix the grading of the highest weight vector of the 
finite dimensional irreducible
$U_{q}(\mathfrak{g})$-module $V_{\lambda}$ with integral dominant highest weight
$\lambda = \sum_{i=1}^{n} \lambda_{i} \epsilon_{i} \in {\cal{P}}^{+}$ to be
even (resp. odd) if $\sum_{i=1}^{n} \lambda_{i}$ is even (resp. odd).

The following two lemmas restate some important results from \cite{zsuper}.
\begin{lemma}
\label{lem:quantumospandquantumsomrepresentations}

Let $\mu \in {\cal{P}}^{+}$ and let $V_{\mu}$ be a finite dimensional irreducible 
$U_{q}(osp(1|2n))$-module with highest weight $\mu$.
Then $V_{\mu}$ is also an irreducible $U_{-q}(so(2n+1))$-module with highest weight $\mu$.
The weight space decompositions of $V_{\mu}$ as a $U_{q}(osp(1|2n))$-module and as a
$U_{-q}(so(2n+1))$-module are the same.
Furthermore, the quantum dimension of $V_{\mu}$ as a $U_{-q}(so(2n+1))$-module
equals the quantum superdimension of $V_{\mu}$ as a $U_{q}(osp(1|2n))$-module.
\end{lemma}

Let $\mu, \nu \in {\cal{P}}^{+}$ and let $V_{\mu}$, $V_{\nu}$ be irreducible
$U_{-q}(so(2n+1))$-modules with highest weights
$\mu$ and $\nu$, respectively, then it is well known that 
$V_{\mu} \otimes V_{\nu}$ is completely reducible into a direct sum
of irreducible $U_{-q}(so(2n+1))$-submodules:
\begin{equation}
\label{eq:decompositionmultiplicitiesofquantumsomodules}
V_{\mu} \otimes V_{\nu} \cong \bigoplus_{\lambda \in {\cal{P}}^{+}} 
\left( \mathbb{C}^{m(\lambda,\mu,\nu)} \otimes V_{\lambda} \right),
\end{equation}
where $m(\lambda, \mu, \nu) \geq 0$ is the number of copies of irreducible $U_{-q}(so(2n+1))$-submodules
of $V_{\mu} \otimes V_{\nu}$ in the decomposition isomorphic to $V_{\lambda}$. Then:

\begin{lemma}
\label{lem:decompositionofthetensorproductoftwoirreduciblequantumospmodules}

Let $\mu, \nu \in {\cal{P}}^{+}$ and let $V_{\mu}$, $V_{\nu}$ be irreducible
$U_{q}(osp(1|2n))$-modules with highest weights
$\mu$ and $\nu$, respectively, then 
$V_{\mu} \otimes V_{\nu}$ is completely reducible into a direct sum
of irreducible $U_{q}(osp(1|2n))$-submodules:
$$V_{\mu} \otimes V_{\nu} \cong \bigoplus_{\lambda \in {\cal{P}}^{+}} 
\left( \mathbb{C}^{m(\lambda,\mu,\nu)} \otimes V_{\lambda} \right),$$
where $m(\lambda, \mu, \nu) \geq 0$ is the number of copies of irreducible $U_{q}(osp(1|2n))$-submodules
of $V_{\mu} \otimes V_{\nu}$ in the decomposition
isomorphic to $V_{\lambda}$, where the constants $m(\lambda,\mu,\nu)$ are the
same as in (\ref{eq:decompositionmultiplicitiesofquantumsomodules}).

\end{lemma}

Let $W$ be a $U_{q}(\mathfrak{g})$-module with homogeneous basis $\{w_{i}\}_{i \in I}$.
Let $\{w_{i}^{*}\}$ be a basis of $W^{*}$ where $[w_{i}^{*}] = [w_{i}]$
and the dual space pairing is $\langle w_{i}^{*}, w_{j} \rangle = \delta_{ij}$.
Then 
$W^{*}$ is the dual $U_{q}(\mathfrak{g})$-module to $W$ with the action of $U_{q}(\mathfrak{g})$ given by
$\langle a w_{i}^{*}, w_{j} \rangle = 
(-1)^{[a][w_{i}^{*}]}\langle w_{i}^{*}, S(a) w_{j} \rangle$, 
$\forall a \in U_{q}(\mathfrak{g})$.

We now introduce the fundamental (or vector) $U_{q}(\mathfrak{g})$-module $V$ 
in the next lemma adapted from \cite{lz}, which we state without proof. 
\begin{lemma}
\label{lem:fundamentaldimensional}
There exists a $(2n+1)$-dimensional irreducible $U_{q}(osp(1|2n))$-module $V=V_{\epsilon_{1}}$
with highest weight $\epsilon_{1}$.
This module admits a basis 
$\{ v_{i} | \ -n \leq i \leq n\}$ with $v_{1}$ being the highest weight vector.
The actions of the generators of $U_{q}(osp(1|2n))$ on the basis elements are 
$$f_{i} v_{i} = v_{i+1},  \hspace{5mm}
f_{n} v_{n} = v_{0},  \hspace{5mm}
f_{n} v_{0} = v_{-n},  \hspace{5mm}
f_{i} v_{-i-1} = v_{-i},  
$$
$$
e_{i} v_{i+1} = v_{i}, \hspace{5mm}
e_{n} v_{0}   = v_{n}, \hspace{5mm}
e_{n} v_{-n}  = -v_{0}, \hspace{5mm}
e_{i} v_{-i}  = v_{-i-1}, $$  
$$
K_{j}^{\pm 1} v_{k} = q^{\pm (\alpha_{j},\epsilon_{k})} v_{k},$$
where $1 \leq i< n$,  $1 \leq j \leq n$, $-n \leq k \leq n$, and we fix
$\epsilon_{0}=0$, and $\epsilon_{-i} = - \epsilon_{i}$.
All remaining actions are zero.
\end{lemma}
Note that the highest weight vector $v_{1}$ of $V$ has an odd grading.

\begin{proposition}
\label{prop:mooV}
There exists a $U_{q}(\mathfrak{g})$-invariant, non-degenerate bilinear form
$\langle \langle \ , \ \rangle \rangle : V \times V \rightarrow \mathbb{C}$.
Thus the dual $U_{q}(\mathfrak{g})$-module of $V$ is isomorphic to $V$.
\end{proposition}
\begin{proof}
Let $\{ v_{i} | \ -n \leq i \leq n \}$ be the basis of $V$ given in Lemma \ref{lem:fundamentaldimensional}.
Now define a non-degenerate $\mathbb{C}$-bilinear form
$\langle \langle \ , \ \rangle \rangle : V \times V \rightarrow \mathbb{C}$ by
$$
\begin{array}{rcll}
\langle \langle v_{1}, v_{-1} \rangle \rangle & = & 1, & \\
\langle \langle v_{i}, v_{-i} \rangle \rangle & = & 
-q^{-1} \langle \langle v_{i-1}, v_{-(i-1)} \rangle \rangle, &  2 \leq i \leq n, \\
\langle \langle v_{0}, v_{0} \rangle \rangle & = & q^{-1} \langle \langle v_{n}, v_{-n} \rangle \rangle, & \\
\langle \langle v_{-n}, v_{n} \rangle \rangle & = & - \langle \langle v_{0}, v_{0} \rangle \rangle, & \\
\langle \langle v_{-j}, v_{j} \rangle \rangle & = & 
-q^{-1} \langle \langle v_{-(j+1)}, v_{j+1} \rangle \rangle, &  1 \leq j \leq n-1, \\
\langle \langle v_{k}, v_{l} \rangle \rangle & = & 0, & \mbox{if } k + l \neq 0.
\end{array}
$$
A direct calculation shows that
$$\langle \langle x v_{i}, v_{j} \rangle \rangle = 
(-1)^{[x][v_{i}]}  \langle \langle v_{i}, S(x) v_{j} \rangle \rangle,  
\hspace{5mm} \forall x \in U_{q}(\mathfrak{g}), \hspace{5mm} v_{i}, v_{j} \in V,$$
thus proving the $U_{q}(\mathfrak{g})$-invariance of the bilinear form.
This form identifies $V$ with its dual module.
\end{proof}

Let us discuss in more detail the dual module of $V$.
Recall the definition of the dual $U_{q}(\mathfrak{g})$-module $V^{*}$ to $V$.
Let $\{v_{i}^{*} | \ -n \leq i \leq n\}$ be a basis of $V^{*}$ such that
$\langle v_{i}^{*}, v_{j} \rangle = \delta_{ij}$ and $[v_{i}^{*}]=[v_{i}]$ where
$\langle \ , \  \rangle: V^{*} \times V \rightarrow \mathbb{C}$ is the dual space pairing.
Now define a homogeneous bijection $T \in \mathrm{Hom}_{\mathbb{C}}(V, V^{*})$ 
of degree $0$ for all $1 \leq i \leq n$ by
\begin{equation}
\label{eq:johnfaulkner(20)}
T:  v_{i}  \mapsto (-1)^{i-1} q^{-(i-1)}v_{-i}^{*}, \hspace{5mm}
       v_{0}  \mapsto (-1)^{n-1} q^{-n} v_{0}^{*}, \hspace{5mm}
       v_{-i} \mapsto (-1)^{i} q^{-(2n-i)} v_{i}^{*}.
\end{equation}
A direct calculation shows that
this map is an element of $\mathrm{Hom}_{U_{q}(\mathfrak{g})}(V, V^{*})$ and that it satisfies
$ \langle T(v_{i}), v_{j} \rangle = \langle \langle v_{i}, v_{j} \rangle \rangle$
for all $v_{i}, v_{j} \in V$.

\end{section}

\begin{section}{$R$-matrices for representations of $U_{q}(osp(1|2n))$}
\label{subsec:RRRR}
%\markright{\text{$R$-matrices for representations of $U_{q}(osp(1|2n))$}}

Drinfel'd's quantum superalgebra $U_{h}(\mathfrak{g})$ 
admits a universal $R$-matrix \cite{kt,y}.
We will show that there does not exist an element of $U_{q}(\mathfrak{g})$
that corresponds to the universal $R$-matrix of $U_{h}(\mathfrak{g})$ in any obvious way.
However, there exists a completion 
$\overline{U}_{q}^{+}(\mathfrak{g})$ of $U_{q}(\mathfrak{g})$ such
that one of the multiplicative factors of the universal $R$-matrix of $U_{h}(\mathfrak{g})$
maps to an element $\widetilde{R}$ of $\overline{U}_{q}^{+}(\mathfrak{g})$.
Although $\widetilde{R}$ is not an element of $U_{q}(\mathfrak{g})$, 
only a finite number of its terms act as non-zero endomorphisms on each 
tensor product of finite dimensional irreducible $U_{q}(\mathfrak{g})$-modules, and so
the action of $\widetilde{R}$ on such tensor products is well-defined.
This allows us to define $R$-matrices for representations of
$U_{q}(\mathfrak{g})$ later in this paper.

\begin{subsection}{The universal $R$-matrix of $U_{h}(\mathfrak{g})$}

Khoroshkin and Tolstoy  \cite{kt} wrote down a universal $R$-matrix of $U_{h}(\mathfrak{g})$ 
using infinite sums of root vectors in $U_{h}(\mathfrak{g})$. 
Yamane also wrote down a universal $R$-matrix of $U_{h}(\mathfrak{g})$ \cite{y} 
but we use Khoroshkin and Tolstoy's work for ease.
The root vectors in \cite{kt} are defined in a different way to the way in which root vectors in universal
$R$-matrices of quantum algebras are defined. 
Khoroshkin and Tolstoy's procedure is general for quantum superalgebras 
and we write it down here for $U_{h}(\mathfrak{g})$.

For the universal $R$-matrix of $U_{h}(\mathfrak{g})$, 
we only define root vectors for the elements of the
{\emph{reduced root system}} $\phi$ of $\mathfrak{g}$.
The reduced root system $\phi$ is the set of all 
positive roots of $\mathfrak{g}$ except those roots $\alpha$ for which $\alpha/2$ is also a positive root:
the reduced root system of $\mathfrak{g}=osp(1|2n)$ 
is $\phi=\overline{\Phi}^{+}_{0} \cup \Phi^{+}_{1} = 
\{\epsilon_{i} \pm \epsilon_{j}, \epsilon_{k} | \ 1 \leq i < j \leq n, 1 \leq k \leq n \}$.

We introduce a total ordering of $\phi$ called a {\emph{normal ordering}}
which we denote by ${\cal{N}}(\phi)$, 
and then recursively define the root vectors of $U_{h}(\mathfrak{g})$ using ${\cal{N}}(\phi)$ 
and a map involving the $q$-bracket we introduce below. 
The way in which the universal $R$-matrix of $U_{h}(\mathfrak{g})$ is formally written down 
explicitly depends on ${\cal{N}}(\phi)$.
A difference between the root vectors in quantum algebras and 
the root vectors in $U_{q}(\mathfrak{g})$ is that
the latter are defined by a map that is not necessarily an algebra automorphism.

A normal ordering of the reduced root system of $\mathfrak{g}$ is defined 
as follows \cite[Def. 3.1]{kt}.
\begin{definition}
A normal ordering ${\cal{N}}(\phi)$ of $\phi =  \overline{\Phi}^{+}_{0} \cup \Phi^{+}_{1}$
is a total order $\prec$ of the elements of $\phi$ such that if $\alpha \prec \beta$ and 
$\alpha + \beta \in \phi$, then $\alpha \prec \alpha + \beta \prec \beta$.
\end{definition}
 
In general, there is more than one normal ordering of $\phi$ \cite{kt}.  
For example, the reduced root system of
$osp(1|4)$ is $\phi = \{ \epsilon_{1}, \epsilon_{2}, \epsilon_{1} \pm \epsilon_{2}\}$ 
and there are two different normal orderings of $\phi$:
$$\alpha_{1} \prec \alpha_{1} + \alpha_{2} \prec \alpha_{1} + 2\alpha_{2} \prec \alpha_{2},$$
$$\alpha_{2} \prec \alpha_{1} + 2\alpha_{2} \prec \alpha_{1} + \alpha_{2} \prec \alpha_{1},$$
writing the elements of $\phi$ as sums of the simple roots.

The universal $R$-matrix of $U_{h}(\mathfrak{g})$ \cite{kt}, 
adapted slightly to take account of the different co-multiplication used in this paper, 
is as follows. Let us write $q=e^{h} \in \mathbb{C}[[h]]$, and
$\displaystyle{\exp_{q}{(x)} = \sum_{k=0}^{\infty} \frac{x^{k}}{\left[k\right]_{q}!} }$. 

Construct the root vectors as follows.  
Firstly fix $E_{\alpha_{i}}=E_{i}$, $F_{\alpha_{i}}=F_{i}$ and $H_{\alpha_{i}}=H_{i}$ 
for each simple root $\alpha_{i}$.
Now recursively construct the non-simple root vectors.
Let $\alpha, \beta, \gamma \in \phi$ be roots such that $\gamma = \alpha + \beta$ and 
$\alpha \prec \beta$, and furthermore let no other roots 
$\alpha', \beta' \in \phi$ exist which  
satisfy $\alpha' + \beta' = \gamma$, 
 $\alpha \prec \alpha' \prec \beta$ and  $\alpha \prec \beta' \prec \beta$.  
Then, if all of the root vectors 
$E_{\alpha}, E_{\beta}, F_{\alpha}, F_{\beta} \in U_{h}(\mathfrak{g})$ 
have already been defined, fix
$$E_{\gamma} = \left[E_{\alpha},E_{\beta} \right]_{q}, \hspace{10mm} 
F_{\gamma} = \left[F_{\beta},F_{\alpha}\right]_{q^{-1}},$$
 where the $q$-bracket
 $[ \cdot, \cdot ]_{q}$
 is defined by
 $$[X_{\alpha}, X_{\beta}]_{q} = 
    X_{\alpha} X_{\beta} - (-1)^{[X_{\alpha}][X_{\beta}]} q^{(\alpha,\beta)} X_{\beta}X_{\alpha},$$
where we replace $X$ with $E$ or $F$ as appropriate.

Now for each $\gamma \in \phi$, fix
\begin{equation}
\label{eq:johnfaulkner(10)}
\widetilde{R}_{\gamma} = 
\exp_{q_{\gamma}}{\big((-1)^{[E_{\gamma}]} (a_{\gamma})^{-1}(q-q^{-1}) E_{\gamma} \otimes F_{\gamma}\big)}
\in U_{h}(\mathfrak{g}) \otimes U_{h}(\mathfrak{g}),
\end{equation}
where $q_{\gamma} = (-1)^{[E_{\gamma}]} q^{-(\gamma, \gamma)}$ and
$a_{\gamma} \in \mathbb{C}[[h]]$ is defined by 
$$E_{\gamma} F_{\gamma}-(-1)^{[E_{\gamma}]}F_{\gamma}E_{\gamma}=
\frac{a_{\gamma}\left(q^{H_{\gamma}}-q^{-H_{\gamma}}\right)}{q-q^{-1}}.$$
It is important to observe that $a_{\gamma}$ is a rational function of $q$.
Now we can write down the universal $R$-matrix of $U_{h}(\mathfrak{g})$  \cite[Thm. 8.1]{kt}.
\begin{theorem}
Define $\mathbf{H}_{i} = \sum_{j=i}^{n} H_{j} \in U_{h}(\mathfrak{g})$ for each $i=1, \ldots, n$.
The universal $R$-matrix of $U_{h}(\mathfrak{g})$ is
\begin{equation}
\label{eq:reneeboyle(1)}
R=\exp{\left(h \sum_{i=1}^{n} \mathbf{H}_{i} \otimes \mathbf{H}_{i} \right)}
\prod_{\gamma \in \phi} \widetilde{R}_{\gamma},
\end{equation}
where the product 
is ordered with respect to the same normal ordering ${\cal{N}}(\phi)$ 
that was used to define the root vectors in $U_{h}(\mathfrak{g})$ so that
 $\prod_{\gamma \in \phi} \widetilde{R}_{\gamma} = 
  \widetilde{R}_{\gamma_{1}} \widetilde{R}_{\gamma_{2}} \cdots \widetilde{R}_{\gamma_{k}}$ 
if the normal ordering ${\cal{N}}(\phi)$ is
 $\gamma_{1} \prec \gamma_{2} \prec \cdots \prec \gamma_{k}$.
 
\end{theorem}

\end{subsection}

\begin{subsection}{$R$-matrices for representations of $U_{q}(\mathfrak{g})$}
\label{subsect:Rmatricesfrorepresofquantumosp}

It is unknown whether Jimbo's quantum algebras admit universal $R$-matrices.
However, it is well-known that there are 
{\emph{$R$-matrices for representations}} of these quantum algebras.
Let $\pi_{\lambda}$ and $\pi_{\mu}$ be finite dimensional irreducible representations 
of the quantum algebra $A$, then there exists an invertible element 
${\cal{R}}_{\lambda \mu} \in \mathrm{End}_{\mathbb{C}}(V_{\lambda} \otimes V_{\mu})$ satisfying  
 \begin{equation}
\label{chap2:FyHendy(1)}
{\cal{R}}_{\lambda \mu} \cdot (\pi_{\lambda} \otimes \pi_{\mu}) \big(\Delta(x)\big) = 
        (\pi_{\lambda} \otimes \pi_{\mu})  \big(\Delta'(x)\big)  \cdot {\cal{R}}_{\lambda \mu}
	\hspace{10mm} \forall x \in A.
\end{equation}
We will show that
for each tensor product $V_{\lambda} \otimes V_{\mu}$
 of finite dimensional irreducible $U_{q}(\mathfrak{g})$-modules, there exists a map
${\cal{R}}_{\lambda \mu} \in \mathrm{End}_{\mathbb{C}}(V_{\lambda} \otimes V_{\mu})$ satisfying 
(\ref{chap2:FyHendy(1)}) for all $x \in U_{q}(\mathfrak{g})$.
We will do this following the method shown in \cite{cp,ks} for representations of quantum algebras.

We firstly define a completion $\overline{U}^{+}_{q}(\mathfrak{g})$ of $U_{q}(\mathfrak{g})$.
%following \cite[Subsec. 6.3.3]{ks}.
Let $U_{q}(\mathfrak{n}_{+})$ (resp. $U_{q}(\mathfrak{n}_{-})$) be the subalgebra  of
$U_{q}(\mathfrak{g})$ generated by $\{ e_{i} | \ 1 \leq i \leq n \}$
(resp. $\{ f_{i} | \ 1 \leq i \leq n \}$).  
We say that a non-zero element $x \in U_{q}(\mathfrak{g})$ has degree 
$\lambda = \sum_{i=1}^{n} m_{i}\alpha_{i}, \ m_{i} \in \mathbb{Z}$, 
if $K_{i}xK_{i}^{-1} = q^{(\lambda,\alpha_{i})}x$ for all $i=1, 2, \ldots, n$. 
We define $\overline{U}^{\pm}_{q}(\mathfrak{g})$ by
$$\overline{U}^{\pm}_{q}(\mathfrak{g}) = 
\prod_{\beta \in Q_{+}} U_{q}(\mathfrak{b}_{\pm})U_{q}^{\mp \beta}(\mathfrak{n}_{\mp}),$$
where $U_{q}^{\mp \beta}(\mathfrak{n}_{\mp})$ is defined by
$$U_{q}^{\pm \beta}(\mathfrak{n}_{\pm}) = 
\{ x \in U_{q}(\mathfrak{n}_{\pm}) |
 \ K_{i} x K_{i}^{-1} = q^{\pm (\alpha_{i},\beta)}x \}, \hspace{10mm} i=1, 2, \ldots, n,$$
 and $Q_{+}$ is defined by
 $Q_{+} = \{ \sum_{i=1}^{n} n_{i} \alpha_{i} | \ n_{i} \in \mathbb{Z}_{+} \}$.

The elements of $\overline{U}^{\pm}_{q}(\mathfrak{g})$ are sequences 
$x = (x_{\beta})_{\beta \in Q_{+}}$ where 
$x_{\beta} \in U_{q}(\mathfrak{b}_{\pm})U_{q}^{\mp \beta}(\mathfrak{n}_{\mp})$.
Let us write this sequence formally as an infinite sum $x = \sum_{\beta} x_{\beta}$.  
Then $U_{q}(\mathfrak{g})$ can be expressed as
$$U_{q}(\mathfrak{g}) = \bigoplus_{\beta \in Q_{+}} U_{q}(\mathfrak{b}_{+}) U_{q}^{-\beta}(\mathfrak{n}_{-}),$$
thus
$U_{q}(\mathfrak{g})$ can be considered as the subspace of
$\overline{U}^{+}_{q}(\mathfrak{g})$ formed by the sums $x=\sum_{\beta} x_{\beta}$ 
for which all but finitely many terms vanish.

The multiplication in $U_{q}(\mathfrak{g})$ extends to multiplications in 
$\overline{U}^{\pm}_{q}(\mathfrak{g})$, each of which is an algebra, as is 
$\overline{U}^{\pm}_{q}(\mathfrak{g}) \overline{\otimes} 
\cdots \overline{\otimes} \overline{U}^{\pm}_{q}(\mathfrak{g})$ ($m$ factors), 
and these algebras, respectively, contain 
$U_{q}(\mathfrak{g})$ and $U_{q}(\mathfrak{g}) \otimes \cdots \otimes U_{q}(\mathfrak{g})$ 
($m$ factors) as subalgebras.

We now construct an element in $\overline{U}^{+}_{q}(\mathfrak{g})$ corresponding to 
$\prod_{\gamma \in \phi} \widetilde{R}_{\gamma}$ in (\ref{eq:reneeboyle(1)}).
Given a normal ordering ${\cal{N}}(\phi)$ for a reduced root system $\phi$, 
we construct root vectors  $E_{\gamma}, F_{\gamma} \in U_{q}(\mathfrak{g})$ 
following the same procedure as in $U_{h}(\mathfrak{g})$ by setting
$E_{\alpha_{i}} = e_{i}$ and $F_{\alpha_{i}} = f_{i}$ and thinking of $q$ a complex number.
Each $\widetilde{R}_{\gamma}$ is then well-defined as an element of
$\overline{U}^{+}_{q}(\mathfrak{g}) \overline{\otimes} \overline{U}^{+}_{q}(\mathfrak{g})$,
and to simplify its expression we normalise the root vectors:
$$e_{\gamma} = E_{\gamma}, \hspace{10mm} f_{\gamma} = F_{\gamma}/a_{\gamma}.$$
(The expression for $f_{\gamma}$
is well-defined as $a_{\gamma} \neq 0$ \cite[Eqs. (8.3)--(8.4)]{kt}.)
Then, we have
$$\widetilde{R}_{\gamma} = \left\{ \begin{array}{ll}
\displaystyle{\sum_{k=0}^{\infty} 
\frac{\left(q-q^{-1}\right)^{k}\left(e_{\gamma} \otimes f_{\gamma}\right)^{k}}
{\left[k\right]_{q^{-2}}!},} & \mbox{if } [e_{\gamma}]=0,  \\
\displaystyle{\sum_{k=0}^{\infty} 
\frac{\left(q^{-1}-q\right)^{k}\left(e_{\gamma} \otimes f_{\gamma}\right)^{k}}
{\left[k\right]_{-q^{-1}}!}, } & \mbox{if } [e_{\gamma}]=1. 
\end{array} \right.$$
Define $\widetilde{R} \in 
\overline{U}^{+}_{q}(\mathfrak{g}) \overline{\otimes} \overline{U}^{+}_{q}(\mathfrak{g})$ by
$\widetilde{R} = \prod_{\gamma \in \phi} \widetilde{R}_{\gamma}$ where the product is ordered 
using the same normal order ${\cal{N}}(\phi)$ that was used to define the root vectors in
$U_{q}(\mathfrak{g})$; ie 
$\widetilde{R} = \widetilde{R}_{\gamma_{1}} \widetilde{R}_{\gamma_{2}} \cdots \widetilde{R}_{\gamma_{k}}$
where ${\cal{N}}(\phi)$ is 
$\gamma_{1} \prec \gamma_{2} \prec \cdots \prec \gamma_{k}$.
Clearly, $\widetilde{R}$ is invertible as 
$\prod_{\gamma \in \phi} \widetilde{R}_{\gamma} \in U_{h}(\mathfrak{g}) \otimes U_{h}(\mathfrak{g})$
is invertible and $q$ is not a root of unity.

\begin{lemma}
\label{lemlem:twenty}
Define an automorphism $\Psi$ of $U_{q}(\mathfrak{g}) \otimes U_{q}(\mathfrak{g})$ by
$$
\begin{array}{lll}
\Psi(K_{i}^{\pm 1} \otimes 1) = K_{i}^{\pm 1} \otimes 1, & & 
\Psi(1 \otimes K_{i}^{\pm 1}) = 1 \otimes K_{i}^{\pm 1}, \\
\Psi(e_{i} \otimes 1) = e_{i} \otimes K_{i}^{-1}, & &
\Psi(1 \otimes e_{i}) = K_{i}^{-1} \otimes e_{i}, \\
\Psi(f_{i} \otimes 1) = f_{i} \otimes K_{i}, & &
\Psi(1 \otimes f_{i}) = K_{i} \otimes f_{i}.
\end{array}
$$
The automorphism $\Psi$ satisfies the following relations:
\begin{itemize}
\item[(i)] $\widetilde{R} \Delta(x) = \Psi \big(\Delta'(x) \big) \cdot \widetilde{R}, 
\hspace{10mm} \mbox{for all } x \in U_{q}(\mathfrak{g})$,
\item[(ii)] $(\Delta \otimes \mathrm{id}) \widetilde{R}=\Psi_{23} (\widetilde{R}_{13}) \cdot \widetilde{R}_{23}$,
\item[(iii)] $(\mathrm{id} \otimes \Delta) \widetilde{R}=\Psi_{12} (\widetilde{R}_{13}) \cdot \widetilde{R}_{12}$,
\end{itemize} 
where $\Psi_{12} = \Psi \otimes \mathrm{id}$ and $\Psi_{23} = \mathrm{id} \otimes \Psi$.
\end{lemma}
\begin{proof}
We prove (i) for each generator of $U_{q}(\mathfrak{g})$.  
Firstly, we wish to prove the following equations:
\begin{eqnarray}
\label{eq:oneoneoneoneone}
\widetilde{R} (e_{i} \otimes K_{i} + 1 \otimes e_{i}) & = & 
(e_{i} \otimes K_{i}^{-1} + 1 \otimes e_{i})\widetilde{R}, \\
\label{eq:twotwo}
\widetilde{R} (f_{i} \otimes 1 + K_{i}^{-1} \otimes f_{i}) & = &
(f_{i} \otimes 1 + K_{i} \otimes f_{i})\widetilde{R}, \\
\label{eq:threethreethree}
\widetilde{R} (K_{i}^{\pm 1} \otimes K_{i}^{\pm 1}) & = & 
(K_{i}^{\pm 1} \otimes K_{i}^{\pm 1}) \widetilde{R}.
\end{eqnarray}
Eq. (\ref{eq:threethreethree}) is true by inspection and
Eqs. (\ref{eq:oneoneoneoneone})--(\ref{eq:twotwo}) follow from the corresponding results in 
$U_{h}(\mathfrak{g})$ \cite[Prop. 6.2]{kt}.
The proof of (i) then follows from the definition of $\Psi$ and 
the proofs of (ii) and (iii) follow similarly from \cite{kt}.
\end{proof}

We now examine the usual approach used to create $R$-matrices for representations
of a  quantum algebra $A$.
For each tensor product $W_{1} \otimes W_{2}$
of finite dimensional integrable $A$-modules, an invertible element 
${\cal{E}}_{W_{1}, W_{2}} \in \mathrm{End}_{\mathbb{C}} (W_{1} \otimes W_{2})$ is constructed
implementing the automorphism $\Psi$, in the sense that ${\cal{E}}_{W_{1}, W_{2}}$ satisfies
$${\cal{E}}_{W_{1}, W_{2}}^{-1} \cdot (\pi_{W_{1}} \otimes \pi_{W_{2}})(x) 
\cdot {\cal{E}}_{W_{1}, W_{2}} = 
(\pi_{W_{1}} \otimes \pi_{W_{2}})\Psi(x), 
\hspace{5mm} \forall x \in A \otimes  A.$$
This ${\cal{E}}_{W_{1}, W_{2}}$ is completely fixed by defining its action to be
$${\cal{E}}_{W_{1}, W_{2}} (w_{\lambda} \otimes w_{\mu}) = 
   q^{(\lambda,\mu)} (w_{\lambda} \otimes w_{\mu}),$$
on all weight vectors $ w_{\lambda} \in W_{1}$, $w_{\mu} \in W_{2}$ 
 with weights $\lambda$ and $\mu$, respectively \cite[Prop. 10.1.19]{cp}.

We could use the same method to construct $R$-matrices for representations of
$U_{q}(\mathfrak{g})$ but we have found a more useful approach.
Above, one needs to know
the weight space decompositions of each of  $W_{1}$ and $W_{2}$
before defining ${\cal{E}}_{W_{1}, W_{2}}$. With this knowledge, we can do something more 
interesting: instead of defining an element of $\mathrm{End}_{\mathbb{C}}(W_{1} \otimes W_{2})$,
we can define an element $E_{W_{2}} \in U_{q}(\mathfrak{g})$ with the property that
$(\pi_{W_{1}} \otimes \pi_{W_{2}}) E_{W_{2}}={\cal{E}}_{W_{1}, W_{2}}$.
A reason for doing this 
is that this allows us to define, for each finite dimensional tensorial irreducible representation 
$V_{\lambda}$ of $U_{q}(\mathfrak{g})$,
a class of even invertible elements of 
$\overline{U}^{+}_{q}(\mathfrak{g})$, each of which acts  
as a non-zero scalar multiple of the identity on each vector in $V_{\lambda}$.
The specific non-zero scalar is $q^{-(\lambda, \lambda + 2\rho)}$. 
We define $E_{W_{2}}$ for each tensor product of 
finite dimensional irreducible $U_{q}(\mathfrak{g})$-modules
following a related idea in \cite{z1}.

For each $i=1, \ldots, n$, set $J_{i} = K_{i}K_{i+1} \cdots K_{n}$.
The action of $J_{i}$ on a weight vector 
$w_{\xi}$ with weight $\xi = \sum_{j=1}^{n} \xi_{j} \epsilon_{j} \in H^{*}$ 
 of a $U_{q}(\mathfrak{g})$-module is
$$J_{i} w_{\xi} = q^{\xi_{i}} w_{\xi}.$$

Consider the weight space decomposition of a 
finite dimensional irreducible $U_{q}(\mathfrak{g})$-module $V_{\mu}$ 
with integral dominant highest weight $\mu$.  
The weight of the weight vector $w_{\xi} \in V_{\mu}$ is 
$\xi = \sum_{i=1}^{n} \xi_{i} \epsilon_{i} \in \bigoplus_{i=1}^{n} \mathbb{Z} \epsilon_{i}$.  
Now define
\begin{equation}
\label{eq:chatswoodtrain}
E_{\mu} = \prod_{a=1}^{n} \sum^{s}_{b=p} (J_{a})^{b} \otimes P_{a}[b], \hspace{10mm} 
P_{a}[b] = \prod_{\stackrel{c = p}{c \neq b}}^{s} \frac{J_{a}-q^{c}}{q^{b}-q^{c}}, \hspace{5mm} c \in \mathbb{Z},
\end{equation}
where $p$ and $s$ are any integers satisfying $p \leq s$ and the further condition that
\begin{itemize}
\item[(i)]
$J_{i} w_{\xi} = q^{\xi_{i}} w_{\xi}$ 
for some $\xi_{i}$ satisfying $p \leq \xi_{i} \leq s$, 
for each weight vector $w_{\xi} \in V_{\mu}$.
\end{itemize}
Once we have any such $p$ and $s$ satisfying these conditions, 
we can use any $p'$ and $s'$ satisfying
$p' \leq p$ and $s' \geq s$ in (\ref{eq:chatswoodtrain}) instead of $p$ and $s$, respectively.

The element $E_{\mu}$ is well-defined and invertible in 
$U_{q}(\mathfrak{g}) \otimes U_{q}(\mathfrak{g})$, and
for all weight vectors $v_{\lambda'} \in V_{\lambda}$ and $v_{\mu'} \in V_{\mu}$ we have
\begin{equation}
\label{chapter2:opera(1)}
E_{\mu} (v_{\lambda'} \otimes v_{\mu'}) = q^{(\lambda',\mu')} \ (v_{\lambda'} \otimes v_{\mu'}),
\end{equation}
where the weights of $v_{\lambda'}$ and $v_{\mu'}$ are $\lambda'$ and $\mu'$, respectively.
The element $E_{\mu}$ is not a {\emph{universal}} element in that it does not 
satisfy (\ref{chapter2:opera(1)}) for all representations of $U_{q}(\mathfrak{g})$;
it would be useful if one could construct such a universal element.

Using this, we obtain $R$-matrices for tensor products of finite dimensional irreducible
$U_{q}(\mathfrak{g})$-modules in the following sense.
Let $V_{\lambda}$ and $V_{\mu}$ be irreducible $U_{q}(\mathfrak{g})$-modules with  
integral dominant highest weights $\lambda$ and $\mu$, respectively.  
Then we have the following important theorem.
\begin{theorem}
\label{th:RRRmatrix}
Define $R_{\mu} = E_{\mu} \widetilde{R} \in 
\overline{U}_{q}^{+}(\mathfrak{g}) \overline{\otimes} \overline{U}_{q}^{+}(\mathfrak{g})$
and  ${\cal{R}}_{\lambda \mu} = (\pi_{\lambda} \otimes \pi_{\mu})R_{\mu}$, then
\begin{equation}
\label{eq:sachanadra(1)}
{\cal{R}}_{\lambda \mu} \cdot (\pi_{\lambda} \otimes \pi_{\mu}) \big(\Delta(x)\big) 
    = (\pi_{\lambda} \otimes \pi_{\mu})\big(\Delta'(x)\big) \cdot {\cal{R}}_{\lambda \mu}, 
    \hspace{10mm} \forall x \in U_{q}(\mathfrak{g}).
\end{equation}
\end{theorem}
\begin{proof}
This is similar to the proof of the corresponding theorem in quantum algebras \cite[Prop. 10.1.19]{cp}.
A direct calculation readily shows that
$$(\pi_{\lambda} \otimes \pi_{\mu}) \Psi \big(\Delta(x)\big)  = 
(\pi_{\lambda} \otimes \pi_{\mu})  \big( E_{\mu}^{-1} \cdot \Delta(x) \cdot E_{\mu} \big),
\hspace{10mm} \forall x \in U_{q}(\mathfrak{g}),$$
 then by using  
$\widetilde{R}\Delta(x) = \Psi \big(\Delta'(x)\big) \cdot \widetilde{R}$ from 
Lemma \ref{lemlem:twenty}, we have
\begin{equation}
\label{chap2:equationannie(alpha)}
(\pi_{\lambda} \otimes \pi_{\mu}) \big(R_{\mu} \Delta(x) \big) = 
(\pi_{\lambda} \otimes \pi_{\mu}) \big( \Delta'(x) R_{\mu} \big),
\end{equation}
which is precisely Eq. (\ref{eq:sachanadra(1)}).
\end{proof}
\noindent

We now determine some useful results involving $R_{\mu}$.
\begin{proposition}
The element $R_{\mu}$ has the following properties:
\begin{equation}
\label{chapt2:annieeq(1)}
(\epsilon \otimes \mathrm{id}) R_{\mu} = (\mathrm{id} \otimes \epsilon) R_{\mu} = 1,
\end{equation}
\begin{equation}
\label{chapt2:annieeq(2)}
(\pi_{\lambda} \otimes \pi_{\mu})\big( (S \otimes \mathrm{id}) R_{\mu}\big) = 
(\pi_{\lambda} \otimes \pi_{\mu})R_{\mu}^{-1}, \ 
  (\pi_{\lambda} \otimes \pi_{\mu})\big((\mathrm{id} \otimes S) R_{\mu}^{-1}\big) = 
   (\pi_{\lambda} \otimes \pi_{\mu})R_{\mu},
  \end{equation}
  \begin{equation}
  \label{chapt2:annieeq(3)}
  (\pi_{\lambda} \otimes \pi_{\mu})\big((S \otimes S) R_{\mu}\big) = 
  (\pi_{\lambda} \otimes \pi_{\mu})R_{\mu}.
\end{equation}
\end{proposition}
\begin{proof}
Eq. (\ref{chapt2:annieeq(1)}) is proved by inspection.  
The proofs of (\ref{chapt2:annieeq(2)})--(\ref{chapt2:annieeq(3)}) are straightforward
and almost identical to the proofs of the corresponding equations 
in $\mathbb{Z}_{2}$-graded quasitriangular Hopf algebras. 
\end{proof}

Let $v_{\lambda}$ and $v_{\nu}$ be weight vectors of $U_{q}(\mathfrak{g})$-modules 
with weights
$\lambda$, $\nu \in \bigoplus_{i=1}^{n} \mathbb{Z} \epsilon_{i}$, respectively
and let $v_{\mu'} \in V_{\mu}$ be a weight vector with weight $\mu'$, then
it can be easily shown that
$$\big[(\Delta \otimes 1)E_{\mu}\big] (v_{\lambda} \otimes v_{\nu} \otimes v_{\mu'}) = 
q^{(\mu',\lambda + \nu)}(v_{\lambda} \otimes v_{\nu} \otimes v_{\mu'}),$$
$$\big[(1 \otimes \Delta)E_{\mu}\big] (v_{\lambda} \otimes v_{\nu} \otimes v_{\mu'}) =
q^{(\lambda, \nu + \mu')}(v_{\lambda} \otimes v_{\nu} \otimes v_{\mu'}),$$
where in  $(1 \otimes \Delta)E_{\mu}$ we change the limits $p$ and $s$ if necessary.

We now consider analogues in $U_{q}(\mathfrak{g})$ of the equations
$(\Delta \otimes 1) R = R_{13} R_{23}$ and $(1 \otimes \Delta)R = R_{13} R_{12}$
of a $\mathbb{Z}_{2}$-graded quasitriangular Hopf algebra.  
By definition, we have
$(\Delta \otimes 1) R_{\mu} = 
\big[ (\Delta \otimes 1) E_{\mu} \big] \cdot \Psi_{23}(\widetilde{R}_{13}) \cdot \widetilde{R}_{23}$.
Let $V_{\lambda}$ and $V_{\nu}$ be finite dimensional irreducible $U_{q}(\mathfrak{g})$-modules, 
then from the properties of $E_{\mu}$ we have
\begin{eqnarray}
(\pi_{\lambda} \otimes \pi_{\nu} \otimes \pi_{\mu}) \big[ (\Delta \otimes 1) R_{\mu} \big] 
& = &  (\pi_{\lambda} \otimes \pi_{\nu} \otimes \pi_{\mu})
        \left[  (\Delta \otimes 1) E_{\mu}  \cdot 
	(E_{\mu}^{23})^{-1} \widetilde{R}_{13} E_{\mu}^{23}\widetilde{R}_{23} \right] \nonumber \\
& = &  (\pi_{\lambda} \otimes \pi_{\nu} \otimes \pi_{\mu})
        \left[ E_{\mu}^{13} \widetilde{R}_{13} E_{\mu}^{23} \widetilde{R}_{23} \right], 
	\label{chapter2:R-matrix(99)}
\end{eqnarray}
where writing $E_{\mu} = \sum_{t} \alpha_{t} \otimes \beta_{t}$ we have
$E_{\mu}^{13} = \sum_{t} \alpha_{t} \otimes \mathrm{id} \otimes \beta_{t}$ and
$E_{\mu}^{23}= \sum_{t}  \mathrm{id} \otimes \alpha_{t} \otimes \beta_{t}$.
Note that (\ref{chapter2:R-matrix(99)}) uses  the result
$$(\pi_{\lambda} \otimes \pi_{\nu} \otimes \pi_{\mu})
        \left[ (\Delta \otimes 1) E_{\mu}  \cdot (E_{\mu}^{23})^{-1} \right] 
	= (\pi_{\lambda} \otimes \pi_{\nu} \otimes \pi_{\mu}) E_{\mu}^{13},$$
rather than an equality in $U_{q}(\mathfrak{g})^{\otimes 3}$.
Similarly, we have
$$(1 \otimes \Delta) R_{\mu} = 
\left[(1 \otimes \Delta) E_{\mu} \right] \cdot \Psi_{12} 
(\widetilde{R}_{13}) \cdot \widetilde{R}_{12}, \hspace{5mm} \mbox{and}$$ 
\begin{equation}
\label{chapter2:R-matrix(100)}
(\pi_{\lambda} \otimes \pi_{\nu} \otimes \pi_{\mu}) \big[ (1 \otimes \Delta) R_{\mu} \big]
= (\pi_{\lambda} \otimes \pi_{\nu} \otimes \pi_{\mu})
      \left[ E^{13}_{\mu} \widetilde{R}_{13} E_{\nu}^{12} \widetilde{R}_{12} \right]. 
      \end{equation}
Together with Theorem \ref{th:RRRmatrix}, this shows that  
$R_{\mu}$ satisfies the defining relations 
(\ref{eq:definingrelationsR1})--(\ref{eq:definingrelationsR3}) 
of the universal $R$-matrix of a $\mathbb{Z}_{2}$-graded quasitriangular Hopf algebra 
in {\emph{each triple tensor product}} of finite dimensional
irreducible $U_{q}(\mathfrak{g})$ representations, 
if we carefully choose the limits in the definition of $E_{\mu}$, which can always be done.  
Furthermore,  
Eqs. (\ref{chap2:equationannie(alpha)}) and 
(\ref{chapter2:R-matrix(99)})--(\ref{chapter2:R-matrix(100)}) imply that
\begin{equation}
\label{chapter2:R-matrix(101)}
(\pi_{\lambda} \otimes \pi_{\lambda} \otimes \pi_{\lambda}) R_{12} R_{13} R_{23} 
 = (\pi_{\lambda} \otimes \pi_{\lambda} \otimes \pi_{\lambda})   R_{23}R_{13}R_{12} ,
\end{equation}
where we have fixed $R = R_{\lambda}$.

For later use, note the following easily proved results.
Define an automorphism 
$\Psi_{m}: U_{q}(\mathfrak{g})^{\otimes m} \rightarrow U_{q}(\mathfrak{g})^{\otimes m}$ 
generalising the automorphism 
$\Psi: U_{q}(\mathfrak{g}) \otimes U_{q}(\mathfrak{g}) \rightarrow
U_{q}(\mathfrak{g}) \otimes U_{q}(\mathfrak{g})$ in Lemma \ref{lemlem:twenty} by
$$\begin{array}{lcl}
\Psi_{m} (1^{\otimes j} \otimes K_{i}^{\pm 1} \otimes 1^{\otimes (m-j-1)}) & = & 
   1^{\otimes j} \otimes K_{i}^{\pm 1} \otimes 1^{\otimes (m-j-1)}, \\
\Psi_{m} (1^{\otimes j} \otimes e_{i} \otimes 1^{\otimes (m-j-1)}) & = & 
   (K_{i}^{-1})^{\otimes j} \otimes e_{i} \otimes (K_{i}^{-1})^{\otimes (m-j-1)}, \\
\Psi_{m} (1^{\otimes j} \otimes f_{i} \otimes 1^{\otimes (m-j-1)}) & = & 
   (K_{i})^{\otimes j}  \otimes f_{i} \otimes (K_{i})^{\otimes (m-j-1)}, 
\end{array}$$
for each $1 \leq i \leq n$ and all $0 \leq j \leq m-1$. Then
\begin{eqnarray}
\big(\Delta \otimes \mathrm{id}^{\otimes t} \big) \Psi_{2,3,\ldots, t+1} (\widetilde{R}_{1(t+1)})
& = & \Psi_{2,3,\ldots, t+2}(\widetilde{R}_{1(t+2)}) \cdot 
   \Psi_{3,4,\ldots, t+2}(\widetilde{R}_{2(t+2)}) \label{chap2:eqmorganspurlock(1)} \\
\big( \mathrm{id}^{\otimes t} \otimes \Delta \big) \Psi_{1,2,\ldots, t} (\widetilde{R}_{1(t+1)})
& = & \Psi_{1,2,\ldots, t+1} (\widetilde{R}_{1(t+2)}) \cdot 
   \Psi_{1,2,\ldots, t} (\widetilde{R}_{1(t+1)}) \label{chap2:eqmorganspurlock(agent92)}
\end{eqnarray}
where  $\Psi_{k, \ldots, m} = \mathrm{id}^{\otimes (k-1)} \otimes \Psi_{m-k+1}$, 
$k \geq 2$, in (\ref{chap2:eqmorganspurlock(1)}) and
$\Psi_{1, \ldots, m} = \Psi_{m} \otimes \mathrm{id}$ in (\ref{chap2:eqmorganspurlock(agent92)}).
Then it may be easily shown that
\begin{eqnarray}
\left( \pi^{\otimes t} \otimes \pi \right)  \big( \Delta^{(t-1)} \otimes \mathrm{id} \big)R 
 & = &
 \left( \pi^{\otimes t} \otimes \pi \right) R_{1(t+1)} R_{2(t+1)} \cdots R_{t(t+1)}, 
 \label{starequationlikeearth} \\
\left( \pi \otimes \pi^{\otimes t} \right)  \big( \mathrm{id} \otimes \Delta^{(t-1)} \big)R 
 & = & 
  \left( \pi \otimes \pi^{\otimes t} \right) R_{1(t+1)} R_{1t} \cdots R_{12}, \nonumber  
\end{eqnarray}
where we fix $R = R_{V}$.

\end{subsection}

\end{section}

\begin{section}{Two useful elements of $\overline{U}^{+}_{q}(\mathfrak{g})$ }
%\markright{\text{Two useful elements of $\overline{U}^{+}_{q}(\mathfrak{g})$ }}
\label{sec:'central'element}

Define a set of elements 
$\{u_{\lambda} \in \overline{U}^{+}_{q}(\mathfrak{g}) | \ \lambda \in {\cal{P}}^{+} \}$ by
$$u_{\lambda} = \sum_{t} S(b_{\lambda_{t}}) a_{\lambda_{t}} (-1)^{[a_{\lambda_{t}}]},$$
where  $R_{\lambda} = \sum_{t} a_{\lambda_{t}} \otimes b_{\lambda_{t}}$. 
The following lemma was proved in \cite{blumen05}.
\begin{lemma}
\label{lemming:overtheygo}
The element $u_{\lambda}$ has the following properties:
\begin{itemize}
\item[(i)] $\epsilon(u_{\lambda})=1$,
\item[(ii)] $\pi_{\lambda}\big( S^{2}(x) u_{\lambda} \big) = \pi_{\lambda}\left( u_{\lambda}x \right)$, 
$\forall x \in U_{q}(\mathfrak{g})$.
\item[(iii)]  $\pi_{\lambda}\big(u_{\lambda} \widetilde{u}_{\lambda}\big) = \pi_{\lambda}(1)= 
\pi_{\lambda}\big(\widetilde{u}_{\lambda}u_{\lambda}\big)$, where 
$\widetilde{u}_{\lambda}$ is defined by
$\widetilde{u}_{\lambda} = \sum_{s} S^{-1}(d_{\lambda_{s}})c_{\lambda_{s}} (-1)^{[c_{\lambda_{s}}]}$,
where
$R_{\lambda}^{-1} = \sum_{s} c_{\lambda_{s}} \otimes  d_{\lambda_{s}}$,
\item[(iv)] $(\pi_{\lambda} \otimes \pi_{\lambda})\big(\Delta(u_{\lambda})\big) = 
(\pi_{\lambda} \otimes \pi_{\lambda})
\left[(u_{\lambda} \otimes u_{\lambda})\big(R^{T}_{\lambda}R_{\lambda}\big)^{-1}\right]$.
\end{itemize} 
\end{lemma}

Now define a set of elements 
$\{v_{\lambda} \in \overline{U}^{+}_{q}(\mathfrak{g}) | \ \lambda \in {\cal{P}}^{+} \}$ by
\begin{equation}
\label{eq:defofthevelements}
v_{\lambda} = u_{\lambda} K_{2 \rho}^{-1}.
\end{equation}
\begin{lemma}
\label{lem:galliano}
The element $v_{\lambda}$ has the following properties:
$$\epsilon(v_{\lambda})=1, \hspace{10mm} \pi_{\lambda}(v_{\lambda}x)=\pi_{\lambda}(x v_{\lambda}), 
\hspace{5mm} \forall x \in U_{q}(\mathfrak{g}),$$
\begin{equation}
\label{eq:laurasia}
(\pi_{\lambda} \otimes \pi_{\lambda})
\Delta(v_{\lambda})=(\pi_{\lambda} \otimes \pi_{\lambda})
\left[(v_{\lambda} \otimes v_{\lambda})\big(R_{\lambda}^{T}R_{\lambda}\big)^{-1}\right].
\end{equation}
\end{lemma}
\begin{proof}
The proofs of $\epsilon(v_{\lambda})=1$ and (\ref{eq:laurasia}) follow from the definition of 
$v_{\lambda}$.  To prove the remaining relation, note that
$S^{2}(e_{i}) = K_{i}e_{i}K_{i}^{-1} = K_{2\rho}e_{i}K_{2\rho}^{-1}$,
$S^{2}(f_{i}) = K_{i}f_{i}K_{i}^{-1} = K_{2\rho}f_{i}K_{2\rho}^{-1}$,
and $S^{2}(K^{\pm 1}_{i}) = K_{2\rho}K_{i}^{\pm 1}K_{2\rho}^{-1}$.
As $S^{2}$ is a homomorphism we have
$S^{2}(x) = K_{2\rho}x K_{2\rho}^{-1}$ for all $x \in U_{q}(\mathfrak{g})$ and then 
\begin{eqnarray*}
\pi_{\lambda}\big(v_{\lambda} x v_{\lambda}^{-1}\big) 
& = & \pi_{\lambda} \big( u_{\lambda} K_{2\rho}^{-1} x K_{2\rho} u_{\lambda}^{-1}\big) \\
& = & \pi_{\lambda} \big(u_{\lambda} S^{-2}(x) u^{-1}_{\lambda} \big)  
 =  \pi_{\lambda} \big( S^{2}(S^{-2}(x)) \big) 
= \pi_{\lambda}(x),
\end{eqnarray*}
completing the proof.
\end{proof}

\begin{lemma}
\label{lem:markhoopkins(a)}
The element $v_{\lambda}$ acts on each vector in the irreducible $U_{q}(\mathfrak{g})$-module
$V_{\lambda}$ with highest weight $\lambda \in {\cal{P}}^{+}$
as the multiplication by the scalar $q^{-(\lambda + 2\rho, \lambda)}$.
\end{lemma}  
\begin{proof}
Note that $v_{\lambda}$ is even and that 
$\pi_{\lambda}(v_{\lambda}) \in \mathrm{End}_{U_{q}(\mathfrak{g})}(V_{\lambda})$.
Write $R_{\lambda} = E_{\lambda} \widetilde{R}$ where
$\widetilde{R}=\sum_{t=0}^{\infty} a_{t} \otimes b_{t}$, 
$a_{t} \in U_{q}(\mathfrak{b}_{+})$, $b_{t} \in U_{q}(\mathfrak{b}_{-})$ 
and $a_{0}=b_{0}=1$, then
$$\pi_{\lambda}(v_{\lambda}) = 
\pi_{\lambda}\left(\sum_{t=0}^{\infty} S(b_{t}) E a_{t} K_{2\rho}^{-1} (-1)^{[a_{t}]}\right),$$
where $E$ is an even element of $U_{q}(\mathfrak{g})$ satisfying
$E w_{\xi} = q^{-(\xi,\xi)} w_{\xi}$
for each weight vector $w_{\xi} \in V_{\lambda}$ of weight 
$\xi \in \bigoplus_{i=1}^{n} \mathbb{Z} \epsilon_{i}$.
Let $w_{\lambda}$ be a non-zero highest weight vector of $V_{\lambda}$, 
then $a_{t} w_{\lambda} = 0$ for all $t > 0$,  yielding 
$$v_{\lambda} \cdot w_{\lambda} = E K_{2\rho}^{-1} w_{\lambda} =
q^{-(\lambda+2\rho,\lambda)}w_{\lambda}.$$

As it is true that $V_{\lambda}$ is irreducible,
$\pi_{\lambda}(v_{\lambda}) \in \mathrm{End}_{U_{q}(\mathfrak{g})}(V_{\lambda})$
and $\pi_{\lambda}(v_{\lambda})$ is a homogeneous map of degree zero,
it follows that
$v_{\lambda}$ acts on each weight vector $w \in V_{\lambda}$ 
as the multiplication by the claimed
scalar from Schur's lemma.
\end{proof} 

We may denote $q^{-(\lambda + 2\rho, \lambda)}$ by $\chi_{\lambda}(v_{\lambda})$.
Note that 
$v_{\mu}$ may act on each weight vector $w \in V_{\lambda}$ as the multiplication by the 
salar $q^{-(\lambda + 2\rho, \lambda)}$ even if  $\mu \neq \lambda$, and in this case 
we also write $\chi_{\lambda}(v_{\mu})$ to denote $q^{-(\lambda + 2\rho, \lambda)}$.

Following \cite{z3}, we define the quantum superdimension of the 
finite dimensional $U_{q}(\mathfrak{g})$-module $W$ to be 
$$sdim_{q}(W) = \mathrm{str}_{q} (\mathrm{id}_{W}).$$
The following lemma, first stated in \cite{zsuper}, was proved in \cite{blumen05} drawing on \cite{k0,k1}.
\begin{lemma}
Let $V_{\lambda}$ be a finite dimensional irreducible 
$U_{q}(\mathfrak{g})$-module with integral dominant highest weight $\lambda$. 
The quantum superdimension of $V_{\lambda}$ is
\begin{equation}
\label{eq:qqqsuperdim}
sdim_{q}(V_{\lambda})
= (-1)^{[\lambda]}q^{-(\lambda,2\rho)}
\prod_{\alpha \in \overline{\Phi}^{+}_{0}}\left(
\frac{q^{2(\lambda+\rho,\alpha)}-1}{q^{2(\rho,\alpha)}-1}\right)
\prod_{\beta \in \Phi_{1}^{+}}\left(
\frac{q^{2(\lambda+\rho,\beta)}+1}{q^{2(\rho,\beta)}+1}\right), 
\end{equation}
where $[\lambda]$ is the grading of the highest weight vector of $V_{\lambda}$.
\end{lemma}

It is easy to calculate that the quantum superdimension of the fundamental ($2n+1$)-dimensional
$U_{q}(osp(1|2n))$-module $V$ is
\begin{equation}
\label{eq:johnfaulkner(12)}
sdim_{q}(V) = 1 - \frac{q^{2n}-q^{-2n}}{q-q^{-1}},
\end{equation}
where we recall that the grading of the highest weight vector of $V$ is odd.

\end{section}

\begin{section}{Spectral decomposition of $\check{\cal{R}}_{V,V}$}
\label{sec:spectam(a)}
%\markright{\text{Spectral decomposition of $\check{\cal{R}}_{V,V}$}}

Let $V_{\lambda}$ and $V_{\mu}$ be finite dimensional irreducible
$U_{q}(\mathfrak{g})$-modules with integral dominant highest weights $\lambda$ and $\mu$, respectively.
Let $R_{\mu}$ be as in Theorem \ref{th:RRRmatrix} and define 
$\check{\cal{R}}_{V_{\lambda}, V_{\mu}} 
\in \mathrm{Hom}_{\mathbb{C}}(V_{\lambda} \otimes V_{\mu}, V_{\mu} \otimes V_{\lambda})$ by
\begin{equation}
\label{eq:bigjimmyboy(2)}
\check{\cal{R}}_{V_{\lambda},V_{\mu}}(v_{\lambda} \otimes v_{\mu}) = 
P \circ \big(R_{\mu}(v_{\lambda} \otimes v_{\mu})\big),
\end{equation}
where $v_{\lambda} \in V_{\lambda}$ and $v_{\mu} \in V_{\mu}$ are weight vectors
and $P$ is the graded permutation operator.
A standard argument proves the following lemma:
\begin{lemma}
Let $V_{\lambda}$ be an irreducible $U_{q}(\mathfrak{g})$-module  
with integral dominant highest weight $\lambda$, then
$\check{\cal{R}}_{V_{\lambda},V_{\lambda}} \in \mathrm{End}_{U_{q}(\mathfrak{g})}(V_{\lambda} \otimes V_{\lambda})$.
\end{lemma}
For $n=1$, $V \otimes V$ decomposes into a direct sum of irreducible 
$U_{q}(\mathfrak{g})$-modules:
\begin{equation}
\label{eq:johnfaulkner(14)}
V \otimes V = V_{2 \epsilon_{1}} \oplus V_{\epsilon_{1}} \oplus V_{0},
\end{equation}
and for $n \geq 2$, we have
\begin{equation}
\label{eq:johnfaulkner(13)}
V \otimes V = V_{2 \epsilon_{1}} \oplus V_{\epsilon_{1} + \epsilon_{2}} \oplus V_{0}.
\end{equation}
\begin{lemma}
\label{lem:ohgremlinsintheworks}
Let $n \geq 2$ and let
$\{ P[\mu] \in \mathrm{End}_{U_{q}(\mathfrak{g})}(V \otimes V) | 
 \ \mu= 2\epsilon_{1}, \epsilon_{1} + \epsilon_{2}, 0 \}$ be a set of even 
$U_{q}(\mathfrak{g})$-linear maps: 
$P[\mu]: V \otimes V \rightarrow V \otimes V$,
where the image of $P[\mu]$ is $V_{\mu}$ and the maps
satisfy  $\big(P[\mu]\big)^{2} = P[\mu]$ and $P[\mu] P[\nu] = \delta_{\mu \nu} P[\mu]$.
Then there is a spectral decomposition of $\check{\cal{R}}_{V,V}$:
$$\check{\cal{R}}_{V,V} = -q P[2\epsilon_{1}] + q^{-1}P[\epsilon_{1} + \epsilon_{2}] + q^{-2n} P[0].$$
\end{lemma}
\begin{proof}
As $\check{\cal{R}}_{V,V} \in \mathrm{End}_{U_{q}(\mathfrak{g})}(V \otimes V)$, 
we can write
$$\check{\cal{R}}_{V,V} = 
\beta_{2 \epsilon_{1}} P[2 \epsilon_{1}] + \beta_{\epsilon_{1}+\epsilon_{2}}P[\epsilon_{1}+\epsilon_{2}] +
\beta_{0}P[0],$$
for some set of constants 
$\{ \beta_{\mu} \in \mathbb{C} | \ \mu = 2\epsilon_{1}, \epsilon_{1}+\epsilon_{2}, 0 \}$,
where $\beta_{\mu}$ is the scalar action of $\check{\cal{R}}_{V,V}$ on 
the irreducible $U_{q}(\mathfrak{g})$-submodule $V_{\mu} \subset V \otimes V$.  
We explicitly calculate each $\beta_{\mu}$ using $R_{V}$.

Let $\{v_{i} | \ -n \leq i \leq n\}$ be the basis of 
weight vectors of $V$ given in Lemma \ref{lem:fundamentaldimensional}. 
The highest weight vector of $V_{2\epsilon_{1}}$ is
$w_{2\epsilon_{1}} = v_{1} \otimes v_{1}$,
the highest weight vector of $V_{\epsilon_{1} + \epsilon_{2}}$ is
$w_{\epsilon_{1} + \epsilon_{2}} = v_{1} \otimes v_{2} - q^{-1} v_{2} \otimes v_{1}$
and the highest weight vector of the trivial module $V_{0} \subset V \otimes V$ is 
$w_{0}=\sum_{i=-n}^{n} c_{i} v_{i} \otimes v_{-i}$,
where $\{c_{i} \in \mathbb{C} | \ -n \leq i \leq n\}$ 
is a set of non-zero constants inductively defined by
$$
\begin{array}{rclrcl}
c_{n} & = & -c_{0}, & c_{-n} & = & q^{-1}c_{0}, \\
c_{n-1} & = & -qc_{n}, & c_{-(n-1)} & = & -q^{-1}c_{-n}, \\
c_{i} & = & -q c_{i+1}, & c_{-i} & = & -q^{-1} c_{-(i+1)},
\end{array} $$
where $i=1, 2, \ldots, n-2$ and we fix $c_{0} \neq 0$.

To study the action of $\check{\cal{R}}_{V,V}$ on the highest weight vectors 
$w_{2 \epsilon_{1}}, w_{\epsilon_{1} + \epsilon_{2}}$ and $w_{0}$, 
we make some observations about $(\pi \otimes \pi) \widetilde{R}$. 
From the weight space decomposition of $V$, we have
\begin{eqnarray*}
 \pi ( f_{\epsilon_{i}})^{3} = \pi ( e_{\epsilon_{i}})^{3}=0, & \hspace{5mm} & \mbox{for all }  
i=1, \ldots, n, \\
\pi ( f_{\gamma})^{2} = \pi ( e_{\gamma})^{2}=0, & \hspace{5mm} & 
\mbox{for all }  \gamma \in \phi \mbox{ where } \gamma \neq \epsilon_{i}, 
\end{eqnarray*}
and thus
\begin{equation}
\label{eq:johnfaulkner(15)}
(\pi \otimes \pi) \widetilde{R} = 
(\pi \otimes \pi) \prod_{\gamma \in \phi}  \widetilde{R}^{V}_{\gamma},
\end{equation}
where 
\begin{eqnarray*}
\widetilde{R}^{V}_{\gamma} & = & \left\{ \begin{array}{lll}
\displaystyle{\sum_{k=0}^{2} \frac{(q^{-1}-q)^{k} (e_{\gamma} \otimes f_{\gamma})^{k}}{ [k]^{-q^{-1}}!}},
& & \mbox{if } \gamma = \epsilon_{i}, \\
\displaystyle{\sum_{k=0}^{1} \frac{ (q-q^{-1})^{k} (e_{\gamma} \otimes f_{\gamma})^{k}}{ [k]^{q^{-2}}!  
} }, & &  \mbox{if } \gamma \neq \epsilon_{i}.
\end{array} \right. \\
& = & \left\{ \begin{array}{lll}
\displaystyle{ 1 \otimes 1 + (q^{-1}-q)(e_{\gamma} \otimes f_{\gamma}) + 
                       \frac{(q^{-1}-q)^{2}(e_{\gamma} \otimes f_{\gamma})^{2}}{(1-q^{-1})} },
 & & \mbox{if } \gamma = \epsilon_{i}, \\
\displaystyle{ 1 \otimes 1 + (q-q^{-1}) (e_{\gamma} \otimes f_{\gamma})   }, 
  & & \mbox{if }  \gamma \neq \epsilon_{i},
\end{array} \right.
\end{eqnarray*}
and where the product in (\ref{eq:johnfaulkner(15)})
is ordered using the same normal ordering ${\cal{N}}(\phi)$ used to construct
the root vectors, ie we fix
$\prod_{\gamma \in \phi} \widetilde{R}^{V}_{\gamma} =  
\widetilde{R}^{V}_{\gamma_{1}} \widetilde{R}^{V}_{\gamma_{2}} \cdots \widetilde{R}^{V}_{\gamma_{k}}$
where ${\cal{N}}(\phi) = \gamma_{1} \prec \gamma_{2} \prec \cdots \prec \gamma_{k}$.
The expression for $(\pi \otimes \pi) \widetilde{R}$ in (\ref{eq:johnfaulkner(15)})
readily assists the use of $(\pi \otimes \pi) R_{V}$ and $\check{\cal{R}}_{V,V}$ in calculations.

Using this, 
we have
$\check{\cal{R}}_{V,V} (w_{2 \epsilon_{1}}) = -q w_{2 \epsilon_{1}}$ and
$\check{\cal{R}}_{V,V} (w_{\epsilon_{1} + \epsilon_{2}}) = q^{-1} w_{\epsilon_{1} + \epsilon_{2}}$.
Calculating $\beta_{0}$ is more difficult: note that
$$
\check{\cal{R}}_{V,V}\left( c_{-1} v_{-1} \otimes v_{1} + \sum_{\stackrel{i=-n}{i \neq -1}}^{n} c_{i} v_{i} \otimes
v_{-i}\right) = -q^{-1} c_{-1} v_{1} \otimes v_{-1} + \sum_{\stackrel{j=-n}{j \neq -1}}^{n} c_{j}'
v_{-j} \otimes v_{j},
$$
for some set of non-zero constants 
$\left\{c_{j}' \in \mathbb{C} | \ -n \leq j \leq n, \ j \neq -1 \right\}$.
Recall that $\check{\cal{R}}_{V,V} (w_{0}) = \beta_{0} w_{0}$, 
so we obtain $\beta_{0}$ by comparing $-q^{-1}c_{-1}$ and $c_{1}$.
Now $c_{-1} = (-1)^{n-1} q^{-n} c_{0}$ and $c_{1} = (-1)^{n} q^{n-1} c_{0}$, thus $\beta_{0} = q^{-2n}$.
\end{proof}

\begin{lemma}
\label{lem:ohgremlinsintheworks(2)}
Let $n=1$ and let
$\{ P[\mu] \in \mathrm{End}_{U_{q}(\mathfrak{g})}(V \otimes V) | \ \mu= 2\epsilon_{1}, \epsilon_{1}, 0 \}$ 
be a set of even 
$U_{q}(\mathfrak{g})$-linear maps:
$P[\mu]: V \otimes V \rightarrow V \otimes V$,
where the image of $P[\mu]$ is $V_{\mu}$ and the maps
satisfy  $\big(P[\mu]\big)^{2} = P[\mu]$ and $P[\mu] P[\nu] = \delta_{\mu \nu} P[\mu]$.
Then there is a spectral decomposition of $\check{\cal{R}}_{V,V}$:
$$\check{\cal{R}}_{V,V} = -q P[2\epsilon_{1}] + q^{-1}P[\epsilon_{1}] + q^{-2} P[0].$$
\end{lemma}
\begin{proof}
The proof is almost identical to the proof of Lemma \ref{lem:ohgremlinsintheworks} 
except for the following minor difference.
The highest weight vector of $V_{\epsilon_{1}}$ 
in the decomposition of $V \otimes V$ in (\ref{eq:johnfaulkner(14)}) is 
$w_{\epsilon_{1}} = v_{1} \otimes v_{0} + q^{-1} v_{0} \otimes v_{1}$.
To complete the proof we note that $\check{\cal{R}}_{V,V} (w_{\epsilon_{1}}) = q^{-1} w_{\epsilon_{1}}$.
\end{proof}

\begin{corollary}
\label{cor:stonemebudgie}
For each $n = 1, 2, \ldots$, $\check{\cal{R}}_{V,V}$ satisfies 
\begin{equation}
\label{eq:RRRRequat}
(\check{\cal{R}}_{V,V}+q)(\check{\cal{R}}_{V,V}-q^{-1})(\check{\cal{R}}_{V,V}-q^{-2n})=0.
\end{equation}
\end{corollary}

\end{section}

\begin{section}{A representation of the Birman-Wenzl-Murakami algebra $BW_{t}(-q^{2n},q)$}
\label{eq:theXfactor(a)}
%\markright{\text{A representation of the Birman-Wenzl-Murakami algebra}}

In this section we recall the 
Birman-Wenzl-Murakami algebra $BW_{t}(r,q)$ from \cite{bw,w2} 
and define a representation of $BW_{t}(-q^{2n},q)$ in
a subalgebra ${\cal{C}}_{t}$ of the centraliser algebra 
$\mathrm{End}_{U_{q}(\mathfrak{g})}(V^{\otimes t})$.  

\begin{definition}
Define ${\cal{C}}_{t}$ to be the subalgebra of 
$\mathrm{End}_{U_{q}(\mathfrak{g})}(V^{\otimes t})$
generated by the elements
$$\left\{\check{R}_{i}^{\pm 1} \in \mathrm{End}_{U_{q}(\mathfrak{g})}(V^{\otimes t}) |
 \ 1 \leq i \leq t-1\right\}, \hspace{5mm} \mbox{where}$$ 
\begin{equation}
\label{eq:tom3}
\check{R}_{i} = \mathrm{id}^{\otimes (i-1)} \otimes \check{\cal{R}}_{V,V}\otimes \mathrm{id}^{\otimes (t-(i+1))}
\in \mathrm{End}_{U_{q}(\mathfrak{g})}(V^{\otimes t}).
\end{equation}
\end{definition}
\noindent

Let us investigate ${\cal{C}}_{t}$.
Let $\{v_{i} | \ -n \leq i \leq n\}$ be the basis of weight vectors
of $V$ given in Lemma \ref{lem:fundamentaldimensional} and let
$\{v^{*}_{i} | \ -n \leq i \leq n\}$ be a basis of $V^{*}$ such that
$\langle v^{*}_{i}, v_{j}\rangle = \delta_{ij}$ and $[v^{*}_{i}]=[v_{i}]$; then
$$av_{i} = \sum_{j} \langle v_{j}^{*}, av_{i} \rangle v_{j}, \hspace{10mm} 
av_{i}^{*} = \sum_{j} \langle av_{i}^{*}, v_{j} \rangle v_{j}^{*}, 
\hspace{10mm} \forall a \in U_{q}(\mathfrak{g}).$$
Define $\check{e} \in \mathrm{End}_{\mathbb{C}}(V \otimes V^{*})$ by
$$\check{e}(x \otimes y^{*}) = 
(-1)^{[y^{*}][x]}\langle y^{*}, v^{-1}u \ x\rangle \sum_{i=-n}^{n} v_{i} \otimes v_{i}^{*},$$ 
where $v$ and $u$ are the elements
$v_{\epsilon_{1}}, u_{\epsilon_{1}} \in \overline{U}_{q}^{+}(\mathfrak{g})$ respectively.

\begin{lemma}
\label{lem:marrickville}
The map $\check{e}$ satisfies
\begin{itemize}
\item[(i)]  $(\check{e})^{2}=sdim_{q}(V)\check{e}$,
\item[(ii)] $a\check{e}=\epsilon(a) \check{e},$  $\hspace{5mm} \forall a \in U_{q}(\mathfrak{g})$,
\item[(iii)] $\check{e} a=\epsilon(a) \check{e}$, $\hspace{5mm} \forall a \in U_{q}(\mathfrak{g})$,
\item[(iv)] $\check{e}_{2} \check{R}_{1} \check{e}_{2} = q^{2n} \check{e}_{2},$
where 
$$\check{e}_{2} = \mathrm{id}_{V} \otimes \check{e}: V \otimes V \otimes V^{*}
\rightarrow V \otimes V \otimes V^{*},$$
$$\check{R}_{1} = \check{\cal{R}}_{V,V} \otimes \mathrm{id}_{V^{*}}: V \otimes V \otimes V^{*}
\rightarrow V \otimes V \otimes V^{*}.$$
\end{itemize}
\end{lemma}
\begin{proof}
\begin{itemize}
\item[(i)]
\begin{eqnarray*}
(\check{e})^{2}(x \otimes y^{*}) & = & 
(-1)^{[y^{*}][x]}\langle y^{*}, v^{-1}u x \rangle
\sum_{i} (-1)^{[v_{i}]}\langle v_{i}^{*}, v^{-1}u \ v_{i} \rangle \sum_{j} v_{j} \otimes v_{j}^{*} \\
& = & sdim_{q}(V) (-1)^{[y^{*}][x]}\langle y^{*}, v^{-1}u \ x \rangle
\sum_{j} v_{j} \otimes v_{j}^{*} 
  = sdim_{q}(V) \check{e} (x \otimes y^{*}).
\end{eqnarray*}
\item[(ii)] By definition, $a \check{e} = a_{V \otimes V^{*}} \circ \check{e}$; we calculate that
\begin{eqnarray*}
a \check{e}(x \otimes y^{*}) 
& = & (-1)^{[y^{*}][x]} \langle y^{*}, v^{-1}u x \rangle
\sum_{(a), i, j, k} \langle v_{j}^{*}, a_{(1)} v_{i} \rangle
\langle v_{i}^{*}, S(a_{(2)}) v_{k} \rangle v_{j} \otimes v_{k}^{*} \\
& = & (-1)^{[y^{*}][x]} \langle y^{*}, v^{-1}u x \rangle
\sum_{(a),j,k} \langle v_{j}^{*}, a_{(1)}S(a_{(2)}) v_{k} \rangle 
v_{j} \otimes v_{k}^{*} \\
& = & \epsilon(a) (-1)^{[y^{*}][x]} \langle y^{*}, v^{-1}u x \rangle 
\sum_{k} v_{k} \otimes v_{k}^{*} 
  =  \epsilon(a) \check{e}(x \otimes y^{*}).
\end{eqnarray*}
\end{itemize}
\noindent
Similar calculations prove (iii) and (iv) (see \cite{lr} for the corresponding calculations in 
ungraded quasitriangular Hopf algebras).
\end{proof}

Define $\hat{e} \in \mathrm{End}_{\mathbb{C}}(V^{*} \otimes V)$ by
$$\hat{e}(x^{*} \otimes y) = 
\langle x^{*}, y \rangle \sum_{i=-n}^{n} (-1)^{[v_{i}]} v_{i}^{*} \otimes vu^{-1} \ v_{i},$$
where $v$ and $u$ are fixed to be $v_{\epsilon_{1}}$ and $u_{\epsilon_{1}}$, respectively.
\begin{lemma}
\label{lemming:sausage2}
The map $\hat{e}$ satisfies
\begin{itemize}
\item[(i)]  $(\hat{e})^{2}=sdim_{q}(V)\hat{e}$,
\item[(ii)] $a \hat{e}=\epsilon(a) \hat{e}$, $\hspace{5mm} \forall a \in U_{q}(\mathfrak{g})$,
\item[(iii)] $\hat{e} a= \epsilon(a) \hat{e}$, $\hspace{5mm} \forall a \in U_{q}(\mathfrak{g})$,
\item[(iv)] $\hat{e}_{2} \check{R}_{1}^{-1} \hat{e}_{2} = q^{-2n} \hat{e}_{2}$
where 
$\hat{e}_{2} = \mathrm{id}_{V^{*}} \otimes \hat{e}: V^{*} \otimes V^{*} \otimes V
\rightarrow V^{*} \otimes V^{*} \otimes V$.
\end{itemize}
\end{lemma}
\begin{proof}
The proofs of (i)--(iv) are similar to the proofs of (i)--(iv) in Lemma \ref{lem:marrickville}.
\end{proof}

\begin{remark}  The maps $\check{e}$ and $\hat{e}$
are $U_{q}(\mathfrak{g})$-invariant maps onto one-dimensional $U_{q}(\mathfrak{g})$-submodules 
in $V \otimes V^{*}$ and $V^{*} \otimes V$, respectively.
\end{remark}

Recall that $V \otimes V$ has the decomposition into irreducible $U_{q}(\mathfrak{g})$-modules  
given in (\ref{eq:johnfaulkner(14)})--(\ref{eq:johnfaulkner(13)})
and that there exists an even $U_{q}(\mathfrak{g})$-invariant 
map $P[0]: V \otimes V \rightarrow V \otimes V$ the image of which is $V_{0} \subset V \otimes V$,
 defined in Lemmas \ref{lem:ohgremlinsintheworks}--\ref{lem:ohgremlinsintheworks(2)}.
Recall that $V^{*} \cong V$ and define 
$$E = (\mathrm{id} \otimes T^{-1}) \circ \check{e} \circ (\mathrm{id} \otimes T)
= (T^{-1} \otimes \mathrm{id}) \circ \hat{e} \circ (T \otimes \mathrm{id}) = sdim_{q}(V)P[0],$$
where $T$ is the isomorphism $T: V \rightarrow V^{*}$ given in (\ref{eq:johnfaulkner(20)}).
Furthermore, define the elements
$$E_{i} =  \mathrm{id}^{\otimes (i-1)} \otimes E \otimes \mathrm{id}^{\otimes (t-(i+1))}
\in \mathrm{End}_{U_{q}(\mathfrak{g})}(V^{\otimes t}), \hspace{5mm}  i=1, 2, \ldots, t-1.$$
\begin{lemma}
\label{lem:10001}
The elements $\check{R}_{i}$, 
$E_{i} \in \mathrm{End}_{U_{q}(\mathfrak{g})}(V^{\otimes t})$ satisfy the 
relations
\begin{itemize}
\item[(i)] $\check{R}_{i} \check{R}_{i+1} \check{R}_{i} = 
\check{R}_{i+1} \check{R}_{i} \check{R}_{i+1}$, $ \ 1 \leq i \leq t-2$,
\item[(ii)] $\check{R}_{i} \check{R}_{j} = \check{R}_{j} \check{R}_{i}$, $ \ |i-j|>1$,
\item[(iii)] $(\check{R}_{i}+q)(\check{R}_{i}-q^{-1})(\check{R}_{i}-q^{-2n})=0$, 
$ \ 1 \leq i \leq t-1$,
\item[(iv)] $-\check{R}_{i} + \check{R}_{i}^{-1} = (q-q^{-1})(1-E_{i})$,
\item[(v)] $E_{i} \check{R}_{i-1}^{\pm 1} E_{i} = q^{\pm 2n} E_{i}$, 
\item[(vi)] $E_{i} \check{R}_{i}^{\pm 1} = \check{R}_{i}^{\pm 1}E_{i}
=q^{\mp 2n} E_{i}$, $ \ 1 \leq i \leq t-1$.
\end{itemize}
\end{lemma}
\begin{proof}
The proofs of (i) and (ii) are standard arguments.
The proof of (iii) follows from Corollary \ref{cor:stonemebudgie}.
The proof of (v) follows from Lemmas \ref{lem:marrickville}--\ref{lemming:sausage2}.
The proofs of (iv) and (vi) follow from the definition of 
$E_{i}$, Eq. (\ref{eq:johnfaulkner(12)})
and the fact that $\check{R}_{1}$ acts on $V_{0} \subset V \otimes V$ as 
$\check{R}_{1}w = q^{-2n}w$ for all $w \in V_{0}$.
\end{proof}

We now give the definition of the Birman-Wenzl-Murakami algebra $BW_{t}(r, q)$ from \cite{rw}.
Let $r$ and $q$ be non-zero complex constants and let $t \geq 2$ be an integer.  
The Birman-Wenzl-Murakami algebra
$BW_{t}(r, q)$ is the algebra over $\mathbb{C}$ generated by the invertible elements
$\{ g_{i} | \ 1 \leq i \leq t-1 \}$ subject to the relations
\newline
$$
\begin{array}{ll}
g_{i} g_{j} = g_{j} g_{i}, 				 & |i-j| >1, 	\\
g_{i} g_{i+1} g_{i} = g_{i+1} g_{i} g_{i+1}, 	 & 1 \leq i \leq t-2, 	\\
e_{i} g_{i} = r^{-1} e_{i},                      & 1 \leq i \leq t-1, 	\\
e_{i}g_{i-1}^{\pm 1}e_{i} = r^{\pm 1} e_{i}, 	 & 1 \leq i \leq t-1,
\end{array} $$
\newline
where $e_{i}$ is defined by
$$(q-q^{-1})(1-e_{i}) = g_{i}-g_{i}^{-1}, \hspace{10mm}  1 \leq i \leq t-1.$$
It can be shown that each $g_{i}$ also satisfies
$$(g_{i} - r^{-1})(g_{i} + q^{-1})(g_{i} - q) = 0.$$
From Lemma \ref{lem:10001} we have the following:
\begin{lemma}
\label{eq:peacockinblackbeansauce}
Let $q \in \mathbb{C}$ be non-zero and not a root of unity.  
The algebra homomorphism 
$\Upsilon: BW_{t}(-q^{2n}, q) \rightarrow {\cal{C}}_{t}$ defined by
$$\Upsilon: g_{i} \mapsto -\check{R}_{i}$$
yields a representation of $BW_{t}(-q^{2n}, q)$ in ${\cal{C}}_{t}$.
\end{lemma}

\end{section}

\begin{section}{Bratteli diagrams and path algebras}
\label{eq:caseydonovan(a)}
%\markright{\text{Bratteli diagrams and path algebras}}

\begin{subsection}{Bratteli diagrams}

To proceed further with the study of the Birman-Wenzl-Murakami algebra
we consider the notions of {\emph{Bratelli diagrams}} and 
{\emph{Path algebras for Bratelli diagrams}}, both of which we take 
from \cite{lr}.  (The reader is also referred to \cite[Chap. 2]{ghj}).

A Bratteli diagram is an undirected graph that encodes information about
a sequence $\mathbb{C} \cong A_{0} \subset A_{1} \subset A_{2} \subset \cdots$ of 
inclusions of finite dimensional semisimple algebras \cite{rw}.  
The properties of a Bratteli diagram, graph-theoretically, are that:
\begin{itemize}
\item[(i)]  its vertices are elements of certain sets $\widetilde{A}_{t}$, $t \in \mathbb{Z}_{+}$, and
\item[(ii)] if $n(a,b) \in \mathbb{Z}_{+}$ denotes the number of edges between the vertices $a$ and $b$, 
then $n(a,b)=0$ for any vertices $a \in \widetilde{A}_{i}$ and $b \in \widetilde{A}_{j}$ where 
$|i-j| \neq 1$.
\end{itemize}
Assume that $\widetilde{A}_{0}$ consists of a unique vertex that we denote by $\emptyset$.  
We call the elements of $\widetilde{A}_{i}$ {\emph{shapes}} and say that $\widetilde{A}_{i}$
is the set of shapes on the $i^{th}$ level of the Bratteli diagram.  
If $\lambda \in \widetilde{A}_{i}$ is
connected to $\mu \in \widetilde{A}_{i+1}$, we write $\lambda \leq \mu$.  

A {\emph{multiplicity free Bratteli diagram}} is a Bratteli diagram in which any two
vertices are connected by at most one edge.  
All Bratteli diagrams considered in this paper are multiplicity free. 

Let $A$ be a Bratteli diagram and let $\lambda \in \widetilde{A}_{i}$ and 
$\mu \in \widetilde{A}_{j}$ for some $0 \leq i < j$.  
We define a {\emph{path from $\lambda$ to $\mu$}} to be a sequence of shapes
$$P = (s_{i}, s_{i+1}, \ldots, s_{j}),$$
where $\lambda = s_{i} \leq s_{i+1} \leq \cdots \leq s_{j-1} \leq s_{j} = \mu$ 
and $s_{k}$ is a shape on the $k^{th}$ level of $A$ for each $k$.

Given a path $T = (\lambda, \ldots, \mu)$ from $\lambda$ to $\mu$ and a path
$S = (\mu, \ldots, \nu)$ from $\mu$ to $\nu$, we define the concatenation 
of $T$ and $S$ to be the path from $\lambda$ to $\nu$ defined by
$$T \circ S = (\lambda, \ldots, \mu, \ldots, \nu).$$

We define a tableau $T$ of shape $\lambda$ to be a path 
from $\emptyset \in \widetilde{A}_{0}$ to $\lambda$ and we write
$shp(T) = \lambda$. We say that the length of $T$ is $t$
if there are $t+1$ shapes in the tableau.

\end{subsection}

\begin{subsection}{Path algebras related to Bratteli diagrams}

We now define the concept of a {\emph{Path algebra}} for a Bratteli diagram $A$.
For each $t \in \mathbb{Z}_{+}$, let ${\cal{T}}^{t}$  
be the set of tableaux of length $t$ in $A$
and let $\Omega^{t} \subset {\cal{T}}^{t} \times {\cal{T}}^{t}$ 
be the set of pairs $(S,T)$ of tableaux where
$shp(S)=shp(T)$, that is both $S$ and $T$ end in the same shape.  
Let us further define an algebra $A_{t}$ over $\mathbb{C}$ generated by  
$\{E_{ST} | \ (S, T) \in \Omega^{t}\}$, 
where the algebra multiplication is defined by 
\begin{equation}
\label{eq:multiplicationmatrixunits}
E_{ST}E_{PQ} = \delta_{TP} E_{SQ}.
\end{equation}
Any set of elements of an associative algebra satisfying (\ref{eq:multiplicationmatrixunits}) 
are called matrix units; matrix units also figure later in this paper.
Note that $A_{0} \cong \mathbb{C}$.  
Each element $a \in A_{t}$ can be written in the form
$$a = \sum_{(S,T) \in \Omega^{t}} a_{ST} E_{ST}, \hspace{10mm} a_{ST} \in \mathbb{C}.$$  
We refer to the collection of algebras $A_{t}$, $t \in \mathbb{Z}_{+}$, as the 
{\emph{tower of path algebras corresponding to the Bratteli diagram $A$}}.

Each of the algebras $A_{t}$ is isomorphic to a direct sum of matrix algebras.
The irreducible representations of $A_{t}$ are indexed by the elements of $\widetilde{A}_{t}$,
which is the set of shapes on the $t^{th}$ level of $A$.
Let ${\cal{T}}^{\lambda}$ denote the set of tableaux of shape $\lambda$, then 
the cardinality $d_{\lambda}$ of
${\cal{T}}^{\lambda} \cap {\cal{T}}^{t}$ is equal to the dimension of the irreducible $A_{t}$-module
indexed by $\lambda \in \widetilde{A}_{t}$.
We record this in the formula
$$A_{t} \cong \bigoplus_{\lambda \in \widetilde{A}_{t}} M_{d_{\lambda}} (\mathbb{C}),$$
where $M_{d} (\mathbb{C})$ denotes the algebra of $d \times d$ matrices with complex entries.

We now define some useful sets.
Let ${\cal{T}}^{\mu}_{\lambda}$ be the set of paths in $A$ from the shape $\lambda$ to the shape $\mu$
and let ${\cal{T}}^{t}_{r}$ be the set of paths 
%in $A$ from any shape in $\widetilde{A}_{r}$ 
%to any shape in $\widetilde{A}_{t}$, that is the set of paths 
starting on the $r^{th}$ level of $A$ and going down to the $t^{th}$ level. 
Furthermore, let ${\cal{T}}^{t}_{\lambda}$ be the set of paths in $A$ from the shape 
$\lambda$ to any shape on the $t^{th}$ level of $A$.

We also define 
$\Omega^{\mu}_{\lambda} \subset {\cal{T}}^{\mu}_{\lambda} \times {\cal{T}}^{\mu}_{\lambda}$ 
to be the set of pairs $(S,T)$ of paths $S, T \in {\cal{T}}^{\mu}_{\lambda}$ and 
$\Omega^{t}_{r} \subset {\cal{T}}^{t}_{r} \times {\cal{T}}^{t}_{r}$ 
to be the set of pairs $(S, T)$ of paths where in both situations we have $shp(S)=shp(T)$.

We define the inclusion of path algebras
$A_{r} \subseteq A_{t}$ for $0 \leq r < t$ as follows: for each
pair $(P, Q) \in \Omega^{r}$ we fix $E_{PQ} \in A_{t}$ by
$$E_{PQ} = \sum_{T \in {\cal{T}}^{t}_{\lambda} \cap {\cal{T}}^{t}_{r} } E_{P \circ T, Q \circ T}, \hspace{5mm}
\mbox{where } \lambda = shp(P) = shp(Q).$$
In particular, we have $A_{s} \subseteq A_{s+1}$ for each $s \in \mathbb{Z}_{+}$.

Let $\lambda \in \widetilde{A}_{t}$ and let ${\cal{V}}_{\lambda}$ be an irreducible representation of
$A_{t}$ indexed by $\lambda$.  The restriction of ${\cal{V}}_{\lambda}$ to 
the subalgebra $A_{t-1} \subseteq A_{t}$ decomposes into irreducible representations of $A_{t-1}$
according to
$${\cal{V}}_{\lambda} \downarrow^{A_{t}}_{A_{t-1}} \cong \bigoplus_{\mu \in \lambda^{-}}
{\cal{V}}_{\mu}, \hspace{10mm} \mbox{where } \lambda^{-} = \{ \nu \in \widetilde{A}_{t-1} | \ \nu \leq \lambda \}.$$
This decomposition is multiplicity free as the Bratteli diagram $A$ is multiplicity free. 

For each $r \in \mathbb{Z}_{+}$ satisfying $r < t$, 
the {\emph{centraliser of $A_{r}$ contained in $A_{t}$}} is defined to be 
$${\cal{L}}(A_{r} \subseteq A_{t}) = \{ a \in A_{t} | \ ab=ba, \ \forall b \in A_{r} \}.$$

Let now $(S,T)$ be a pair of paths each starting on the $r^{th}$ level of $A$ at the shape 
$\lambda$ and ending on the $t^{th}$ level of $A$ at the shape $\mu$.  
For each such pair we define $E_{ST} \in A_{t}$ by
$$E_{ST} = \sum_{P \in {\cal{T}}^{\lambda} \cap {\cal{T}}^{r}} E_{ P \circ S, P \circ T},$$
which we can think of as the sum of all pairs of paths 'ending' in $(S,T)$.
We then have the
following lemma, stated in \cite[Prop. $(1.4)$]{lr} and proved in \cite[Sect. 2.3]{ghj}.
\begin{lemma}
\label{lemonelem:one}
A basis of ${\cal{L}}(A_{r} \subseteq A_{t})$ is given by the elements 
$$\left\{ E_{ST} | \ (S, T) \in \Omega_{\lambda}^{\mu} \cap \Omega^{t}_{r}, 
 \ \lambda \in \widetilde{A}_{r}, \ \mu \in \widetilde{A}_{t}\right\}.$$
\end{lemma}

\end{subsection}

\begin{subsection}{Centraliser algebras}
\label{subsec:brattelidiagramsandcentraliser}

Let $U$ be a $\mathbb{Z}_{2}$-graded Hopf algebra over $\mathbb{C}$.
Let $V$ be a finite dimensional $U$-module with the property that
$V^{\otimes t}$ is completely reducible for each $t \in \mathbb{Z}_{+}$.
We now define the concepts of a {\emph{Bratteli diagram for tensor powers of $V$}} and
the {\emph{Bratteli diagram for $V^{\otimes t}$}}; the purpose of this subsection is to show that the centraliser 
${\cal{L}}_{t}$ of $U$ in $\mathrm{End}_{\mathbb{C}}(V^{\otimes t})$ defined by
${\cal{L}}_{t} = \mathrm{End}_{U}(V^{\otimes t})$
is isomorphic to the path algebra $A_{t}$ of the Bratteli diagram for $V^{\otimes t}$.

In this subsection we regard all modules as being graded.
By convention $V^{\otimes 0} \cong \mathbb{C}$ and thus ${\cal{L}}_{0} = \mathbb{C}$.
If $V$ is an irreducible $U$-module then ${\cal{L}}_{1} \cong \mathbb{C}$ by Schur's lemma.
For all $0 \leq r < t$ we define the inclusion ${\cal{L}}_{r} \subseteq {\cal{L}}_{t}$
by $a \mapsto a \otimes \mathrm{id}^{\otimes (t-r)}$ for all $a \in {\cal{L}}_{r}$.
Now ${\cal{L}}_{t}$ acts on $V^{\otimes t}$ in the obvious way.  
Since $U$ and ${\cal{L}}_{t}$ commute, $V^{\otimes t}$ has
a natural ${\cal{L}}_{t} \otimes U$-module structure.

Let $\{\Lambda_{\lambda} | \ \lambda \in I \}$ be the set of 
non-isomorphic finite dimensional irreducible $U$-modules.
Then by the double centraliser theorem there exists a finite subset 
$\widetilde{\cal{L}}_{t}$ of $I$ such that
$$V^{\otimes t} \cong \bigoplus_{ \lambda \in \widetilde{\cal{L}}_{t}}
{\cal{L}}^{\lambda} \otimes \Lambda_{\lambda},$$
where each ${\cal{L}}^{\lambda}$ is an irreducible ${\cal{L}}_{t}$-module such that 
${\cal{L}}^{\lambda} \not\cong {\cal{L}}^{\mu}$ if $\lambda \neq \mu$.

We now assume that $V$ is an irreducible $U$-module and 
continue to assume that all tensor powers of $V$ are completely reducible.
We will consider the Bratteli diagram for tensor powers of $V$.
Let $\lambda \in \widetilde{\cal{L}}_{t}$ for some $t$.  
Then we have the decomposition
\begin{equation}
\label{eq:100}
\Lambda_{\lambda} \otimes V = 
\bigoplus_{\mu \in \widetilde{\cal{L}}_{t+1}} \left(\Lambda_{\mu}\right)^{\oplus n_{\lambda}(\mu)},
\hspace{10mm} n_{\lambda}(\mu) \in \mathbb{Z}_{+},
\end{equation}
of $\Lambda_{\lambda} \otimes V$ into a direct sum of irreducible $U$-modules.  
The non-negative integer $n_{\lambda}(\mu)$ 
is the multiplicity of $\Lambda_{\mu}$ in the decomposition.
We say that the decomposition of $\Lambda_{\lambda} \otimes V$ is
{\emph{multiplicity free}} if $n_{\lambda}(\mu) \leq 1$ for all $\mu \in \widetilde{\cal{L}}_{t+1}$.   

The {\emph{branching rule for inclusion}} ${\cal{L}}_{t} \subseteq {\cal{L}}_{t+1}$ describes the decomposition
of the ${\cal{L}}_{t+1}$-module ${\cal{L}}^{\nu}$ into ${\cal{L}}_{t}$-modules
\begin{equation}
\label{eq:101}
{\cal{L}}^{\nu} = \bigoplus_{\lambda \in \widetilde{\cal{L}}_{t}}
\big( {\cal{L}}^{\lambda} \big)^{\oplus n_{\lambda}(\nu)},
\hspace{10mm} n_{\lambda}(\nu) \in \mathbb{Z}_{+}.
\end{equation}
Note that the positive integers $n_{\lambda}(\nu)$ appearing in (\ref{eq:100}) and (\ref{eq:101}) are the same \cite{lr}.

The {\emph{Bratteli diagram for tensor powers of $V$}} is defined as follows:
for each $t \in \mathbb{Z}_{+}$ fix the vertices on the $t^{th}$ level of the Bratteli diagram
to be the elements of 
$\widetilde{\cal{L}}_{t}$.  
Then a vertex $\lambda \in \widetilde{\cal{L}}_{t}$
is connected to a vertex $\mu \in \widetilde{\cal{L}}_{t+1}$ by $n_{\lambda}(\mu)$ edges.

For a fixed $t$, the {\emph{Bratteli diagram for $V^{\otimes t}$}} is an undirected graph with
vertices given by the elements of $\bigcup_{i=0}^{t} \widetilde{\cal{L}}_{i}$, 
and the edges are such that a vertex $\lambda \in \widetilde{\cal{L}}_{i}$ is connected
to a vertex $\mu \in \widetilde{\cal{L}}_{i+1}$ by $n_{\lambda}(\mu)$ edges 
for each $0 \leq i \leq t-1$.

Let $V$ be a finite dimensional irreducible $U$-module 
with the property that for every irreducible $U$-module $W$, 
the decomposition of the tensor product $W \otimes V$ is multiplicity free.  
In this case, we say that tensoring by $V$ is multiplicity free.  
We will show that the centraliser algebra
${\cal{L}}_{t} = 
\mathrm{End}_{U}(V^{\otimes t})$ is isomorphic to the path algebra $A_{t}$ associated with 
the Bratteli diagram for $V^{\otimes t}$.

We construct an algebra isomorphism $A_{t} \rightarrow {\cal{L}}_{t}$ inductively.  
Assume that there is an identification of ${\cal{L}}_{t}$ with 
the path algebra $A_{t}$ for some $t \geq 0$, so that
$$V^{\otimes t} = \bigoplus_{\lambda \in \widetilde{\cal{L}}_{t}}
\left( \bigoplus_{T \in {\cal{T}}^{\lambda} \cap {\cal{T}}^{t} } E_{TT} V^{\otimes t}\right)$$
is a decomposition of $V^{\otimes t}$ into irreducible $U$-modules $\Lambda_{\lambda}$
where the $U$-submodule $E_{TT} V^{\otimes t}$ is isomorphic to
$\Lambda_{\lambda}$ given $shp(T)=\lambda$.
The map $E_{TT}$ is a $U$-invariant map from $V^{\otimes t}$ onto a $U$-submodule
isomorphic to $\Lambda_{\lambda}$.

Let $T=( \emptyset, s_{1}, \ldots, \lambda) \in {\cal{T}}^{\lambda}$ be a tableau
of length $t$
and let $E_{TT} V^{\otimes t} \cong \Lambda_{\lambda}$ for some $\lambda \in \widetilde{\cal{L}}_{t}$.
As tensoring by $V$ is multiplicity free, the decomposition
\begin{equation}
\label{eq:equationa}
(E_{TT}V^{\otimes t}) \otimes V = 
\bigoplus_{\stackrel{\nu \in \widetilde{\cal{L}}_{t+1}}{\lambda \leq \nu}}V_{T \circ \nu},
\end{equation}
is multiplicity free and thus unique, where $T \circ \nu$ is the tableau
$$T \circ \nu = (\emptyset, s_{1}, \ldots, \lambda, \nu), \hspace{10mm} \lambda \leq \nu,$$
and $V_{T \circ \nu} \cong \Lambda_{\nu}$.

The next step is to identify
$E_{T \circ \nu, T \circ \nu}$ with the unique $U$-invariant projection operator mapping 
$\left(E_{TT}V^{\otimes t}\right) \otimes V$ onto $V_{T \circ \nu}$.  
This way we identify each element $E_{SS}$ of the path algebra
$A_{t+1}$, where $S \in {\cal{T}}^{t+1}$, with an element of ${\cal{L}}_{t+1}$.  
Thus we have the decomposition
$$V^{\otimes (t+1)} = \bigoplus_{\nu \in \widetilde{\cal{L}}_{t+1}}
\left( \bigoplus_{S \in {\cal{T}}^{\nu} \cap {\cal{T}}^{t+1}} E_{SS} V^{\otimes (t+1)} \right),$$
of $V^{\otimes (t+1)}$ into irreducible $U$-modules 
$E_{SS}V^{\otimes (t+1)} = V_{S} \cong \Lambda_{\nu}$, where
$\nu \in \widetilde{\cal{L}}_{t+1}$ and $S \in {\cal{T}}^{\nu} \cap {\cal{T}}^{t+1}$.

We now identify the other elements in the basis $\{ E_{PQ} \in A_{t+1} | \ (P,Q) \in \Omega^{t+1}\}$ 
with elements of ${\cal{L}}_{t+1}$.  
For each pair of paths $(P,Q) \in \Omega^{t+1}$ we choose non-zero elements
$$E_{PQ} \in E_{PP} {\cal{L}}_{t+1}E_{QQ}, 
\hspace{10mm} E_{QP} \in E_{QQ}{\cal{L}}_{t+1}E_{PP},$$
normalised in such a way that $E_{PQ}E_{QP} = E_{PP}$ and $E_{QP}E_{PQ}=E_{QQ}$.
Thus there is an algebra isomorphism $A_{t+1} \rightarrow {\cal{L}}_{t+1}$.

We then have the following theorem.
\begin{theorem}
\label{th:onelast}
Let $V$ be a finite dimensional irreducible
$U$-module such that $V^{\otimes t}$ is completely reducible for each $t \in \mathbb{Z}_{+}$ and
such that tensoring by $V$ is multiplicity free.
Then for any $t \in \mathbb{Z}_{+}$, the centraliser algebra 
${\cal{L}}_{t} = \mathrm{End}_{U} (V^{\otimes t})$
is isomorphic to the path algebra $A_{t}$ corresponding to the Bratteli diagram for $V^{\otimes t}$.
\end{theorem}

\end{subsection}

\end{section}

\begin{section}{Projections from $V^{\otimes t}$ onto its irreducible $U_{q}(\mathfrak{g})$-submodules}
\label{subsec:projectontome}
%\markright{\text{Projections onto finite dimensional irreducible $U_{q}(osp(1|2n))$-modules}}

Recall from Lemma \ref{lem:decompositionofthetensorproductoftwoirreduciblequantumospmodules}
that all tensor products of irreducible $U_{q}(\mathfrak{g})$-modules with integral dominant highest weights
at generic $q$ are completely reducible.
In this section we define projections from $V^{\otimes t}$ 
onto all the irreducible $U_{q}(\mathfrak{g})$-submodules 
$V_{\lambda} \subset V^{\otimes t}$, $\lambda \in {\cal{P}}^{+}$, using elements of ${\cal{C}}_{t}$.
Recall from Section \ref{subsec:100}
that ${\cal{P}}^{+}$ is the set of integral dominant weights.
No substantially new results appear in this section, however, we are not aware of 
this specific formulation of the projections in the literature. 

Let $V_{\mu}$ be a finite dimensional
irreducible $U_{q}(\mathfrak{g})$-module with highest weight $\mu \in {\cal{P}}^{+}$.
Since each weight space of $V$ is one-dimensional, $V_{\mu} \otimes V$ is multiplicity free, and from 
Lemma \ref{lem:decompositionofthetensorproductoftwoirreduciblequantumospmodules} 
we know the highest weights of the irreducible $U_{q}(\mathfrak{g})$-submodules in $V_{\mu} \otimes V$.
\begin{definition}
\label{def:twopointthreepointten}
We define ${\cal{P}}^{+}_{\mu} \subset {\cal{P}}^{+}$ to be the set such that for each
$\lambda \in {\cal{P}}^{+}_{\mu}$, $V_{\lambda}$ appears in $V_{\mu} \otimes V$
as an irreducible $U_{q}(\mathfrak{g})$-submodule.
\end{definition}
Now each $\lambda \in {\cal{P}}^{+}_{\mu}$ can only have one of the following three forms: 
$\mu, \mu+\epsilon_{i}, \mu - \epsilon_{i}$ for some $i$.
Thus
$${\cal{P}}^{+}_{\mu} \subseteq {\cal{P}}^{0}_{\mu} = \{ \mu, \mu \pm \epsilon_{i} \in {\cal{P}}^{+} | 
\ 1 \leq i \leq n \}.$$

\begin{definition}
Let $V_{\mu}$ be an irreducible $U_{q}(\mathfrak{g})$-module with highest weight
$\mu \in {\cal{P}}^{+}$, then 
$V_{\mu} \otimes V = \bigoplus_{\lambda \in {\cal{P}}^{+}_{\mu}} V_{\lambda}$.
Let $\big\{p_{\mu}[\lambda] \in \mathrm{End}_{U_{q}(\mathfrak{g})}(V_{\mu} \otimes V) 
| \ \lambda \in {\cal{P}}^{+}_{\mu} \big\}$ 
be a set of even maps
$$p_{\mu}[\lambda]: V_{\mu} \otimes V \rightarrow V_{\mu} \otimes V$$
such that
\begin{itemize}
\item[(i)] the image of $p_{\mu}[\lambda]$ is $V_{\lambda}$,
\item[(ii)]  $\big(p_{\mu}[\lambda]\big)^{2} = p_{\mu}[\lambda]$,
\item[(iii)] $p_{\mu}[\lambda] \cdot p_{\mu}[\nu] = \delta_{\lambda \nu} p_{\mu}[\lambda]$.
\end{itemize}
We call each such $p_{\mu}[\lambda]$ a projection.
\end{definition}

Recall that for each integral dominant $\lambda$, there exists an element
$v_{\lambda} \in \overline{U}_{q}^{+}(\mathfrak{g})$ defined in
 (\ref{eq:defofthevelements}) that acts on each
 vector in the finite dimensional irreducible $U_{q}(\mathfrak{g})$-module $V_{\lambda}$ as the
 multiplication by the scalar $q^{-(\lambda + 2\rho, \lambda)}$.
 
For each $\mu \in {\cal{P}}^{+}$ and each $\lambda \in {\cal{P}}^{+}_{\mu}$, define 
$p_{\mu}[\lambda] \in \mathrm{End}_{U_{q}(\mathfrak{g})}(V_{\mu} \otimes V)$ by
\begin{equation}
\label{eq:onetwofive}
p_{\mu}[\lambda] = (\pi_{\mu} \otimes \pi) 
\left( \prod_{\stackrel{\nu \in {\cal{P}}^{+}_{\mu}}{\nu \neq \lambda}} 
\frac{\Delta(v_{\xi})- q^{-(\nu + 2\rho, \nu)}}
%\chi_{\nu}(v_{\xi})
{ q^{-(\lambda + 2\rho, \lambda)} - q^{-(\nu + 2\rho, \nu)}    }
%{\chi_{\lambda}(v_{\xi})-\chi_{\nu}(v_{\xi})} 
\right),
\end{equation}
where $v_{\xi}$ is the element $v_{\lambda} \in \overline{U}_{q}^{+}(\mathfrak{g})$ with $\lambda = \xi$, 
for some integral dominant $\xi$
which is chosen so that $v_{\xi}$ 
acts as the multiplication by the scalar $q^{-(\nu + 2\rho, \nu)}$ on each vector in the
irreducible $U_{q}(\mathfrak{g})$-module $V_{\nu}$, for each $\nu \in {\cal{P}}^{+}_{\mu}$.  
For each integral dominant  $\mu$ there always exists at least one such $\xi$.
To see this, all we need is some
$E_{\xi}$ given by Eq. (\ref{eq:chatswoodtrain}):
$\displaystyle{E_{\xi} = \prod_{a=1}^{n} \sum_{b=p}^{s} (J_{a})^{b} \otimes P_{a}[b]},$
such that the element
$\displaystyle{E = \prod_{a=1}^{n} \sum_{b=p}^{s} P_{a}[b](J_{a})^{-b}}$
acts as the multiplication by the scalar $q^{-(\zeta, \zeta)}$
on each weight vector $w_{\zeta} \in V_{\nu} \subseteq V_{\mu} \otimes V$,
where $w_{\zeta}$ has the weight $\zeta$, 
and this is true for each $\nu \in {\cal{P}}^{+}_{\mu}$.

The element $E$ has this action whenever $s$ and $|p|$ are sufficiently large enough, 
and so all we need do is to choose some $\xi$ for which this is true. 
To do this, let $I_{\nu}$ be the set of distinct weights of the weight vectors of $V_{\nu}$
for each $\nu \in {\cal{P}}^{+}_{\mu}$,
then $(\zeta_{\nu}, \epsilon_{i}) \in \mathbb{Z}$ for each weight $\zeta_{\nu} \in I_{\nu}$ and 
each $i=1, \ldots, n$.
Let
$$m = \max{ \Big\{ |(\zeta_{\nu}, \epsilon_{i})| \in \mathbb{Z}_{+} \big| \ \zeta_{\nu} \in I_{\nu},
 \nu \in {\cal{P}}^{+}_{\mu}, i=1, \ldots, n  \Big\}  },$$
 then fixing $\xi = \sum_{i=1}^{n} m \epsilon_{i}$ yields elements $E_{\xi}$ and $E$ 
 with the desired properties.

Note that  $(\pi_{\mu} \otimes \pi) \Delta(v_{\xi})$ in (\ref{eq:onetwofive}) is diagonalisable as
$V_{\mu} \otimes V$ is completely reducible and
$\Delta(v_{\xi})$ acts on each irreducible $U_{q}(\mathfrak{g})$-submodule 
$V_{\nu} \subset V_{\mu} \otimes V$ as the 
multiplication by the scalar $q^{-(\nu + 2\rho, \nu)}$.

\begin{lemma}
\label{lem:lemmings}
The maps $p_{\mu}[\lambda]$ are well-defined and satisfy
\begin{itemize}
\item[(i)] $\big(p_{\mu}[\lambda]\big)^{2}= p_{\mu}[\lambda]$,
\item[(ii)] $p_{\mu}[\lambda] \cdot p_{\mu}[\nu] = \delta_{\lambda \nu} p_{\mu}[\lambda]$,
\item[(iii)] $\displaystyle{\sum_{\lambda \in {\cal{P}}^{+}_{\mu}} p_{\mu}[\lambda] = 
              \mathrm{id}_{V_{\mu} \otimes V}}$.
\end{itemize}
\end{lemma}
\begin{proof} 
If $\alpha$ and $\beta$ are the highest weights of irreducible 
$U_{q}(\mathfrak{g})$-submodules in $V_{\mu} \otimes V$, then
$(\alpha + 2\rho, \alpha) = (\beta + 2\rho, \beta)$ implies that $\alpha = \beta$.
Then $p_{\mu}[\lambda]$ is well defined as tensoring by $V$ is multiplicity free and $q$ is not a root of unity.  
The proof of (i) follows from the result that
$\big(p_{\mu}[\lambda]\big)^{2}(V_{\mu} \otimes V) = p_{\mu}[\lambda] (V_{\lambda}) = V_{\lambda}$.
For (ii) the case $\lambda = \nu$ reduces to (i), and for $\lambda \neq \nu$ we have
\begin{eqnarray*}
\lefteqn{
p_{\mu}[\lambda] \cdot p_{\mu}[\nu] } \\
& &  =
(\pi_{\mu} \otimes \pi)  \left(
\prod_{\stackrel{\lambda' \in {\cal{P}}^{+}_{\mu}}{\lambda' \neq \lambda}}
\frac{\Delta(v_{\xi})-q^{-(\lambda' + 2\rho, \lambda')}}
{q^{-(\lambda+ 2\rho, \lambda)}-q^{-(\lambda'+2\rho, \lambda')}}
\prod_{\stackrel{\nu' \in {\cal{P}}^{+}_{\mu}}{\nu' \neq \nu}}
\frac{\Delta(v_{\xi})-q^{-(\nu' + 2\rho, \nu')}}
{q^{-(\nu + 2\rho, \nu)}- q^{-(\nu' + 2\rho, \nu')}} \right) = 0. 
\end{eqnarray*}
\begin{itemize}
\item[(iii)] $\displaystyle{\sum_{\lambda \in{\cal{P}}^{+}_{\mu}}}
p_{\mu}[\lambda] \left(V_{\mu} \otimes V \right)=
\displaystyle{\bigoplus_{\lambda \in {\cal{P}}^{+}_{\mu}} } V_{\lambda} = V_{\mu} \otimes V$.
\end{itemize}
\end{proof}
\noindent
Note that 
$(\pi_{\mu} \otimes \pi) \displaystyle{\prod_{\lambda \in {\cal{P}}^{+}_{\mu}} 
\Big( \Delta(v_{\xi})-q^{-(\lambda + 2\rho, \lambda)} \Big) } = 0$.

We introduce some notation.
Let ${\cal{T}}^{t}$ be the set of tableaux of length $t$ derived from the Bratteli diagram for
$V^{\otimes t}$.  Let
$$i^{t} = (0, s_{1}, \ldots, s_{t}) \in {\cal{T}}^{t}.$$
We write $\lambda_{i}^{t} = i^{t}$ where $\lambda = s_{t}$.

Let $i^{t} \in {\cal{T}}^{t}$ and $s_{j}$, $s_{j+1} \in i^{t}$. Define a map 
$$p^{t-(j+1)}_{s_{j}}[s_{j+1}]:\left(V_{s_{j}} \otimes V\right)\otimes V^{\otimes
(t-(j+1))} \rightarrow V_{s_{j+1}} \otimes V^{\otimes (t-(j+1))}$$
 by
$$p^{t-(j+1)}_{s_{j}}[s_{j+1}] = p_{s_{j}}[s_{j+1}] \otimes \mathrm{id}^{\otimes
(t-(j+1))}.$$
\begin{lemma}
\label{lem:velvetvelvet}
The map   $p^{t-(j+1)}_{s_{j}}[s_{j+1}]$ satisfies
\begin{itemize}
\item[(i)]
$\left(p_{s_{j}}^{t-(j+1)}[s_{j+1}]\right)^{2} = p_{s_{j}}^{t-(j+1)}[s_{j+1}]$,
\item[(ii)]
$p_{s_{j}}^{t-(j+1)}[s_{j+1}] \cdot p_{s_{j}}^{t-(j+1)}[r_{j+1}]
= \delta_{s_{j+1},{r_{j+1}}} p_{s_{j}}^{t-(j+1)}[s_{j+1}],$
\item[(iii)]
$\displaystyle{\sum_{s_{j+1} \in {\cal{P}}^{+}_{s_{j}}}} p_{s_{j}}^{t-(j+1)}[s_{j+1}] 
= \mathrm{id}_{V_{s_{j}} \otimes V^{\otimes (t-j)}}$.
\end{itemize}
\end{lemma}
\begin{proof}
The proofs of parts (i) and (ii) follow from Lemma \ref{lem:lemmings} (i) and (ii), 
respectively.
The proof of (iii) follows from  Lemma \ref{lem:lemmings} (iii): explicitly, we have
$$\displaystyle{\sum_{s_{j+1} \in {\cal{P}}^{+}_{s_{j}}}
p_{s_{j}}^{t-(j+1)}[s_{j+1}]\cdot \left(V_{s_{j}}\otimes V \right)
\otimes V^{\otimes (t-(j+1))}
 = V_{s_{j}} \otimes V \otimes V^{\otimes t-(j+1)}}.$$
\end{proof}

\begin{definition}
\label{def:sartor}
Let $\tilde{p}_{i}^{t}[\lambda] \in \mathrm{End}_{U_{q}(\mathfrak{g})}(V^{\otimes t})$ be a map
$\tilde{p}_{i}^{t}[\lambda]: V^{\otimes t} \rightarrow V_{\lambda} \subset V^{\otimes t}$ 
defined by
$$\tilde{p}_{i}^{t}[\lambda] =
p^{0}_{s_{t-1}}[\lambda]p^{1}_{s_{t-2}}[s_{t-1}]\cdots
p^{t-2}_{\epsilon_{1}}[s_{2}],$$
where $\lambda_{i}^{t} \in {\cal{T}}^{t}$. 
We say that {\emph{$\tilde{p}_{i}^{t}[\lambda]$ projects from $V^{\otimes t}$ onto
$V_{\lambda}$ by the path $\lambda_{i}^{t} \in {\cal{T}}^{t}$}} and 
we call $\tilde{p}_{i}^{t}[\lambda]$ a {\emph{path projection of length $t$}}.
\end{definition}
\begin{lemma}
\label{lem:greatgnats}
%Let $\tilde{p}_{i}^{t}[\lambda] \in \mathrm{End}_{U_{q}(\mathfrak{g})}(V^{\otimes t})$ be a map from Definition
%\ref{def:sartor}.  
The map $\tilde{p}_{i}^{t}[\lambda]$ satisfies
\begin{itemize}
\item[(i)] $\big(\tilde{p}_{i}^{t}[\lambda]\big)^{2} = \tilde{p}_{i}^{t}[\lambda]$,
\item[(ii)] $\tilde{p}_{i}^{t}[\lambda] \cdot \tilde{p}_{j}^{t}[\lambda] = 
\left\{  
\begin{array}{ll}
0,                          & \mbox{if } i^{t} \neq j^{t}, \\
\tilde{p}_{i}^{t}[\lambda], & \mbox{if } i^{t} = j^{t},
\end{array}
\right. $
\item[(iii)]  $\tilde{p}_{i}^{t}[\lambda]\cdot \tilde{p}_{j}^{t}[\mu]=0$ if $\lambda \neq \mu$.
\end{itemize}
Furthermore, the map $\displaystyle{P_{t} = \sum_{i^{t} \in {\cal{T}}^{t}} \tilde{p}_{i}^{t}[\lambda]}$
is the identity on $V^{\otimes t}$.
\end{lemma}
\begin{proof}
\noindent
\begin{itemize}
\item[(i)] This follows from Lemma \ref{lem:velvetvelvet} (i).
\item[(ii)] For $i^{t} = j^{t}$ the case reduces to (i), let $i^{t} \neq j^{t}$ where
\begin{eqnarray*}
i^{t} & = & (0, s_{1}, s_{2}, \ldots, s_{k-1}, s_{k}, r_{k+1}, \ldots, r_{t-2}, r_{t-1}, \lambda), \\
j^{t} & = & (0, s_{1}, s_{2}, \ldots, s_{k-1}, s_{k}, s_{k+1}, \ldots, s_{t-2}, s_{t-1}, \lambda).
\end{eqnarray*}
Now $i^{t}$, $j^{t} \in {\cal{T}}^{t}$ and $r_{k+1} \neq s_{k+1}$ for some $2 \leq k+1 \leq t$, then
\begin{eqnarray*}  
\tilde{p}_{i}^{t}[\lambda]\cdot \tilde{p}_{j}^{t}[\lambda] & = &
p^{0}_{r_{t-1}}[\lambda]p^{1}_{r_{t-2}}[r_{t-1}]\cdots 
p^{t-k-1}_{s_{k}}[r_{k+1}]p^{t-k}_{s_{k-1}}[s_{k}]
p^{t-k+1}_{s_{k-2}}[s_{k-1}]\cdots
p^{t-2}_{s_{1}}[s_{2}] \\
&  &  \times
p^{0}_{s_{t-1}}[\lambda]p^{1}_{s_{t-2}}[s_{t-1}]\cdots 
p^{t-k-1}_{s_{k}}[s_{k+1}]p^{t-k}_{s_{k-1}}[s_{k}]
p^{t-k+1}_{s_{k-2}}[s_{k-1}] \cdots
p^{t-2}_{s_{1}}[s_{2}]  \\
& = & p^{0}_{r_{t-1}}[\lambda]p^{0}_{s_{t-1}}[\lambda]
p^{1}_{r_{t-2}}[r_{t-1}]p^{1}_{s_{t-2}}[s_{t-1}]
\cdots p^{t-k-1}_{s_{k}}[r_{k+1}]p^{t-k-1}_{s_{k}}[s_{k+1}]  \\
& &  \times  p^{t-k}_{s_{k-1}}[s_{k}]p^{t-k}_{s_{k-1}}[s_{k}] 
p^{t-k+1}_{s_{k-2}}[s_{k-1}]p^{t-k+1}_{s_{k-2}}[s_{k-1}]\cdots
p^{t-2}_{s_{1}}[s_{2}]p^{t-2}_{s_{1}}[s_{2}]  \\
& = & 0,
\end{eqnarray*}
as $p^{t-k-1}_{s_{k}}[r_{k+1}] \cdot p^{t-k-1}_{s_{k}}[s_{k+1}]=0$.
\item[(iii)]  Assume that
\begin{eqnarray*}
i^{t} & = & (0, s_{1}, s_{2}, \ldots, s_{k-1}, s_{k}, r_{k+1}, \ldots, r_{t-2}, r_{t-1}, \lambda), \\
j^{t} & = & (0, s_{1}, s_{2}, \ldots, s_{k-1}, s_{k}, s_{k+1}, \ldots, s_{t-2}, s_{t-1}, \mu),
\end{eqnarray*}
where $r_{k+1}\neq s_{k+1}$ for some $2 \leq k+1 \leq t$.
The calculations are similar to those of (ii) and we have
$\tilde{p}^{t}_{i}[\lambda]\cdot \tilde{p}^{t}_{j}[\mu]= 0$.
\end{itemize}
The last claim
follows inductively from the result that
$\displaystyle{\sum_{\lambda \in {\cal{P}}^{+}_{\mu}} p_{\mu}[\lambda] = \mathrm{id}_{V_{\mu} \otimes V}}$.
\end{proof}

Recall that ${\cal{C}}_{t}$ is the algebra over $\mathbb{C}$ generated by the elements
$$\left\{\check{R}_{i}^{\pm 1} \in \mathrm{End}_{U_{q}(\mathfrak{g})}(V^{\otimes t}) |
 \ 1 \leq i \leq t-1\right\}.$$ 
\begin{proposition}
For each path $\lambda_{i}^{t} \in {\cal{T}}^{t}$, $\tilde{p}_{i}^{t}[\lambda] \in {\cal{C}}_{t}$.
\end{proposition}
\begin{proof}
We prove the proposition inductively.  
Firstly, for some appropriately chosen integral dominant weight $\xi$, 
$$(\pi \otimes \pi)\Delta(v_{\xi}) = 
(\pi \otimes \pi)\left[(v_{\xi} \otimes v_{\xi})\left(R^{T}_{\xi} R_{\xi}\right)^{-1}\right]
 = q^{-2(\epsilon_{1} + 2\rho, \epsilon_{1})} \check{\cal{R}}^{-2} \in {\cal{C}}_{2}.$$
Now assume that $\tilde{p}^{(t-1)}_{i}[\mu] \in {\cal{C}}_{(t-1)}$ where
$\tilde{p}^{(t-1)}_{i}[\mu]$ is a path projection 
$\tilde{p}^{(t-1)}_{i}[\mu]: V^{\otimes (t-1)} \rightarrow V_{\mu}$ and
$V_{\mu}$ is an irreducible $U_{q}(\mathfrak{g})$-submodule of $V^{\otimes (t-1)}$.
We will show that  $(\pi_{\mu} \otimes \pi)\Delta(v_{\zeta})$ is an element of ${\cal{C}}_{t}$ 
for some appropriately chosen $\zeta$.
Let $\zeta$ be an integral dominant weight such that the element 
$v_{\zeta} \in \overline{U}_{q}^{+}(\mathfrak{g})$ acts as the multiplication by the scalar
$q^{-(\lambda + 2\rho, \lambda)}$ on each vector in the irreducible $U_{q}(\mathfrak{g})$-submodule
$V_{\lambda} \subset V_{\mu} \otimes V$ for each $\lambda \in {\cal{P}}^{+}_{\mu}$.
Now 
\begin{eqnarray*}
\lefteqn{
(\pi_{\mu} \otimes \pi) \Delta(v_{\zeta}) } \\
 & = & (\pi_{\mu} \otimes \pi) \left[ (v_{\zeta} \otimes v_{\zeta}) \left(R^{T}_{\zeta}R_{\zeta}\right)^{-1} \right] \\
& = & q^{-(\mu + 2\rho, \mu)-(\epsilon_{1} + 2\rho, \epsilon_{1})} 
      (\tilde{p}^{(t-1)}_{i}[\mu] \otimes \mathrm{id}) \left(\pi^{\otimes (t-1)} \otimes \pi \right) 
      \left(\Delta^{(t-2)} \otimes \mathrm{id}\right) 
      \left(R^{T}_{\zeta}R_{\zeta}\right)^{-1} \\
& = & q^{-(\mu + 2\rho, \mu)-(\epsilon_{1} + 2\rho, \epsilon_{1})}  
      (\tilde{p}^{(t-1)}_{i}[\mu] \otimes \mathrm{id}) 
      \check{R}_{t-1}^{-1} \check{R}_{t-2}^{-1} \cdots \check{R}_{1}^{-1} \check{R}_{1}^{-1} \cdots
       \check{R}_{t-2}^{-1}\check{R}_{t-1}^{-1},
\end{eqnarray*}
where we have used the following result from (\ref{starequationlikeearth}): (writing $R = R_{\zeta}$)
$$ \left( \pi^{\otimes (t-1)} \otimes \pi \right) \left(\Delta^{(t-2)} \otimes \mathrm{id} \right)R = 
\left(\pi^{\otimes (t-1)} \otimes \pi \right) R_{1t} R_{2t} \cdots R_{(t-1)t}.$$ 

\end{proof}

\end{section}

\begin{section}{Matrix units for ${\cal{C}}_{t}$}
\label{subsec:rhapsodyinred}
%\markright{\text{Matrix units for ${\cal{C}}_{t}$}}

It is clear that the Bratteli diagram for $V^{\otimes t}$ is multiplicity free,
as tensoring by the fundamental $U_{q}(\mathfrak{g})$-module $V$ is multiplicity free.
It follows then from Theorem \ref{th:onelast} that the centraliser algebra 
${\cal{L}}_{t} = \mathrm{End}_{U_{q}(\mathfrak{g})}(V^{\otimes t})$ is isomorphic to 
the path algebra $A_{t}$ obtained from the Bratteli diagram for $V^{\otimes t}$.
Clearly, we have the inclusion ${\cal{C}}_{t} \subseteq {\cal{L}}_{t}$.  
The aim of this section is to show that ${\cal{C}}_{t}$ and ${\cal{L}}_{t}$ are in fact equal:
\begin{theorem}
\label{th:areallybigthe(199)}
The centraliser algebra ${\cal{L}}_{t} = \mathrm{End}_{ U_{q}(\mathfrak{g}) }\big(V^{\otimes t}\big)$ 
is generated by the elements
$$\big\{ \check{R}_{i}^{\pm 1} \in \mathrm{End}_{U_{q}(\mathfrak{g})}\big(V^{\otimes t}\big) \big| \ i=1, 2,
\ldots, t-1 \big\}.$$
\end{theorem}

To prove this theorem we firstly partition the matrix units in $A_{t}$ 
into two groups:
the {\emph{projectors}}  $\{E_{SS} \in A_{t} | \ (S,S) \in \Omega^{t}\}$ and the 
{\emph{intertwiners}}  $\{E_{ST} \in A_{t} | \ (S,T) \in \Omega^{t}, S \neq T \}$
and we use an invertible homomorphism to map matrix units in $A_{t}$ 
to matrix units in ${\cal{C}}_{t}$.

Recall that $V^{\otimes t}$ is completely reducible.
Each matrix unit in ${\cal{C}}_{t}$ corresponding to a
projector in $A_{t}$ projects down from $V^{\otimes t}$ onto an irreducible 
$U_{q}(\mathfrak{g})$-submodule $V_{\lambda} \subseteq V^{\otimes t}$. 
Each matrix unit in ${\cal{C}}_{t}$ corresponding to an
intertwiner in $A_{t}$ maps between isomorphic irreducible 
$U_{q}(\mathfrak{g})$-submodules of $V^{\otimes t}$.

Recall that the homomorphism $\Upsilon: g_{i} \mapsto -\check{R}_{i}$ 
given in Lemma \ref{eq:peacockinblackbeansauce} yields a representation of
$BW_{t}(-q^{2n},q)$ in ${\cal{C}}_{t}$.
In Subsection \ref{subsect:matrixbirmanwenzlstuff} we will write down the matrix units in a semisimple
quotient of $BW_{t}(-q^{2n},q)$ that map via $\Upsilon$ onto the 
projectors and intertwiners in ${\cal{C}}_{t}$.
We will do this for the intertwiners, 
but we choose to define the projectors more straightforwardly using our previous work.

The projections $E_{SS}$ that project down
from $V^{\otimes t}$ onto irreducible $U_{q}(\mathfrak{g})$-submodules 
$V_{shp(S)} \subset V^{\otimes t}$ that we defined in Section \ref{subsec:projectontome}, 
are elements of ${\cal{C}}_{t}$, and satisfy $(E_{SS})^{2} = E_{SS}$ and
$\sum_{S \in {\cal{T}}^{t}} E_{SS} = \mathrm{id}_{V^{\otimes t}}$. 
We fix the projectors in ${\cal{C}}_{t}$ by  
mapping the projector $E_{SS} \in A_{t}$ to $\tilde{p}_{i}^{t}[\lambda] \in {\cal{C}}_{t}$, where
$\lambda_{i}^{t} = S \in {\cal{T}}^{t}$ is a path of length $t$:
$E_{SS} \leftrightarrow \tilde{p}_{i}^{t}[\lambda]$.

All we need do now is
construct the matrix units in ${\cal{C}}_{t}$ corresponding to the intertwiners in $A_{t}$.
We denote the matrix unit in ${\cal{C}}_{t}$ corresponding to 
$E_{MP} \in A_{t}$ also by $E_{MP}$.  

\begin{subsection}{Matrix units in $BW_{t}(-q^{2n},q)$}
\label{subsect:matrixbirmanwenzlstuff}

In this subsection, we say that an algebra $B$ is semisimple if it is isomorphic to a direct sum of
matrix algebras, ie $B \cong \bigoplus_{i \in I} M_{b_{i}}(\mathbb{C})$, 
where $M_{b_{i}}(\mathbb{C})$ is the algebra of $b_{i} \times b_{i}$ matrices with complex entries.
The algebra $BW_{t}(-q^{2n},q)$ is not semisimple at generic $q$ \cite[Cor. 5.6]{w2} but
Ram and Wenzl have constructed matrix units
for the semisimple Birman-Wenzl-Murakami algebra $BW_{t}$ defined over $\mathbb{C}(r,q)$
(the field of rational functions in $r$ and $q$) for indeterminates $r$ and $q$ \cite{rw}.
By replacing the indeterminates $r$ and $q$ with the complex numbers $-q^{2n}$ and $q$,
respectively, we obtain matrix units in some appropriate semisimple quotient of $BW_{t}(-q^{2n},q)$.
By then applying the homomorphism $\Upsilon$ to these matrix units, we obtain matrix units in
${\cal{C}}_{t}$.

Before doing this, let us recall the definition of a Young diagram 
and discuss a relation between certain Young diagrams
and integral dominant highest weights of irreducible $U_{q}(osp(1|2n))$-modules.
For each non-negative integer $m$, there exists 
a Young diagram for each partition of $m$.
Let $m = m_{1} + m_{2} + \cdots + m_{l}$ be a partition of $m$, where 
$m_{i} - m_{i+1} \in \mathbb{Z}_{+}$ for each
$i=1, 2, \ldots, l-1$ and $m_{l} \in \mathbb{Z}_{+}$.
The Young diagram representing this partition is a collection of $m$ boxes arranged in $l$
left-aligned rows where the $i^{th}$ row from the top  contains exactly $m_{i}$ boxes. 
For $m \geq 1$, 
let $c_{i}$, $i = 1, 2, \ldots, m_{1}$, 
be the number of boxes in the $i^{th}$ column from the left in the Young diagram,
then $c_{i} - c_{i+1} \in \mathbb{Z}_{+}$ for each $i=1, 2, \ldots, m_{1}-1$
and $c_{m_{1}} \in \{1, 2, \ldots, l\}$.

Recall that an integral dominant highest weight $\lambda$
of an irreducible $U_{q}(\mathfrak{g})$-module $V_{\lambda}$ has
the form $\lambda = \sum_{i=1}^{n} \lambda_{i} \epsilon_{i} \in {\cal{P}}^{+}$ where 
$\lambda_{i} - \lambda_{i+1} \in \mathbb{Z}_{+}$ for each $i=1, 2, \ldots, n-1$ and
$\lambda_{n} \in \mathbb{Z}_{+}$. 
We can use a Young diagram to label $\lambda$:
this Young diagram consists of $\sum_{i=1}^{n} \lambda_{i}$ boxes arranged in $n$ left-aligned rows,
where the $i^{th}$ row contains exactly $\lambda_{i}$ boxes.

Let $\mu$ be a Young diagram containing no more than $n$ rows of boxes and let
$\mu_{i}$ be the number of boxes in the $i^{th}$ row from the top 
for each $i=1, 2, \ldots, n$.
We can use the Young diagram $\mu$ to label the integral dominant highest weight 
$\sum_{i=1}^{n} \mu_{i} \epsilon_{i} \in {\cal{P}}^{+}$ of an irreducible
representation of $U_{q}(osp(1|2n))$.

\begin{subsection}{The algebra $BW_{t}$}

The Birman-Wenzl-Murakami algebra $BW_{t}$, 
with $r$ and $q$ indeterminates, is semisimple \cite[Thm. 3.5]{w2}.
To discuss the structure of $BW_{t}$, we introduce the Young lattice.
For later purposes we note that $BW_{t}$ is equipped with a functional 
$\mathrm{tr}: BW_{t} \rightarrow \mathbb{C}(r,q)$ which satisfies, 
amongst other relations \cite[Lem. 3.4 (d)]{w2},
\begin{equation}
\label{eq:tracefunctionalnew}
\mathrm{tr}(a \chi b) = \mathrm{tr}(\chi) \mathrm{tr}(ab), 
\hspace{10mm} \forall a, b \in BW_{t-1}, \hspace{5mm} \chi \in \{g_{t-1}, e_{t-1}\},
\end{equation}
where we regard each element of $BW_{t-1}$ as an element of $BW_{t}$
under the canonical inclusion. 

\begin{remark}
The functional $\mathrm{tr}$ in (\ref{eq:tracefunctionalnew}) can be used to define 
the Kauffman link invariant following \cite[p. 404]{w2}.
\end{remark}

The Young lattice is the following infinite graph  \cite[Sec. 1]{w2}.  
The vertices of the Young lattice are the Young diagrams;
the vertices are grouped into levels so that each Young diagram with exactly $t$ boxes
labels a vertex on the $t^{th}$ level of the Young lattice.
The edges of the Young lattice are completely determined as follows: 
a vertex $\lambda$ on the $t^{th}$ level is connected to a vertex $\mu$ on the 
$(t+1)^{st}$ level by one edge if and only if $\lambda$ and $\mu$ differ by exactly one box.
We show the Young lattice up to the $4^{th}$ level in Figure \ref{fig:young}, 
where the circle represents the Young diagram with no boxes.
We say that the level containing the Young diagram with no boxes is the $0^{th}$ level.
Note that the Young lattice is (apart from the $0^{th}$ level) 
identical to the Bratteli diagram for the 
sequence of inclusions of group algebras of the symmetric group:
 $\mathbb{C}S_{1} \subset \mathbb{C}S_{2} \subset  \mathbb{C}S_{3} \subset \cdots $.

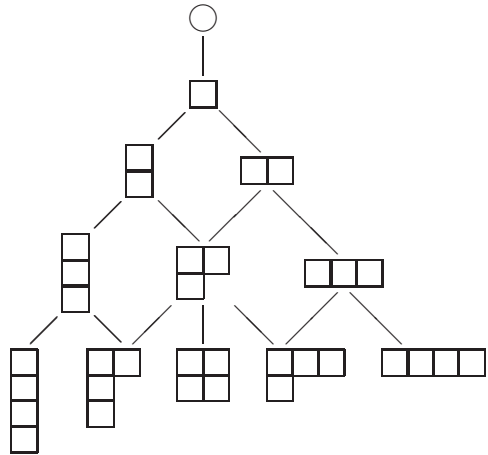
\begin{figure}[hbt]
\begin{center}
\setlength{\unitlength}{0.17mm}
\begin{picture}(350,400) \large\sf
\thinlines
%\linethickness{0.05mm}

\put(150,370){\circle{20}} 
\put(150,325){\line(0,1){30}} 

\put(140,300){\framebox(20,20){}}

\put(90,230){\framebox(20,20){}}
\put(90,250){\framebox(20,20){}}

\put(180,240){\framebox(20,20){}}
\put(200,240){\framebox(20,20){}}

\put(130,170){\framebox(20,20){}}
\put(150,170){\framebox(20,20){}}
\put(130,150){\framebox(20,20){}}

\put(230,160){\framebox(20,20){}}
\put(250,160){\framebox(20,20){}}
\put(270,160){\framebox(20,20){}}

\put(40,180){\framebox(20,20){}}
\put(40,140){\framebox(20,20){}}
\put(40,160){\framebox(20,20){}}

\put(0,90){\framebox(20,20){}}
\put(0,70){\framebox(20,20){}}
\put(0,50){\framebox(20,20){}}
\put(0,30){\framebox(20,20){}}
\put(15,115){\line(1, 1){21}} 

\put(60,90){\framebox(20,20){}}
\put(60,70){\framebox(20,20){}}
\put(60,50){\framebox(20,20){}}
\put(80,90){\framebox(20,20){}}
\put(65,137){\line(1,-1){21}} 
\put(95,115){\line(1,1){30}} 

\put(130,90){\framebox(20,20){}}
\put(130,70){\framebox(20,20){}}
\put(150,70){\framebox(20,20){}}
\put(150,90){\framebox(20,20){}}
\put(150,115){\line(0,1){30}} 

\put(200,90){\framebox(20,20){}}
\put(200,70){\framebox(20,20){}}
\put(220,90){\framebox(20,20){}}
\put(240,90){\framebox(20,20){}}
\put(175,145){\line(1,-1){30}}

\put(290,90){\framebox(20,20){}}
\put(310,90){\framebox(20,20){}}
\put(330,90){\framebox(20,20){}}
\put(350,90){\framebox(20,20){}}
\put(265,155){\line(1,-1){40}} 
\put(215,115){\line(1, 1){40}}

\put(115,275){\line(1, 1){21}} 

\put(163,297){\line(1,-1){32}} 
\put(115,227){\line(1,-1){32}} 
\put(155,195){\line(1, 1){40}} 
\put(65,205){\line(1, 1){21}} 

\put(205,235){\line(1,-1){50}} 

\end{picture}
\caption{The Young lattice up to the $4^{th}$ level}  \label{fig:young}
\end{center}
\end{figure}

For each $t$, let $Y_{t}$ be the set of vertices on the $t^{th}$ level of the Young lattice and define
$$\Gamma_{t} = \bigcup_{\stackrel{k \in \mathbb{Z}_{+}}{t - 2k \geq 0}} Y_{t-2k}.$$ 
Then $BW_{t}$ is isomorphic to a direct sum of matrix algebras 
\cite[Thm. 3.5]{w2}:
$$BW_{t} \cong \bigoplus_{\mu \in \Gamma_{t}} M_{b_{\mu}}(\mathbb{C}).$$
Ram and Wenzl defined matrix units for $BW_{t}$  \cite{rw} which we explicitly write down below
for completeness.

To label the matrix units of $BW_{t}$, we need to discuss the Bratteli diagram of
$BW_{t}$, which is the following graph. 
The vertices of the Bratteli diagram for $BW_{t}$ 
are divided into levels; for each $s=0, 1, \ldots, t$, the vertices on the $s^{th}$ level 
are precisely the elements of $\Gamma_{s}$.
The edges are as follows: a vertex
$\mu$ on the $s^{th}$ level is connected to a vertex $\lambda$ on the $(s+1)^{st}$ level if 
and only if $\mu$ and $\lambda$ differ by exactly one box.
We show the Bratteli diagram for $BW_{t}$ up to the $4^{th}$ level in 
Figure \ref{fig:BWbrattelidiagram}.

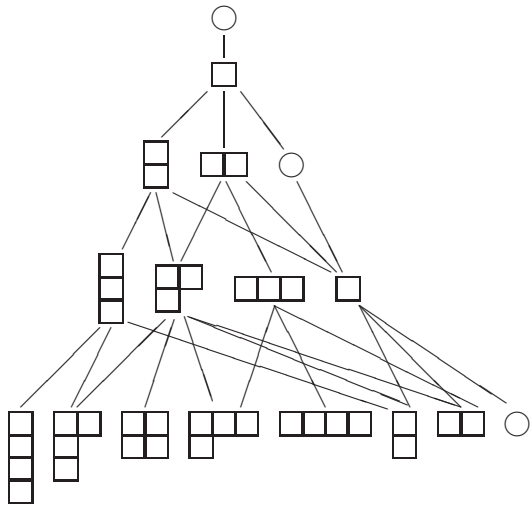
\begin{figure}[hbt]
\label{fig:BWbrattelidiagram}
\begin{center}
\setlength{\unitlength}{0.15mm}
\begin{picture}(460,440) \large\sf
\thinlines
%\linethickness{0.05mm}

\put(190,430){\circle{20}} 

\put(190,395){\line(0,1){20}} 

\put(180,370){\framebox(20,20){}}

\put(120,300){\framebox(20,20){}}
\put(120,280){\framebox(20,20){}}

\put(170,290){\framebox(20,20){}}
\put(190,290){\framebox(20,20){}}

\put(250,300){\circle{20}} 

\put(80,160){\framebox(20,20){}}
\put(80,180){\framebox(20,20){}}
\put(80,200){\framebox(20,20){}}
\put(100,225){\line(1,2){25}} 
\put(130,275){\line(1,-4){15}} 

\put(130,170){\framebox(20,20){}}
\put(130,190){\framebox(20,20){}}
\put(150,190){\framebox(20,20){}}
\put(152,215){\line(1,2){35}} 
\put(192,285){\line(1,-2){40}}

\put(200,180){\framebox(20,20){}}
\put(220,180){\framebox(20,20){}}
\put(240,180){\framebox(20,20){}}

\put(290,180){\framebox(20,20){}}
\put(210,285){\line(1,-1){80}} 
\put(255,285){\line(1,-2){40}} 
\put(145,275){\line(2,-1){140}}

\put(0,0){\framebox(20,20){}}
\put(0,20){\framebox(20,20){}}
\put(0,40){\framebox(20,20){}}
\put(0,60){\framebox(20,20){}}

\put(10,85){\line(1,1){70}} 
\put(135,325){\line(1,1){40}} 
\put(190,315){\line(0,1){50}} 
\put(205,365){\line(3,-4){38}}

\put(60,60){\framebox(20,20){}}
\put(40,20){\framebox(20,20){}}
\put(40,40){\framebox(20,20){}}
\put(40,60){\framebox(20,20){}}
\put(55,85){\line(1,2){35}} 
\put(60,85){\line(1,1){78}}

\put(100,60){\framebox(20,20){}}
\put(100,40){\framebox(20,20){}}
\put(120,60){\framebox(20,20){}}
\put(120,40){\framebox(20,20){}}
\put(120,85){\line(1,3){26}}

\put(160,60){\framebox(20,20){}}
\put(160,40){\framebox(20,20){}}
\put(180,60){\framebox(20,20){}}
\put(200,60){\framebox(20,20){}}
\put(155,165){\line(1,-3){25}} 
\put(205,85){\line(1,3){30}} 
\put(105,160){\line(3,-1){230}}

\put(240,60){\framebox(20,20){}}
\put(260,60){\framebox(20,20){}}
\put(280,60){\framebox(20,20){}}
\put(300,60){\framebox(20,20){}}
\put(235,175){\line(1,-2){45}}

\put(340,60){\framebox(20,20){}}
\put(340,40){\framebox(20,20){}}
\put(158,165){\line(5,-2){195}} 
\put(310,175){\line(1,-1){90}} 
\put(310,175){\line(1,-2){45}} 
\put(310,175){\line(3,-2){130}} 

\put(380,60){\framebox(20,20){}}
\put(400,60){\framebox(20,20){}}
\put(160,165){\line(3,-1){240}} 
\put(235,175){\line(2,-1){180}} 

\put(450,70){\circle{20}} 

\end{picture}
\caption{The Bratteli diagram for
$BW_{t}$ up to the $4^{th}$ level}  
\end{center}
\end{figure}

We say that $R$ is a path of length $t$ in
the Bratteli diagram of $BW_{t}$ 
if $R$ is a sequence of $t+1$ Young diagrams:
$R = ([0], [1], \ldots, r_{t})$ where $r_{s} \in \Gamma_{s}$ for each $s=0, 1, \ldots, t$
and where $r_{i}$ is connected to $r_{i+1}$ for each $0 \leq i \leq t-1$.
We say that $shp(R) = r_{t}$.
Let $\Omega_{t}$ be the set of pairs $(R,S)$ of paths of length $t$ 
in the Bratteli diagram of $BW_{t}$ where $r_{t} = s_{t}$, ie those paths
$R$ and $S$ of length $t$ where $shp(R)=shp(S)$.

Ram and Wenzl \cite{rw} wrote down a basis  
$\{e_{ST} \in BW_{t} | \ (S,T) \in \Omega_{t} \}$ of $BW_{t}$ 
which is also a set of matrix units. We recall this set of matrix units below.

Let us fix some notation.  
Given a sequence $T =(0, s_{1}, \ldots, s_{t})$,
we fix $T' = (0, s_{1}, \ldots, s_{t-1})$.
If $T$ is a path of length $t$, then $T'$ is the path of length $t-1$ obtained by removing the last
vertex and edge of $T$.

Before defining the matrix units of $BW_{t}$ we
define some `pre-matrix units' that we will employ in defining the matrix units.
Let $T$ be a path of length $t$ in the Bratteli diagram for $BW_{t}$
such that $shp(T)$ has $t$ boxes.
We can identify $T$ with a standard tableau containing the numbers $1, 2, \ldots, t$
in a canonical way.
We do this by  placing the number $1$ in the top left hand box of  
$shp(T)$ and we then fill each box of $shp(T)$ with increasing numbers
according to the path $T$ \cite[Sec. 4.2]{tw}.

For each path $T$ of length $t$ in the Bratteli diagram for $BW_{t}$,
 we define the number $d(T,i)$ for each $i=1, 2, \ldots, t-1$ by
\begin{equation}
\label{eq:brontethewhalecatslothdogmouse}
d(T,i) = c(i+1)-c(i)-r(i+1)+r(i),
\end{equation}
where $c(j)$ and $r(j)$ denote the column and row, respectively, 
of the box containing the number $j$ in the standard tableau corresponding to $T$.
For each $d \in \mathbb{Z} \backslash \{0\}$, we define
$$b_{d}(q) = \frac{q^{d}(1-q)}{1-q^{d}}.$$

Let $T$ be a path of length $t$ in the Bratteli diagram for $BW_{t}$. 
Firstly fix the pre-matrix unit $o_{[1]}=1 \in BW_{t}$. 
Let $R$ be a path of length $t-1$ defined by $R = T'$ and inductively define
$$o_{T} = \prod_{S} \frac{o_{R} g_{t-1} o_{R} - b_{d(S,t-1)}(q^{2}) o_{R} }
         {b_{d(T,t-1)}(q^{2}) - b_{d(S,t-1)}(q^{2})} \in  BW_{t},$$
where the product is over all paths $S$ of length $t$ where $shp(S)$ contains $t$ boxes such that
$S \neq T$ and $S' = R$. We write $o_{TT}=o_{T}$.

If $M$ and $P$ are paths of length $t$ in the Bratteli diagram for $BW_{t}$
where $(M,P) \in \Omega_{t}$ and $shp(M) = shp(P)$ has exactly $t$ boxes
and $shp(M') = shp(P')$, then we define
$$o_{MP} = o_{M' P'} o_{PP}.$$

Now, if $M$ and $P$ are paths of length $t$ where $(M,P) \in \Omega_{t}$ 
and $shp(M) = shp(P)$ has exactly $t$ boxes and
$shp(M') \neq shp(P')$, then the pre-matrix unit $o_{MP}$ is defined more intricately.
To define $o_{MP}$, choose paths $\overline{M}$ and $\overline{P}$ of length $t$ that satisfy
$shp(\overline{M}) = shp(M)$, $shp(\overline{P}) = shp(P)$ and the following three conditions:
\begin{itemize}
\item[(i)] $\overline{M}'' = \overline{P}''$,
\item[(ii)] $shp(\overline{M}') = shp(M')$, 
\item[(iii)] $shp(\overline{P}') = shp(P')$.
\end{itemize}
It may appear that these conditions cannot always be satisfied.  
However,  paths $\overline{M}$ and $\overline{P}$ 
satisfying these conditions can always be constructed as follows \cite{rw}.
 
By considering $\overline{M}$ and $\overline{P}$ as standard tableaux,  
we obtain the desired paths $\overline{M}$ and $\overline{P}$ by ensuring the following
is true.
Firstly, fix $t$ to be in the same box in $\overline{M}$ (resp. $\overline{P}$) 
that $t$ is in $M$ (resp. $P$).
Then, fix $(t-1)$ to be in the same box in $\overline{M}$
(resp. $\overline{P}$) that $t$ is in $P$ (resp. $M$).
Lastly, for each $i = 1, 2, \ldots, t-2$, 
fix $i$ to be in the same box in $\overline{M}$ that it is in  $\overline{P}$.

We then define
$$o_{MP} = \frac{1-q^{2d}}{\sqrt{(1-q^{2(d+1)})(1-q^{2(d-1)})}} 
           o_{M'\overline{M}'} g_{t-1} o_{\overline{P}'P'}o_{PP},$$
where $d = d(\overline{M},t-1)$ is as given in (\ref{eq:brontethewhalecatslothdogmouse}).

This completes the definition of the `pre-matrix units'; now we
define the matrix units proper for $BW_{t}$.  

Assume that the matrix units are known for $BW_{t-1}$.
Let $M$ and $P$ be paths of length $t$
in the Bratteli diagram for $BW_{t}$
where $shp(M) = \lambda = shp(P)$ 
and $\lambda$ contains strictly fewer than $t$ boxes, then we define
$$e_{MP} = \frac{Q_{\lambda}(r,q)}{\sqrt{Q_{\mu}(r,q)Q_{\widetilde{\mu}}(r,q)}}
                 e_{M'S} e_{t-1} e_{TP'},$$
where $S$ and $T$ are paths of length $t-1$ satisfying 
\begin{itemize}
\item[(i)]   $shp(S)=shp(M') = \mu$, and
\item[(ii)]  $shp(T)=shp(P') = \widetilde{\mu}$, and
\item[(iii)] $S'=T'$, and
\item[(iv)]  $shp(S')= \lambda = shp(T')$.
\end{itemize}
It may appear that these conditions cannot always be satisfied.  
However, there always exists a pair of paths $S$ and $T$ of length $t-1$
satisfying these conditions for the following reasons.
Firstly, by examining the relevant Bratteli diagrams, it is clear that 
there are no intertwiner matrix units in $BW_{1}$ and $BW_{2}$.  
Now for each $t \geq 3$, a shape $\lambda$ that has at most $t-2$ boxes and which labels a vertex
on the $t^{th}$ level of the Bratteli diagram for $BW_{t}$ 
also labels a vertex on the $(t-2)^{nd}$ level of the Bratteli diagram.
Hence there always exists at least one path of length $t-2$ 
in the Bratteli diagram for $BW_{t}$ ending at the vertex 
$shp(M)$ on the $(t-2)^{nd}$ level, as $shp(M)$ contains no more than $t-2$ boxes 
(note that this shows that (iii) and (iv) might be satisfied).

In the Bratteli diagram for $BW_{t}$,
two vertices $\lambda$ and $\mu$ are connected by an edge 
only if their shapes differ by exactly one box.
Now the vertices $shp(M')$ and $shp(P')$ on the $(t-1)^{st}$ level
are connected to the vertex $shp(M)$ on the $t^{th}$ level by one edge each, and they
are also connected to the vertex $shp(M)$ on the $(t-2)^{nd}$ level by one edge each.
It follows, then, that by fixing $S$ and $T$ to be paths of length $t-1$ that coincide on the first
$t-2$ levels of the Bratteli diagram 
and that pass through the vertex $shp(M)$ on the $(t-2)^{nd}$ level, and also fixing
$shp(S) = shp(M')$ and $shp(T)=shp(P')$ (which is always possible), 
we obtain the desired paths $S$ and $T$.

Let $M$ and $P$ be paths of length $t$
in the Bratteli diagram for $BW_{t}$, where $(M,P) \in \Omega_{t}$, and where
$shp(M)$ contains $t$ boxes. 
Then we define
$$e_{MP} = (1-z_{t})o_{MP},$$
where $z_{t} = \sum_{P} e_{PP}$ with the summation going over all paths $P$ of length $t$ such that
$shp(P)$ contains fewer than $t$ boxes.

The following fact is important  \cite[Lem. 4.2]{w2}: 
let $M$ be a path of length $t$ in the Bratteli diagram for $BW_{t}$
where $shp(M) = \lambda$,
then 
\begin{equation}
\label{eq:gregredinumber3}
\mathrm{tr}(e_{MM}) = Q_{\lambda}(r,q)/x^{t}, 
\end{equation}
where 
$\displaystyle{x = \frac{r-r^{-1}}{q-q^{-1}} + 1}$ and $Q_{\lambda}(r,q)$ is the function
given below in (\ref{eq:themagixQpolynomial}).

It is interesting to note that the quantum superdimension of the fundamental irreducible
$U_{q}(osp(1|2n))$-module $V$ is $\displaystyle{(-q^{2n}+q^{-2n})/(q-q^{-1}) + 1}$, which is just
the expression $x$ 
with the indeterminates $r$ and $q$ replaced with the complex numbers $-q^{2n}$ and $q$, respectively.
Recall that we grade the highest weight vector of $V$ to be odd.

\end{subsection}

\begin{subsection}{The algebra $BW_{t}(r,q)$}

The algebra $BW_{t}(r,q)$, with $r, q \in \mathbb{C}$, is equipped with a functional 
$\mathrm{tr}: BW_{t}(r,q) \rightarrow \mathbb{C}$ which satisfies, 
amongst other relations,
\begin{equation}
\label{eq:tracefunctional}
\mathrm{tr}(a \chi b) = \mathrm{tr}(\chi) \mathrm{tr}(ab), 
\hspace{10mm} \forall a, b \in BW_{t-1}(r,q), \hspace{5mm} \chi \in \{g_{t-1}, e_{t-1}\},
\end{equation}
where we regard each element of $BW_{t-1}(r,q)$ as an element of $BW_{t}(r,q)$
under the canonical inclusion.

Define the annihilator ideal $J_{t}(r,q) \subset BW_{t}(r,q)$ with respect to 
$\mathrm{tr}$ by
$$J_{t}(r,q) = \left\{b \in BW_{t}(r,q) | \ 
\mathrm{tr}(ab) = 0, \ \forall a \in BW_{t}(r,q)\right\}.$$
If $q$ is not a root of unity and $r \neq \pm q^{k}$ for all integers $k$, then
$J_{t}(r,q)=0$ and $BW_{t}(r,q)$ is semisimple  \cite[Cor. 5.6]{w2}.
If $r = \pm q^{k}$ for some $k \in \mathbb{Z}$, then
$J_{t}(\pm q^{k},q) \neq 0$ and
the quotient $BW_{t}(\pm q^{k},q)/J_{t}(\pm q^{k},q)$ is semisimple
\cite[Cor. 5.6]{w2}. 

Let us now fix $k=2n$ and $r=-q^{2n}$;
recall that the homomorphism $\Upsilon: g_{i} \mapsto -\check{R}_{i} \in {\cal{C}}_{t}$ 
yields a representation of $BW_{t}(-q^{2n},q)$ in ${\cal{C}}_{t}$.
The next task is to determine the structure of 
$BW_{t}(-q^{2n},q)/J_{t}(-q^{2n},q)$, which we do 
in the following work from \cite{w2}.

Let us now introduce a subgraph $\Gamma(-q^{2n},q)$ of the Young lattice that we will use
in describing the structure of
$BW_{t}(-q^{2n},q)/J_{t}(-q^{2n},q)$.
We inductively obtain the vertices of $\Gamma(-q^{2n},q)$ as follows.
Firstly fix the Young diagram with no boxes to belong to $\Gamma(-q^{2n},q)$.
The inductive step is that if the Young diagram $\mu$ belongs to $\Gamma(-q^{2n},q)$, 
the Young diagram $\lambda$ also belongs to $\Gamma(-q^{2n},q)$
if $\lambda$ differs from $\mu$ by exactly one box and if 
$Q_{\lambda}(-q^{2n},q) \neq 0$, where
$Q_{\lambda}(r,q)$ is given in the next paragraph.

Given a Young diagram $\lambda$, let
$(i,j)$  denote the box in the $i^{th}$ row and the $j^{th}$ column of $\lambda$, 
and let $\lambda_{i}$ (resp. $\lambda_{j}'$) 
denote the number of boxes in the $i^{th}$ row (resp. $j^{th}$ column) of $\lambda$.
We introduce some useful notation: we may denote the Young diagram $\lambda$ by 
$\lambda = [\lambda_{1}, \lambda_{2}, \ldots, \lambda_{k}]$ where the $i^{th}$ row contains 
$\lambda_{i}$ boxes for each $i=1, 2, \ldots, k$, 
and the $l^{th}$ row contains no boxes for each $l > k$.
The function $Q_{\lambda}(r,q)$ is
\begin{eqnarray}
Q_{\lambda}(r,q) & = & 
\prod_{(j,j) \in \lambda} 
\frac{rq^{\lambda_{j} - \lambda_{j}'}-r^{-1}q^{-\lambda_{j} +\lambda_{j}'}+
q^{\lambda_{j} + \lambda_{j}'-2j+1}-q^{-\lambda_{j} -\lambda_{j}'+2j-1}}
{q^{h(j,j)}-q^{-h(j,j)}} \nonumber \\
& & \hspace{5mm} \times
\prod_{(i,j) \in \lambda, i \neq j} \frac{r q^{d(i,j)}-r^{-1}q^{-d(i,j)}}{q^{h(i,j)}-q^{-h(i,j)}},
\label{eq:themagixQpolynomial}
\end{eqnarray}
where the hooklength $h(i,j)$ is defined by
$h(i,j) = \lambda_{i} - i + \lambda_{j}' - j +1$, and where
$$d(i,j) = \left\{ \begin{array}{ll}
 \lambda_{i} + \lambda_{j}-i-j+1, & \hspace{5mm} \mbox{if } i \leq j, \\
-\lambda_{i}' - \lambda_{j}' + i + j-1, & \hspace{5mm} \mbox{if } i > j.
\end{array} \right.$$
Intuitively, 
the hooklength $h(i,j)$ is the number of boxes 
in the hook going through the box $(i,j)$, ie the number of boxes
below the $(i,j)$ box in the $j^{th}$ column 
plus the number of boxes to the right of the $(i,j)$ box in the $i^{th}$ row, plus one.

Now $h(i,j) \geq 1$ for all $(i,j) \in \lambda$, so $Q_{\lambda}(-q^{2n},q)$ is well-defined for all $\lambda$. Also, for each $(j,j) \in \lambda$ we have
\begin{eqnarray*}
\lefteqn{ -q^{2n + \lambda_{j} - \lambda_{j}'}+q^{-2n-\lambda_{j} +\lambda_{j}'}+
q^{\lambda_{j} + \lambda_{j}'-2j+1}-q^{-\lambda_{j} -\lambda_{j}'+2j-1} } \\
& & \hspace{20mm} = (q^{-n+\lambda_{j}'-j+1/2}-q^{n-\lambda_{j}'+j-1/2})
      (q^{n+\lambda_{j}-j+1/2}+q^{-n-\lambda_{j}+j-1/2}),
\end{eqnarray*}
and so $Q_{\lambda}(-q^{2n},q) = 0$ if and only if at least one of the following conditions is satisfied:
\begin{itemize}
\item[(a)] $q^{4n + 2d(i,j)}=1$ for some $(i,j) \in \lambda$ where $i \neq j$, 
\item[(b)] $q^{2n-2\lambda_{j}'+2j-1} = 1$ or $q^{2n+2\lambda_{j}-2j+1}=-1$ for some $j$.
\end{itemize}
Now $Q_{\lambda}(-q^{2n},q) = Q_{\lambda}\big((-q)^{2n},-q\big)$
from Lemma \ref{lem:JCBsong}, thus 
\newline
$\Gamma(-q^{2n},q) = \Gamma\big((-q)^{2n},-q\big)$. 
Wenzl completely determined $\Gamma\big((-q)^{2n},-q\big)$ in \cite[Cor. 5.6 (c),(c1)]{w2}, 
which we restate here with a correction.
Firstly, \cite[Cor. 5.6 (c)]{w2} reads
\begin{itemize}
\item[(c)] `If $r = q^{n}$ and $q$ is not a root of unity, 
$\Gamma(r, q)$ consists of all Young diagrams $\lambda$ for which'
\end{itemize} 
and \cite[Cor. 5.6 (c1)]{w2} is presented with a slight error; it should read
\begin{itemize}
\item[(c1)] `$\lambda_{1}' + \lambda_{2}' \leq n+1$ for $n>0$.'
\end{itemize}
Then both $\Gamma\big((-q)^{2n},-q\big)$ and $\Gamma(-q^{2n},q)$ consist of all Young diagrams $\lambda$
satisfying $\lambda_{1}' + \lambda_{2}' \leq 2n+1$.
Let us determine $\Gamma(-q^{2n},q)$ independently of \cite[Cor. 5.6 (c),(c1)]{w2}.

Now $q$ is non-zero and not a root of unity, so (b) above is never satisfied for any $\lambda$, and
(a) is only satisfied if $d(i,j) = -2n$.
We now determine the circumstances in which $d(i,j) = -2n$. 
If $i > j$, we can see that
$\mathrm{min}(d(i,j)) = d(2,1) = -\lambda_{1}'-\lambda_{2}'+2$
from the constraints on the lengths of the columns of a Young diagram and
it follows that
$Q_{\lambda}(-q^{2n},q) = 0$ if $\lambda_{1}' + \lambda_{2}' = 2n+2$.
Let us call a Young diagram $\lambda$ {\emph{allowable}} if
$\lambda_{1}' + \lambda_{2}' \leq 2n+1$.

Across all the allowable Young diagrams, let us calculate $\mathrm{min}(d(i,j))$ where $i < j$.
If the first column of the allowable diagram $\lambda$ contains $2n+1$ boxes, 
ie $\lambda_{1}'=2n+1$, then all the other columns must contain no boxes
from the definition of an allowable diagram.
For such a $\lambda$, there does not exist any box $(i,j)$ in the $i^{th}$ row and the
$j^{th}$ column with $i < j$ and so there is nothing more to consider in this case. 
Now if the first column of $\lambda$ contains strictly fewer than $2n+1$ boxes, ie
$\lambda'_{1} \leq 2n$, then the following relations hold:
$i \leq 2n$, $\lambda_{i} - j \geq 0$ and $\lambda_{j} \geq 0$. 
Then $d(i,j)=\lambda_{i} + \lambda_{j}-i-j+1 \geq -2n+1$,
which means that $d(i,j) \neq -2n$ for all $i < j$.

It follows that 
$Q_{\lambda}(-q^{2n},q) = 0$ if $\lambda'_{1} + \lambda'_{2} = 2n+2$ and that
$Q_{\lambda}(-q^{2n},q) \neq 0$ for all allowable Young diagrams $\lambda$.
Consequently, the vertices of $\Gamma(-q^{2n},q)$ are all the allowable Young diagrams, 
that is, all the Young diagrams $\lambda$ satisfying $\lambda_{1}' + \lambda_{2}' \leq 2n+1$, 
and thus $\Gamma(-q^{2n},q)$ is indeed identical to $\Gamma\big((-q)^{2n},-q\big)$. 
In Figure \ref{fig:subgraphgamma} we show the graph
$\Gamma(-q^{2n},q)$ with $n=1$ up to the $4^{th}$ level.
\begin{figure}[hbt]
\label{fig:subgraphgamma}
\begin{center}
\setlength{\unitlength}{0.17mm}
\begin{picture}(350,400) \large\sf
\thinlines
%\linethickness{0.05mm}

\put(150,370){\circle{20}} 
\put(150,325){\line(0,1){30}} 

\put(140,300){\framebox(20,20){}}

\put(90,230){\framebox(20,20){}}
\put(90,250){\framebox(20,20){}}

\put(180,240){\framebox(20,20){}}
\put(200,240){\framebox(20,20){}}

\put(130,170){\framebox(20,20){}}
\put(150,170){\framebox(20,20){}}
\put(130,150){\framebox(20,20){}}

\put(230,160){\framebox(20,20){}}
\put(250,160){\framebox(20,20){}}
\put(270,160){\framebox(20,20){}}

\put(40,180){\framebox(20,20){}}
\put(40,140){\framebox(20,20){}}
\put(40,160){\framebox(20,20){}}

\put(200,90){\framebox(20,20){}}
\put(200,70){\framebox(20,20){}}
\put(220,90){\framebox(20,20){}}
\put(240,90){\framebox(20,20){}}
\put(165,155){\line(1,-1){40}}

\put(290,90){\framebox(20,20){}}
\put(310,90){\framebox(20,20){}}
\put(330,90){\framebox(20,20){}}
\put(350,90){\framebox(20,20){}}
\put(265,155){\line(1,-1){40}} 
\put(215,115){\line(1, 1){40}} 

\put(115,275){\line(1, 1){21}} 

\put(163,297){\line(1,-1){32}} 
\put(115,227){\line(1,-1){32}} 
\put(155,195){\line(1, 1){40}} 
\put(65,205){\line(1, 1){21}} 

\put(205,235){\line(1,-1){50}} 

\end{picture}
\caption{The graph $\Gamma(-q^{2n},q)$ with $n=1$ up to the $4^{th}$ level}
\end{center}
\end{figure}
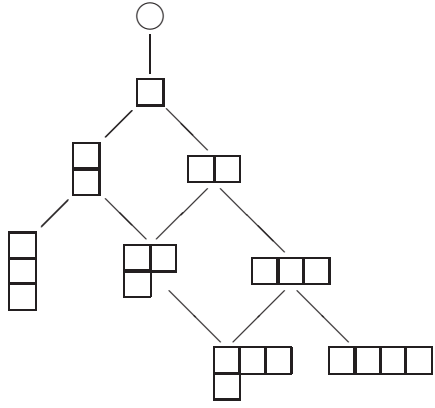

Now $J_{t}(-q^{2n},q) \neq 0$ and $BW_{t}(-q^{2n},q)$ is not semisimple.
However, the quotient $BW_{t}(-q^{2n},q)/J_{t}(-q^{2n},q)$ 
is semisimple:
$$BW_{t}(-q^{2n},q)/J_{t}(-q^{2n},q) \cong
 \bigoplus_{\lambda \in \Gamma_{t}(-q^{2n},q)} M_{b_{\lambda}}(\mathbb{C}),$$
where $\Gamma_{t}(-q^{2n},q)$ is the set of Young diagrams belonging to $\Gamma(-q^{2n},q)$
with $t-2k \geq 0$ boxes, where $k$ ranges over all of $\mathbb{Z}_{+}$  \cite[Cor. 5.6]{w2}.

We obtain the matrix units in 
$BW_{t}(-q^{2n},q)/J_{t}(-q^{2n},q)$ 
by taking a certain proper subset of the matrix units in $BW_{t}$
and replacing the indeterminates $r$ and $q$ with the complex numbers
$-q^{2n}$ and $q$, respectively.

To label the matrix units for $BW_{t}(-q^{2n},q)/J_{t}(-q^{2n},q)$,
we use the Bratteli diagram for $BW_{t}(-q^{2n},q)/J_{t}(-q^{2n},q)$,
which we define in the same way as we defined the Bratteli diagram for
$BW_{t}$ but we replace each $\Gamma_{s}$ with the set $\Gamma_{s}(-q^{2n},q)$ as follows.
Recall that $\Gamma(-q^{2n},q)$ is a subgraph of the Young lattice the vertices of which are all
the allowable Young diagrams.
Then, for each $s=0, 1, \ldots, t$, we fix $\Gamma_{s}(-q^{2n},q)$ 
to be the set of all Young diagrams that are vertices of
$\Gamma(-q^{2n},q)$ that contain exactly $s-2k \geq 0$ boxes, where $k$ ranges over all of 
$\mathbb{Z}_{+}$. 
In Figure \ref{fig:brattelidiagramforquotientalgebra}
we show the Bratteli diagram for $BW_{t}(-q^{2n},q)/J_{t}(-q^{2n},q)$ with
$n=1$ up to the $4^{th}$ level.
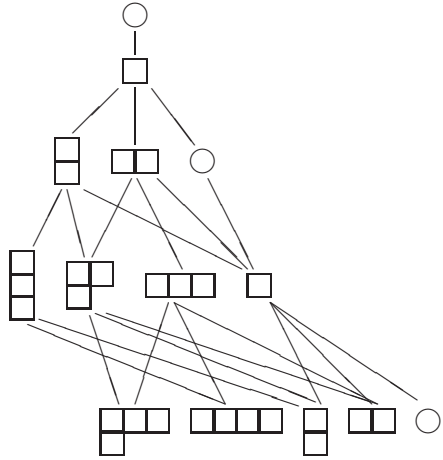
\begin{figure}[hbt]
\label{fig:brattelidiagramforquotientalgebra}
\begin{center}
\setlength{\unitlength}{0.15mm}
\begin{picture}(460,440) \large\sf
\thinlines
%\linethickness{0.05mm}

\put(190,430){\circle{20}} 

\put(190,395){\line(0,1){20}}

\put(180,370){\framebox(20,20){}}

\put(120,300){\framebox(20,20){}}
\put(120,280){\framebox(20,20){}}

\put(170,290){\framebox(20,20){}}
\put(190,290){\framebox(20,20){}}

\put(250,300){\circle{20}}

\put(80,160){\framebox(20,20){}}
\put(80,180){\framebox(20,20){}}
\put(80,200){\framebox(20,20){}}
\put(100,225){\line(1,2){25}} 
\put(130,275){\line(1,-4){15}}

\put(130,170){\framebox(20,20){}}
\put(130,190){\framebox(20,20){}}
\put(150,190){\framebox(20,20){}}
\put(152,215){\line(1,2){35}} 
\put(192,285){\line(1,-2){40}}

\put(200,180){\framebox(20,20){}}
\put(220,180){\framebox(20,20){}}
\put(240,180){\framebox(20,20){}}

\put(290,180){\framebox(20,20){}}
\put(210,285){\line(1,-1){80}} 
\put(255,285){\line(1,-2){40}} 
\put(145,275){\line(2,-1){140}}

\put(135,325){\line(1,1){40}} 
\put(190,315){\line(0,1){50}} 
\put(205,365){\line(3,-4){38}} 

\put(160,60){\framebox(20,20){}}
\put(160,40){\framebox(20,20){}}
\put(180,60){\framebox(20,20){}}
\put(200,60){\framebox(20,20){}}
\put(150,163){\line(1,-3){26}} 
\put(190,85){\line(1,3){30}} 
\put(105,160){\line(3,-1){230}}

\put(240,60){\framebox(20,20){}}
\put(260,60){\framebox(20,20){}}
\put(280,60){\framebox(20,20){}}
\put(300,60){\framebox(20,20){}}
\put(225,175){\line(1,-2){45}} 
\put(95,155){\line(5,-2){175}} 

\put(340,60){\framebox(20,20){}}
\put(340,40){\framebox(20,20){}}
\put(155,165){\line(5,-2){195}} 
\put(310,175){\line(1,-1){90}} 
\put(310,175){\line(1,-2){45}} 
\put(310,175){\line(3,-2){130}} 

\put(380,60){\framebox(20,20){}}
\put(400,60){\framebox(20,20){}}
\put(165,165){\line(3,-1){240}} 
\put(225,175){\line(2,-1){180}} 

\put(450,70){\circle{20}} 

\end{picture}
\caption{The Bratteli diagram for $BW_{t}(-q^{2n},q)/J_{t}(-q^{2n},q)$ with
$n=1$ up to the $4^{th}$ level}  
\end{center}
\end{figure}

We say that $T=(0, s_{1}, \ldots, s_{t})$
is a path of length $t$ in the Bratteli diagram for
$BW_{t}(-q^{2n},q)/J_{t}(-q^{2n},q)$ 
 if $s_{i} \in \Gamma_{i}(-q^{2n},q)$ for each $i$ and if 
 $s_{j}$ is joined by an edge to $s_{j+1}$ for each $j=0, 1, \ldots, t-1$.
Note that $s_{j}$ is joined to $s_{j+1}$ by an edge 
if and only if $s_{j}$ differs from $s_{j+1}$ by exactly one box.
 
Let $\Omega_{t}(-q^{2n},q)$ be the set of pairs $(R,S)$ 
of paths of length $t$ in the Bratteli diagram for
$BW_{t}(-q^{2n},q)/J_{t}(-q^{2n},q)$ 
where $r_{t} = s_{t}$, that is $shp(R)=shp(S)$.
The matrix units 
$$\{ e_{RS} \in BW_{t} | \ (R,S) \in \Omega_{t}(-q^{2n},q) \}$$ 
are all well-defined and non-zero upon the indeterminates
$r$ and $q$ being replaced with the complex numbers $-q^{2n}$ and $q$, respectively.
Henceforth we write $e_{RS}$ to mean the matrix unit 
$e_{RS}(-q^{2n},q) \in BW_{t}(-q^{2n},q)/J_{t}(-q^{2n},q)$.
It is very important to note that $\mathrm{tr}(e_{SS}) \neq 0$ 
for all $(S,S) \in \Omega_{t}(-q^{2n},q)$ and that 
$e_{RS} \notin J_{t}(-q^{2n},q)$ for all  $(R,S) \in \Omega_{t}(-q^{2n},q)$.

\end{subsection}

\begin{subsection}{Matrix units in $BW_{t}(-q^{2n},q)/J_{t}(-q^{2n},q)$
and ${\cal{C}}_{t}$}

We now relate the idempotent matrix units in 
$BW_{t}(-q^{2n},q)/J_{t}(-q^{2n},q)$
to the projectors in $\mathcal{C}_{t}$ we defined at the start of this section.
Let $\widetilde{BW}_{t}(-q^{2n},q)$ be the semisimple subalgebra of
$BW_{t}(-q^{2n},q)$ spanned by 
the matrix units in $BW_{t}(-q^{2n},q)/J_{t}(-q^{2n},q)$, ie
$\{ e_{RS} | \ (R,S) \in \Omega_{t}(-q^{2n},q) \}$.  

Firstly, we will show that
$BW_{t}(-q^{2n},q) = 
\widetilde{BW}_{t}(-q^{2n},q) \oplus J_{t}(-q^{2n},q)$.
Any $f \in \widetilde{BW}_{t}(-q^{2n},q)$ can be written as
$$f = \sum_{(S,T) \in \Omega_{t}(-q^{2n},q)} f_{ST} e_{ST}, \hspace{10mm} f_{ST} \in \mathbb{C},$$
where $f_{ST} \neq 0$ for at least one pair $(S,T)$ of paths.  
Fix $(A,B)$ to be such a pair, then
$$\mathrm{tr}(e_{BA}f) = \mathrm{tr}(f_{AB} e_{BA} e_{AB}) = f_{AB}\mathrm{tr}(e_{BB}) \neq 0,$$
as $\mathrm{tr}(e_{BB}) \neq 0$.  
Thus any non-zero $f$ belonging to 
$\widetilde{BW}_{t}(-q^{2n},q)$ does not belong to 
$J_{t}(-q^{2n},q)$, yielding
$$BW_{t}(-q^{2n},q) = 
\widetilde{BW}_{t}(-q^{2n},q) \oplus J_{t}(-q^{2n},q).$$
Then we can write $a = \widetilde{a} + a_{j}$ for each $a \in BW_{t}(-q^{2n},q)$, 
  where $\widetilde{a} \in \widetilde{BW}_{t}(-q^{2n},q)$ and 
  $a_{j} \in J_{t}(-q^{2n},q)$.

Now define 
$$P_{t} = \sum_{(S,S) \in \Omega_{t}(-q^{2n},q)} e_{SS} \in 
\widetilde{BW}_{t}(-q^{2n},q),$$
then $P_{t} a P_{t} = \widetilde{a}$, which can be seen by regarding
$BW_{t}(-q^{2n},q)$ as a matrix algebra.

\end{subsection}

Let us now turn our attention to ${\cal{C}}_{t}$.
Define ${\cal{J}}_{t} \subset {\cal{C}}_{t}$ to be the annihilator ideal of 
${\cal{C}}_{t}$ with respect to the quantum supertrace:
$${\cal{J}}_{t} = \{ b \in {\cal{C}}_{t} | \ \mathrm{str}_{q}(ab)=0, \ \forall a \in {\cal{C}}_{t} \}.$$
Now define a map $\psi: {\cal{C}}_{t} \rightarrow \mathbb{C}$ by
$$\psi(X) = \mathrm{str}_{q}(X)/ \big( sdim_{q}(V) \big)^{t},$$
then $\psi(X) = 0$ if and only if $\mathrm{str}_{q}(X) = 0$, and furthermore,
\begin{equation}
\label{eq:thegreaans(1)}
\psi \big(\Upsilon(a) \big) = \mathrm{tr}(a), \hspace{10mm} \forall a \in BW_{t}(-q^{2n},q),
\end{equation}
from Lemma \ref{lemma:equalityoftraces}.
Thus we can regard ${\cal{J}}_{t}$ as the annihilator ideal of ${\cal{C}}_{t}$ with respect to $\psi$.

\begin{remark}
The equality of the traces in (\ref{eq:thegreaans(1)}) independently confirms that link invariants can be created from representations of $U_{q}(osp(1|2n))$ as Zhang did in \cite{z1}.  
\end{remark}

Now we will use Eq. (\ref{eq:thegreaans(1)}) to show that
\begin{equation}
\label{eq:shonauy(99)}
\Upsilon\big(J_{t}(-q^{2n},q)\big) = {\cal{J}}_{t}.
\end{equation}
We firstly show that $\Upsilon\big(J_{t}(-q^{2n},q)\big) \subseteq {\cal{J}}_{t}$.
Let $b$ be an arbitrary element of $J_{t}(-q^{2n},q)$, then $\mathrm{tr}(ab)=0$ for all 
$a \in BW_{t}(-q^{2n},q)$, and the surjectivity of 
$\Upsilon$, in addition to the fact that $\psi\big(\Upsilon(ab)\big) = \mathrm{tr}(ab)$, means that
$\Upsilon(b) \in {\cal{J}}_{t}$.

Now let $B$ be an arbitrary element of ${\cal{J}}_{t}$, 
then there exists some $b \in BW_{t}(-q^{2n},q)$
satisfying  $B = \Upsilon(b)$. Furthermore, $b \in J_{t}(-q^{2n},q)$ as
$\mathrm{tr}(ab) = \psi\big(\Upsilon(a)\Upsilon(b) \big) = 0$ for all 
$a \in BW_{t}(-q^{2n},q)$. Then
$\Upsilon\big(J_{t}(-q^{2n},q)\big) \supseteq {\cal{J}}_{t}$, proving
 (\ref{eq:shonauy(99)}).

The surjectivity of $\Upsilon$ implies that 
$${\cal{C}}_{t} = \Upsilon\big(\widetilde{BW}_{t}(-q^{2n},q) \big) + {\cal{J}}_{t},$$
and we will show that this sum is in fact direct.
To see this, assume that there exists some non-zero element $F$ of ${\cal{C}}_{t}$ belonging to
$\Upsilon\big(\widetilde{BW}_{t}(-q^{2n},q) \big)$ and also belonging to ${\cal{J}}_{t}$, then
$\mathrm{str}_{q}(XF) = 0$ for all $X \in {\cal{C}}_{t}$.  
However, $F$ is the image of a linear combination of matrix units: 
$$F = \sum_{(S,T) \in \Omega_{t}(-q^{2n},q)}f_{ST}\Upsilon(e_{ST}),
\hspace{10mm} f_{ST} \in \mathbb{C},$$ 
where $f_{ST} \neq 0$  for at least one pair $(S,T)$.  
Assume that $(A,B)$ is such a pair, then by similar reasoning as
previously, $\mathrm{str}_{q}\left(\Upsilon(e_{BA})F\right) \neq 0$ 
which contradicts the assumption that 
$F \in {\cal{J}}_{t}$.
Thus
$\Upsilon\big(\widetilde{BW}_{t}(-q^{2n},q) \big) \cap {\cal{J}}_{t} = 0$,
and 
\begin{equation} 
\label{eq:milindaandmilinda}
{\cal{C}}_{t} = \Upsilon\big(\widetilde{BW}_{t}(-q^{2n},q) \big) \oplus {\cal{J}}_{t}.
\end{equation}

It is clear that the image of the collection of all
matrix units $e_{ST} \in \widetilde{BW}_{t}(-q^{2n},q)$ in ${\cal{C}}_{t}$
under the map $\Upsilon$ is again a collection of matrix units.
Each matrix unit $e_{SS} \in \widetilde{BW}_{t}(-q^{2n},q)$ is an idempotent, thus
each $\Upsilon(e_{SS})$ is an idempotent that is also a $U_{q}(\mathfrak{g})$-linear map.
The idempotents 
$\{\Upsilon(e_{SS}) | \ (S,S) \in \Omega_{t}(-q^{2n},q)\}$ are all orthogonal
as the matrix units $\{e_{SS} | \ (S,S) \in \Omega_{t}(-q^{2n},q)\}$ are all orthogonal.  

Let $S$ be a path of length $t$ in the Bratteli diagram for $\widetilde{BW}_{t}(-q^{2n},q)$
and let $T$ be a path of length $t$ in the Bratteli diagram for $V^{\otimes t}$.
Both $\Upsilon(e_{SS})$ and the projection $E_{TT} \in {\cal{C}}_{t}$ 
project down onto the same decomposition of $V^{\otimes t}$.
To see this, let us write $S = (0, s_{1}, \ldots, s_{t})$ and $T = (0, t_{1}, \ldots, t_{t})$
and fix $S_{j} = (0, s_{1}, \ldots, s_{j})$ and $T_{j} = (0, t_{1}, \ldots, t_{j})$ for each
$j=1, 2, \ldots, t$. For each $j$, $\Upsilon(e_{S_{j}S_{j}}) \in {\cal{C}}_{j}$ and
$E_{T_{j}T_{j}} \in {\cal{C}}_{j}$ are embedded into ${\cal{C}}_{t}$ via
$Z \mapsto Z \otimes \mathrm{id}^{\otimes (t-j)}$, and
they are recursively defined so that
$\Upsilon(e_{S_{j-1}S_{j-1}})$ and $E_{T_{j-1}T_{j-1}}$ explicitly appear in
the definitions of 
$\Upsilon(e_{S_{j}S_{j}})$ and $E_{T_{j}T_{j}}$, respectively.
Lastly, note that both $\Upsilon(e_{S_{j}S_{j}})$ and $E_{T_{j}T_{j}}$ project down onto
the $j$ left-most tensor powers of $V$ in $V^{\otimes t}$.

Fix $e_{SS} \in \widetilde{BW}_{t}(-q^{2n},q)$ to be an idempotent matrix unit where
$(S,S) \in \Omega_{t}(-q^{2n},q)$.
Now $\Upsilon(e_{SS}) \big( V^{\otimes t} \big) \neq 0$ as 
$\mathrm{str}_{q}\big(\Upsilon(e_{SS})\big) = (sdim_{q}(V))^{t}\mathrm{tr}(e_{SS}) \neq 0$, and
as $V^{\otimes t}$ is completely reducible, $\Upsilon(e_{SS})$ projects down from $V^{\otimes t}$
onto a direct sum of irreducible $U_{q}(\mathfrak{g})$-submodules $W_{shp(S)}$ of $V^{\otimes t}$.

Let us write $\lambda = shp(S)$ and define the Young diagram $\lambda^{*}$ 
corresponding to $\lambda$ by
\begin{equation}
\label{eq:lambdastar}
\lambda^{*} = \left\{ \begin{array}{ll}
\lambda, & \mbox{if } \lambda_{1}' \leq n, \\
\widetilde{\lambda}, & \mbox{if } \lambda_{1}' \geq n+1,
\end{array}
\right.
\end{equation}
where $\widetilde{\lambda}$ is given in Lemma \ref{lem:reflectedyoungdiagram}.
We will inductively show that:
\begin{lemma}
\label{lem:highestweightofWlambda}
$W_{shp(S)}$ is an irreducible 
$U_{q}(\mathfrak{g})$-submodule of $V^{\otimes t}$ with integral dominant highest weight 
$\lambda^{*}$, where
the Young diagram $\lambda^{*}$ is taken to be a highest weight
as discussed in the third paragraph of 
Subsection \ref{subsect:matrixbirmanwenzlstuff}.
\end{lemma}
\begin{proof}
Firstly, let $R=(0, \epsilon_{1})$ be a path of length $1$; $e_{RR} = 1$ and 
$\Upsilon(e_{RR})$ acts as the identity on $V$, the fundamental
irreducible $U_{q}(\mathfrak{g})$-module with integral dominant highest weight $\epsilon_{1}$.

We now do the inductive step.
Assume that $R$ is a path of length $t-1$
where $(R, R) \in \Omega_{t-1}(-q^{2n},q)$ and such that the idempotent
$\Upsilon(e_{RR}) \in {\cal{C}}_{t-1}$ projects down from $V^{\otimes (t-1)}$ onto an irreducible
$U_{q}(\mathfrak{g})$-submodule $V_{\mu^{*}} \subseteq V^{\otimes (t-1)}$
with integral dominant highest weight $\mu^{*}$, where
$\mu^{*}$ is the weight corresponding to $\mu = shp(R)$ as given in (\ref{eq:lambdastar}).
Let $S$ be a path of length $t$ where $(S, S) \in \Omega_{t}(-q^{2n},q)$ and 
such that $S' = R$ and let $\lambda = shp(S)$, then
$\Upsilon(e_{SS}) \in {\cal{C}}_{t}$ projects down from $V^{\otimes t}$ onto exactly one irreducible
$U_{q}(\mathfrak{g})$-submodule $W_{\lambda}$ of $V_{\mu^{*}} \otimes V$ from 
Lemma \ref{lem:projectionprojectsontoasingleirreducible}. 
We now show that the highest weight of
$W_{\lambda}$ is the integral dominant weight $\lambda^{*}$.

From (\ref{eq:gregredinumber3}) and (\ref{eq:thegreaans(1)}), 
$sdim_{q}(W_{\lambda}) = Q_{\lambda}(-q^{2n},q)$, and from
Lemma \ref{lem:JCBsong}, $Q_{\lambda}(-q^{2n},q) = Q_{\lambda}(q^{2n},-q)$.
If $\lambda_{1}' \geq n+1$, then
$Q_{\widetilde{\lambda}}(q^{2n},-q) = Q_{\lambda}(q^{2n},-q)$ from \cite[p. 422]{w2}, and it
follows that $Q_{\lambda^{*}}(q^{2n},-q) = Q_{\lambda}(q^{2n},-q)$.

It is well known that $Q_{\lambda^{*}}(q^{2n},-q)$ is the quantum dimension of
a finite dimensional irreducible $U_{-q}(so(2n+1))$-module $V_{\lambda^{*}}^{so(2n+1)}$
with highest weight $\lambda^{*}$:
$$dim_{-q}\left(V^{so(2n+1)}_{\lambda^{*}}\right) = Q_{\lambda^{*}}(q^{2n},-q).$$
From \cite{zsuper}, there exists a finite dimensional irreducible 
$U_{q}(osp(1|2n))$-module $V_{\lambda^{*}}$ with highest weight $\lambda^{*}$
with the same weight space decomposition as $V_{\lambda^{*}}^{so(2n+1)}$. 
The quantum superdimension of $V_{\lambda^{*}}$ from 
\cite[Eq. (15)]{zsuper} satisfies
$$sdim_{q}(V_{\lambda^{*}}) = dim_{-q}\left(V^{so(2n+1)}_{\lambda^{*}}\right),$$
and thus
\begin{equation}
\label{eq:equalityofsuperdimensions}
sdim_{q}(W_{\lambda}) = Q_{\lambda^{*}}(-q^{2n},q) = 
Q_{\lambda^{*}}(q^{2n},-q) = sdim_{q}(V_{\lambda^{*}}).
\end{equation}

Note that we have not yet proved that the highest weight of the irreducible
$U_{q}(\mathfrak{g})$-module $W_{\lambda}$ is $\lambda^{*}$. 
Recall that $V_{\mu^{*}} \otimes V$ is completely reducible; 
from Lemma \ref{lem:differentquantumsuperdimensions}, 
the quantum superdimensions of the irreducible $U_{q}(\mathfrak{g})$-submodules of 
$V_{\mu^{*}} \otimes V$ are all different, thus
the highest weight of $W_{\lambda}$ is 
indeed $\lambda^{*}$ from (\ref{eq:equalityofsuperdimensions}).
This completes the proof.
\end{proof}

If $R$ and $S$ are paths of length $t$ in the Bratteli diagram for
$\widetilde{BW}_{t}(-q^{2n},q)$ satisfying
$shp(R)=\lambda$ and $shp(S)=\widetilde{\lambda}$,
where $\widetilde{\lambda}$ is given in Lemma \ref{lem:reflectedyoungdiagram}, then
$\Upsilon(e_{RR})$ and $\Upsilon(e_{SS})$ project down from $V^{\otimes t}$ onto 
isomorphic irreducible $U_{q}(\mathfrak{g})$-submodules of $V^{\otimes t}$.
However, we show in the next paragraph that no such paths exist in the Bratteli diagram for
$\widetilde{BW}_{t}(-q^{2n},q)$ using an easy even/odd counting argument.

If the number of boxes in $\lambda$ is even (resp. odd), then the number of boxes in
$\widetilde{\lambda}$ is odd (resp. even), as 
$$\widetilde{\lambda}_{1}' \bmod{2} = (2n+1-\lambda_{1}') \bmod{2} = (1-\lambda_{1}') \bmod{2},
\hspace{5mm} \mbox{and} \hspace{5mm} \widetilde{\lambda}_{j}' = \lambda_{j}',  \ \ j \geq 2.$$
Now let $r$ be an even (resp. odd) number satisfying $0 \leq r \leq t$,
then the vertices on the $r^{th}$ level of the Bratteli diagram for
$\widetilde{BW}_{t}(-q^{2n},q)$
are all the Young diagrams in  $\Gamma(-q^{2n},q)$ with  $k \geq 0$ boxes where $k \leq r$ is
an even (resp. odd) number. 
If $|\lambda|$ denotes the number of boxes in the Young diagram $\lambda$, then the fact that
 $|\lambda| \bmod{2} = \big( |\widetilde{\lambda}|+1 \big) \bmod{2}$ means that
$\lambda$ and $\widetilde{\lambda}$ cannot both be vertices on the same level of 
the Bratteli diagram for $\widetilde{BW}_{t}(-q^{2n},q)$, and it follows that at most only one
of the paths $R$ and $S$ can exist.

Let $R$ be a path of length $t$ in the Bratteli diagram for $\widetilde{BW}_{t}(-q^{2n},q)$.
Let $S$ be the same path of length $t$ as $R$ except that we replace each 
Young diagram $\lambda$ on the path $R$ 
with more than $n$ rows of boxes with the Young diagram 
$\widetilde{\lambda}$ given in (\ref{eq:ACER99}).
Then $S$ is a path of length $t$ in the Bratteli diagram for $V^{\otimes t}$,
and we have $E_{SS} = \Upsilon(e_{RR})$, ie the projectors $E_{SS}$ 
and $\Upsilon(e_{RR})$ project from $V^{\otimes t}$ down onto the same $U_{q}(\mathfrak{g})$-submodule.

We now present some technical lemmas we used in this work.

\begin{lemma}
\label{lem:projectionprojectsontoasingleirreducible}
Let $R$ be a path of length $t-1$ where $(R, R) \in \Omega_{t-1}(-q^{2n},q)$ 
such that the idempotent
$\Upsilon(e_{RR})$ projects down from $V^{\otimes (t-1)}$ onto an irreducible
$U_{q}(\mathfrak{g})$-submodule $V_{\mu} \subseteq V^{\otimes (t-1)}$
with integral dominant highest weight $\mu$.
Let $S$ be a path of length $t$ where $(S, S) \in \Omega_{t}(-q^{2n},q)$ 
such that $S' = R$, then
$\Upsilon(e_{SS})$ projects down from $V^{\otimes t}$ onto exactly one irreducible
$U_{q}(\mathfrak{g})$-submodule of $V_{\mu} \otimes V$.
\end{lemma}
\begin{proof}
Recall that $V_{\mu} \otimes V$ is completely reducible.
Assume that $\Upsilon(e_{SS})$ projects down onto at least two irreducible
$U_{q}(\mathfrak{g})$-submodules of $V_{\mu} \otimes V$:
$$\Upsilon(e_{SS}) (V_{\mu} \otimes V) = 
\bigoplus_{\nu \in I} V_{\nu},$$
for some index set $I$. Then, by construction, there exists a projection
$E_{MM_{\nu}} \in {\cal{C}}_{t}$ for each $\nu \in I$ such that
$E_{MM_{\nu}}(V_{\mu} \otimes V) = V_{\nu} \subseteq V_{\mu} \otimes V$.

From (\ref{eq:milindaandmilinda}), each $E_{MM_{\nu}}$ can be written as a sum
$E_{MM_{\nu}} = \widetilde{E}_{MM_{\nu}} + E_{MM_{\nu}}^{j}$
where 
$$\widetilde{E}_{MM_{\nu}} = \sum_{(S,T) \in \Omega_{t}(-q^{2n},q)} c_{ST} \Upsilon(e_{ST})
\in \Upsilon\big(\widetilde{BW}_{t}(-q^{2n},q)\big),$$
and $E_{MM_{\nu}}^{j} \in \Upsilon\big(J_{t}(-q^{2n},q)\big)$.
Note that $\widetilde{E}_{MM_{\nu}} \neq 0$ as
$\mathrm{str}_{q}(E_{MM_{\nu}}) \neq 0$ and $\mathrm{str}_{q}\big(E_{MM_{\nu}}^{j}\big)=0$.
From Lemma  \ref{lem:theredoesnotexist}, for each projection
$E_{SS} \in {\cal{C}}_{t}$ there exists at least one idempotent
$e_{RR} \in \widetilde{BW}_{t}(-q^{2n},q)$ such that
$$\big(E_{SS} V^{\otimes t}\big) \cap \big(\Upsilon(e_{RR})V^{\otimes t}\big) \neq 0.$$

Let $E_{SS}(V_{\mu} \otimes V) = V_{\nu} \subseteq V^{\otimes t}$
and fix $e_{RR}$ to be such an idempotent, then
$$\big(E_{SS}V^{\otimes t}\big) \cap \big(\Upsilon(e_{RR}) V^{\otimes t}\big) = 
V_{\mu},$$
and $E_{SS} = \Upsilon(e_{RR}) + E_{SS}^{j}$ where 
$E_{SS}^{j} \in \Upsilon\big(J_{t}(-q^{2n},q)\big)$.

Now assume that there exists some non-zero
projection $E_{TT} \in {\cal{C}}_{t}$ orthogonal to 
$E_{SS}$ and satisfying 
$$\big(E_{TT}V^{\otimes t}\big) \cap \big(\Upsilon(e_{RR}) V^{\otimes t}\big) \neq 0,$$
then $E_{TT} = \Upsilon(e_{RR}) + E^{j}_{TT}$ where 
$E^{j}_{TT} \in \Upsilon\big(J_{t}(-q^{2n},q)\big)$.
Now 
$$E_{TT}E_{SS} = \Upsilon(e_{RR}) + \Upsilon(e_{RR})E^{j}_{SS} + 
E^{j}_{TT}\Upsilon(e_{RR}) + E^{j}_{TT}E^{j}_{SS} \neq 0,$$
which is non-vanishing as the first term in the sum is an element of
$\Upsilon \big(\widetilde{BW}_{t}(-q^{2n},q) \big)$
and the last three terms in the sum are elements of
$\Upsilon\big(J_{t}(-q^{2n},q)\big)$.
However, the fact that $E_{TT}E_{SS} \neq 0$
contradicts the assumption that $E_{TT}$ and $E_{SS}$ are orthogonal, 
and thus such a non-zero $E_{TT}$ cannot exist. 
However, such a non-zero $E_{TT}$ must exist if 
$\Upsilon(e_{RR})$ projects onto the direct sum of at least two irreducible modules, thus 
$\Upsilon(e_{RR})$ projects only onto the zero vector or onto one irreducible
$U_{q}(\mathfrak{g})$-module. 
As $\mathrm{str}_{q}(\Upsilon(e_{RR})) \neq 0$, $\Upsilon(e_{RR})$ 
projects onto one irreducible $U_{q}(\mathfrak{g})$-module,
completing the proof.
\end{proof}

\begin{lemma}
\label{lem:differentquantumsuperdimensions}
Let $V_{\mu}$ be a finite-dimensional irreducible $U_{q}(\mathfrak{g})$-module with
integral dominant highest weight $\mu \in {\cal{P}}^{+}$ and let
$$V_{\mu} \otimes V \cong \bigoplus_{\lambda \in {\cal{P}}^{+}_{\mu} } V_{\lambda}$$
be the decomposition of $V_{\mu} \otimes V$ into irreducible 
$U_{q}(\mathfrak{g})$-submodules $V_{\lambda}$.
Then, for all $\lambda, \nu \in {\cal{P}}^{+}_{\mu}$, 
$sdim_{q}\left(V_{\lambda}\right) \neq sdim_{q}\left(V_{\nu}\right)$ 
if $\lambda \neq \nu$.
\end{lemma}
\begin{proof}
Recall that, by definition, 
$$sdim_{q}(V_{\lambda}) = \mathrm{str}(\pi_{\lambda}(K_{2 \rho})) = 
             \sum_{\phi} (-1)^{[w_{\phi}]} \mathrm{mul}(\phi) q^{(2\rho, \phi)},$$
where the sum is over all weights $\phi$ of $V_{\lambda}$, 
$\mathrm{mul}(\phi)$ is the multiplicity of the weight space $\phi$ in $V_{\lambda}$ and $w_{\phi} \in V_{\lambda}$ is a weight vector of weight $\phi$. 
By inspection, $(2\rho, \phi) \leq (2\rho, \lambda)$ for all weights $\phi$ of $V_{\lambda}$.
To complete the proof, all we need show is that $(2\rho,\lambda) \neq (2\rho,\nu)$
for all integral dominant weights 
$\lambda, \nu \in {\cal{P}}^{+}_{\mu}$ where $\lambda \neq \nu$.

Recall that ${\cal{P}}^{+}_{\mu} \subseteq {\cal{P}}^{0}_{\mu} = 
\{\mu, \mu \pm \epsilon_{i} \in {\cal{P}}^{+} | \ i=1, 2, \ldots n \}$.
It is then easy to show that $(2\rho, \beta) \neq (2\rho, \gamma)$ for all
$\beta, \gamma \in {\cal{P}}^{0}_{\mu}$ where $\beta \neq \gamma$, 
which completes the proof.

\end{proof}

\begin{lemma}
\label{lem:reflectedyoungdiagram}
Let $\lambda$ be a Young diagram satisfying $\lambda_{1}' \geq n+1$ and
$\lambda_{1}' + \lambda_{2}' \leq 2n+1$.
Define a diagram as a Young diagram is defined except that the length of the 
$(j+1)^{st}$ row (resp. column) of the diagram can be greater than the length of the 
$j^{th}$ row (resp. column), for each $j$.  
Define the diagram $\widetilde{\lambda}$ corresponding to $\lambda$ by
\begin{equation}
\label{eq:ACER99}
\widetilde{\lambda}_{1}' = 2n+1-\lambda_{1}', \ \mbox{ and } \ \widetilde{\lambda}_{j}' = \lambda_{j}' \ \ 
\mbox{ for all } \ \ j \geq 2.
\end{equation}
Then $\widetilde{\lambda}$ is a Young diagram.
\end{lemma}
\begin{proof}
Let $\lambda$ be a Young diagram as given in the lemma and
write $\lambda_{1}' = \lambda_{2}' + k$ where $k \geq 1$, then
$\widetilde{\lambda}_{1}' = 2n+1 - (\lambda_{2}' + k)$.  
Now $\widetilde{\lambda}$ is a Young diagram if
$\widetilde{\lambda}_{1}' \geq \widetilde{\lambda}_{2}'$, and this condition is just 
$2n+1 - k \geq 2\lambda_{2}'$, which is true as 
$\lambda_{1}' + \lambda_{2}' = 2\lambda_{2}' + k \leq 2n+1$.
This completes the proof.
\end{proof}

\begin{lemma}
\label{lem:JCBsong}
For each Young diagram $\lambda$, 
\begin{equation}
\label{eq:hideandseek}
Q_{\lambda}(-q^{2n},q) = Q_{\lambda}(q^{2n},-q).
\end{equation}
\end{lemma}
\begin{proof}
Simple calculations show that (\ref{eq:hideandseek}) is true if and only if
\begin{eqnarray*}
 \prod_{\stackrel{(i,j) \in \lambda}{i \neq j}} -q^{2n+d(i,j)}+q^{-2n-d(i,j)}
& = & 
\left(\prod_{\stackrel{(i,j) \in \lambda}{i < j}} (-1)^{\lambda_{j}'+\lambda_{j}}(q^{2n+d(i,j)}-q^{-2n-d(i,j)})\right) \\ 
& & \times \left(\prod_{\stackrel{(i,j) \in \lambda}{i > j}} (-1)^{\lambda_{i}+\lambda_{i}'}(q^{2n+d(i,j)}-q^{-2n-d(i,j)})\right)
\end{eqnarray*}
and this last equation is true if
\begin{equation}
\label{eq:makethingsright}
\prod_{\stackrel{(i,j) \in \lambda}{i \neq j}}(-1) = 
\left( \prod_{\stackrel{(i,j) \in \lambda}{i<j}} (-1)^{\lambda_{j}' + \lambda_{j} }  	\right) 
\left( \prod_{\stackrel{(i,j) \in \lambda}{i>j}} (-1)^{\lambda_{i} + \lambda_{i}'}   \right).
\end{equation}
We now show that (\ref{eq:makethingsright}) is true.
Define the following sets:
\begin{eqnarray*}
\mathrm{Hor}_{k} & = & \left\{ (k,j) \in \lambda | \ j=1, 2, \ldots, \min{ \{k-1, \lambda_{k}   \} } \right\} \\
\mathrm{Ver}_{k} & = & \left\{ (i,k) \in \lambda | \ i=1, 2, \ldots, \min{ \{ k-1, \lambda_{k}' \} } \right\}.
\end{eqnarray*}
Noting that $|\mathrm{Ver}_{k} \cap \mathrm{Hor}_{l}| = 0$ for all $k$ and $l$ and that
$|\mathrm{Hor}_{k} \cap \mathrm{Hor}_{l}| = 0 = |\mathrm{Ver}_{k} \cap \mathrm{Ver}_{i}|$ for all $k \neq i$,  
(\ref{eq:makethingsright}) is true if the following equation is true for each $k$: 
\begin{equation}
\label{eq:politesociety}
(-1)^{|\mathrm{Hor}_{k} \cup \mathrm{Ver}_{k}|} = 
\left( \prod_{(i,k) \in \mathrm{Ver}_{k}} (-1)^{\lambda_{k}' + \lambda_{k} }  	\right) 
\left( \prod_{(k,j) \in \mathrm{Hor}_{k}} (-1)^{\lambda_{k} + \lambda_{k}'}   \right).
\end{equation}
%The left hand side of (\ref{eq:politesociety}) equals $1$ (resp. $-1$) if
%$|\mathrm{Hor}_{k} \cup \mathrm{Ver}_{k}|$ is even (resp. odd).
If $|\mathrm{Hor}_{k} \cup \mathrm{Ver}_{k}|$ is even, 
the right hand side of (\ref{eq:politesociety}) clearly equals $1$
as $\mathrm{Hor}_{k}$ and $\mathrm{Ver}_{k}$ are disjoint.
Alternatively, if $|\mathrm{Hor}_{k} \cup \mathrm{Ver}_{k}|$ is odd, 
then $\lambda_{k} \leq k-2$ and/or $\lambda_{k}' \leq k-2$.
If $\lambda_{k}' \leq k-2$, then $\lambda_{k} \leq k-1$ as $\lambda$ is a Young diagram.
Similarly, if $\lambda_{k} \leq k-2$, then $\lambda_{k}' \leq k-1$ as $\lambda$ is a Young diagram.
In both cases it follows that $\lambda_{k}' = |\mathrm{Ver}_{k}|$ and $\lambda_{k} = |\mathrm{Hor}_{k}|$.
If $|\mathrm{Hor}_{k} \cup \mathrm{Ver}_{k}|$ is odd, $\lambda_{k} + \lambda_{k}'$ is also odd as
$\lambda_{k} + \lambda_{k}' = |\mathrm{Hor}_{k}| + |\mathrm{Ver}_{k}| = 
|\mathrm{Hor}_{k} \cup \mathrm{Ver}_{k}|$, and clearly
the right hand side of (\ref{eq:politesociety}) equals $-1$.
Thus (\ref{eq:politesociety}) is true for each $k$, 
from which it follows that (\ref{eq:makethingsright}) is true,
which completes the proof of the lemma.
\end{proof}

\begin{lemma}
\label{lem:theredoesnotexist}
There does not exist a non-zero projector $E_{SS} \in {\cal{C}}_{t}$ that projects onto one
irreducible $U_{q}(\mathfrak{g})$-summand of $V^{\otimes t}$ with the following property:
\begin{equation}
\label{eq:ESS1} 
\left(E_{SS}V^{\otimes t}\right) \cap \big(\Upsilon(e_{RR})V^{\otimes t}\big) = 0,
\end{equation}
for all idempotent matrix units $e_{RR} \in \widetilde{BW}_{t}(-q^{2n},q)$. 
\end{lemma}
\begin{proof} 
Suppose that such a projector did exist, then 
$E_{SS} \Upsilon(e_{RR})=0$ for each idempotent matrix unit   
$e_{RR} \in \widetilde{BW}_{t}(-q^{2n},q)$, and $E_{SS} \Upsilon(e_{RT})=0$ for each matrix unit 
$e_{RT} \in \widetilde{BW}_{t}(-q^{2n},q)$ where $R \neq T$ and
$(R,T) \in \Omega_{t}(-q^{2n},q)$ as
$$E_{SS} \Upsilon(e_{RT}) = E_{SS} \Upsilon(e_{RR}e_{RT}e_{TT}) = 0
 = \Upsilon(e_{RR}e_{RT}e_{TT})E_{SS} = \Upsilon(e_{RT})E_{SS},$$ 
 as $\Upsilon$ is a homomorphism.
We will show that it then follows that $E_{SS} \in {\cal{J}}_{t}$. 
However, the assumption that $E_{SS}$ is a non-zero projector and that each irreducible 
$U_{q}(\mathfrak{g})$-summand of $V^{\otimes t}$ has non-zero quantum superdimension means that
 $\mathrm{str}_{q}(E_{SS}) = sdim_{q}(E_{SS}V^{\otimes t}) \neq 0$ 
and thus that $E_{SS} \notin {\cal{J}}_{t}$, proving the lemma.

We now show that any such projector $E_{SS}$ belongs to ${\cal{J}}_{t}$.
To see this, assume the contrary, then
\begin{equation}
\label{eq:john_singleton(19)}
E_{SS} = \widetilde{E}_{SS} + E_{j},
\end{equation}
 where
$\widetilde{E}_{SS} \in \Upsilon\big(\widetilde{BW}_{t}(-q^{2n},q) \big)$
and $E_{j}$ is some (potentially vanishing) element of ${\cal{J}}_{t}$ from (\ref{eq:milindaandmilinda}).
Explicitly,
\begin{equation}
\label{eq:illbetheretosupportyou}
\widetilde{E}_{SS} = \sum_{(R,T) \in \Omega_{t}(-q^{2n},q)} c_{RT} \Upsilon(e_{RT}), 
\hspace{5mm} c_{RT} \in \mathbb{C},
\end{equation}
where the assumption that $E_{SS} \notin {\cal{J}}_{t}$ means that
$c_{RT} \neq 0$ for at least one pair $(R,T) \in \Omega_{t}(-q^{2n},q)$.
Assume that $(A, B) \in \Omega_{t}(-q^{2n},q)$ is a pair where $c_{AB} \neq 0$, then
\begin{equation}
\label{eq:constructivist}
\mathrm{str}_{q} \left(\Upsilon(e_{BA})\widetilde{E}_{SS} \right) = 
\mathrm{str}_{q} \left( \sum_{T} c_{AT} \Upsilon(e_{BT}) \right) 
 =  c_{AB} \mathrm{str}_{q} \left( \Upsilon(e_{BB}) \right) \neq  0,
\end{equation}
where the sum is over all paths $T$ of length $t$ in the Bratteli diagram for
$\widetilde{BW}_{t}(-q^{2n},q)$, as $c_{AB} \neq 0$ and 
$\mathrm{str}_{q}(\Upsilon(e_{BB})) \neq 0$.

However, (\ref{eq:constructivist}) is not true. To see this, recall that  
$E_{SS}$ satisfies
\begin{equation}
\label{eq:caseychamberscat}
\Upsilon(e_{BA}) E_{SS} = \Upsilon(e_{BA})\big(\widetilde{E}_{SS} + E_{j} \big) = 0,
\end{equation} 
and note that
\begin{equation}
\label{eq:eBA1}
\Upsilon(e_{BA})\widetilde{E}_{SS} \in \Upsilon \big( \widetilde{BW}_{t}(-q^{2n},q) \big)
\end{equation}
and that $\big(\Upsilon(e_{BA})E_{j}\big) \in {\cal{J}}_{t}$,
which we recall arises from the definition of ${\cal{J}}_{t}$. 
Now from (\ref{eq:eBA1}) and (\ref{eq:milindaandmilinda}),
$\Upsilon(e_{BA})\widetilde{E}_{SS} \neq -\Upsilon(e_{BA})E_{j}$ 
unless  $\Upsilon(e_{BA})\widetilde{E}_{SS} = 0 = -\Upsilon(e_{BA})E_{j}$
and then (\ref{eq:caseychamberscat}) implies that 
$\Upsilon(e_{BA})\widetilde{E}_{SS} = \Upsilon(e_{BA})E_{j} = 0$.

It follows from the fact that $\Upsilon(e_{BA})\widetilde{E}_{SS} = 0$ that 
$\mathrm{str}_{q} \big( \Upsilon(e_{BA})\widetilde{E}_{SS} \big) = 0$, however, this 
contradicts (\ref{eq:constructivist}).
Then the assumption in (\ref{eq:illbetheretosupportyou}) that
$c_{RT} \neq 0$ for at least one pair $(R,T) \in \Omega_{t}(-q^{2n},q)$ is false, thus
it must be true that $\widetilde{E}_{SS} = 0$.

It then follows that $E_{SS} = E_{j} \in {\cal{J}}_{t}$ from (\ref{eq:john_singleton(19)}). 
However, this is false as
$\mathrm{str}_{q}(E_{SS}) = sdim_{q}(E_{SS}V^{\otimes t}) \neq 0$.
Thus, our original assumption that there exists 
a projector $E_{SS} \in {\cal{C}}_{t}$ with the property that 
$\left(E_{SS}V^{\otimes t}\right) \cap \big( \Upsilon(e_{RR})V^{\otimes t} \big) = 0$ 
for all idempotent matrix units $e_{RR} \in \widetilde{BW}_{t}(-q^{2n},q)$ is not true.
This completes the proof.

\end{proof}

\end{subsection}

\begin{subsection}{Matrix units in $\mathcal{C}_{t}$}
\label{subsection:genericcalmatrixunitsinCt}

We have not yet proved that ${\cal{L}}_{t} = {\cal{C}}_{t}$. 
We will complete its proof in this subsection by defining a complete set of 
intertwiners in ${\cal{C}}_{t}$
that we obtain by applying the map $\Upsilon$ to the intertwiner matrix units in
$\widetilde{BW}_{t}(-q^{2n},q)$.

There is a canonical bijection between paths in the Bratteli diagram for 
$\widetilde{BW}_{t}(-q^{2n},q)$
and paths in the Bratteli diagram for $V^{\otimes t}$.
Recall that each Young diagram on the $k^{th}$ level of the Bratteli diagram for 
$\widetilde{BW}_{t}(-q^{2n},q)$ contains an even (resp. odd) number of boxes
if $k$ is an even (resp. odd) number.
Each vertex $\lambda$ on the $k^{th}$ level of the Bratteli diagram for 
$\widetilde{BW}_{t}(-q^{2n},q)$ also appears on the $k^{th}$ level of the Bratteli diagram for
$V^{\otimes t}$ unless $\lambda$ has more than $n$ rows of boxes, in which case the Young diagram
$\widetilde{\lambda}$, defined in Lemma \ref{lem:reflectedyoungdiagram}, appears instead.

Given a path $\widetilde{S}$ of length $t$ in the Bratteli diagram for $V^{\otimes t}$, 
we now determine the corresponding path $S$ of length $t$ in the Bratteli diagram for 
$\widetilde{BW}_{t}(-q^{2n},q)$. 
Let $\widetilde{S} = (0, s_{1}, \ldots, s_{t})$ where $s_{i}$ is a Young diagram on the
$i^{th}$ level of the Bratteli diagram for $V^{\otimes t}$. 
If $i$ is even (resp. odd) and $s_{i}$ contains an even (resp. odd) number of boxes, then 
$s_{i}$ is also a vertex on the $i^{th}$ level of the Bratteli diagram for 
$\widetilde{BW}_{t}(-q^{2n},q)$.
If, however, $i$ is even (resp. odd) and $s_{i}$ contains an odd (resp. even) number of
boxes, then $s_{i}=\widetilde{\lambda}$ is the vertex obtained from a vertex $\lambda$ on
the $i^{th}$ level of the Bratteli diagram for  $\widetilde{BW}_{t}(-q^{2n},q)$ given by
$\widetilde{\lambda}'_{1} = 2n+1 - \lambda'_{1}$ 
and $\widetilde{\lambda}'_{j} = \lambda_{j}'$ for $j \geq 2$.
This allows us to define the path $S$ of length $t$ in the Bratteli diagram for
$\widetilde{BW}_{t}(-q^{2n},q)$ corresponding to 
$\widetilde{S}$ and it also
gives rise to the
canonical bijection between the set of paths of length $t$ in the Bratteli 
diagram for $\widetilde{BW}_{t}(-q^{2n},q)$ and the set of paths of length $t$ in the Bratteli diagram
for $V^{\otimes t}$.

In Figure \ref{fig:isomorphismsofbrattelidiagrams}
we show the isomorphism between the Bratteli diagram for $\widetilde{BW}_{t}(-q^{2},q)$ 
up to the $4^{th}$ level (on the left) with the Bratteli diagram for $V^{\otimes t}$
up to the $4^{th}$ level (on the right), where
$V$ is the $3$-dimensional irreducible $U_{q}(osp(1|2))$-module.
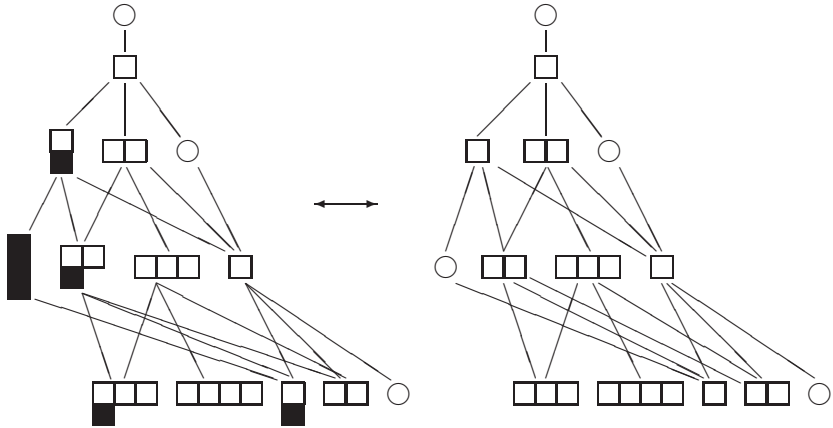
\begin{figure}[hbt]
\label{fig:isomorphismsofbrattelidiagrams}
\begin{center}
\setlength{\unitlength}{0.14mm}
\begin{picture}(960,440) \large\sf
\thinlines
%\linethickness{0.05mm}

\put(400, 250){\vector(1, 0){30}} 
\put(400, 250){\vector(-1, 0){30}}

%first diagram

\put(190,430){\circle{20}} 

\put(190,395){\line(0,1){20}} 

\put(180,370){\framebox(20,20){}}

\put(120,300){\framebox(20,20){}}
\linethickness{1.5mm}
\put(130,290){\framebox(0,0){}}
\thinlines

\put(170,290){\framebox(20,20){}}
\put(190,290){\framebox(20,20){}}

\put(250,300){\circle{20}} 

\linethickness{1.5mm}
\put(90,170){\framebox(0,0){}}
\put(90,190){\framebox(0,0){}}
\put(90,210){\framebox(0,0){}}
\thinlines

\put(100,225){\line(1,2){25}} 
\put(130,275){\line(1,-4){15}} 

\linethickness{1.5mm}
\put(140,180){\framebox(0,0){}}
\thinlines
\put(130,190){\framebox(20,20){}}
\put(150,190){\framebox(20,20){}}
\put(152,215){\line(1,2){35}} 
\put(192,285){\line(1,-2){40}} 

\put(200,180){\framebox(20,20){}}
\put(220,180){\framebox(20,20){}}
\put(240,180){\framebox(20,20){}}

\put(290,180){\framebox(20,20){}}
\put(215,285){\line(1,-1){79}} 
\put(260,285){\line(1,-2){40}} 
\put(145,275){\line(2,-1){140}} 

\put(135,325){\line(1,1){40}} 
\put(190,315){\line(0,1){50}} 
\put(205,365){\line(3,-4){38}} 

\linethickness{1.5mm}
\put(170,50){\framebox(0,0){}}
\thinlines
\put(160,60){\framebox(20,20){}}
\put(180,60){\framebox(20,20){}}
\put(200,60){\framebox(20,20){}}
\put(150,165){\line(1,-3){27}} 
\put(190,85){\line(1,3){30}} 
\put(105,160){\line(3,-1){230}} 

\put(240,60){\framebox(20,20){}}
\put(260,60){\framebox(20,20){}}
\put(280,60){\framebox(20,20){}}
\put(300,60){\framebox(20,20){}}
\put(220,175){\line(1,-2){45}} 

\put(340,60){\framebox(20,20){}}
\linethickness{1.5mm}
\put(350,50){\framebox(0,0){}}
\thinlines
\put(150,165){\line(5,-2){196}} 
\put(305,175){\line(1,-1){90}} 
\put(305,175){\line(1,-2){45}} 
\put(305,175){\line(3,-2){138}} 

\put(380,60){\framebox(20,20){}}
\put(400,60){\framebox(20,20){}}
\put(150,165){\line(3,-1){240}} 
\put(220,175){\line(2,-1){180}} 

\put(450,70){\circle{20}}

%second diagram - 400 points across

\put(590,430){\circle{20}} 

\put(590,395){\line(0,1){20}} 

\put(580,370){\framebox(20,20){}}

\put(515,290){\framebox(20,20){}}
%\put(520,300){\framebox(20,20){}}
%\put(520,280){\framebox(20,20){}}

\put(570,290){\framebox(20,20){}}
\put(590,290){\framebox(20,20){}}

\put(650,300){\circle{20}} 

\put(495,190){\circle{20}} 
%\put(480,160){\framebox(20,20){}}
%\put(480,180){\framebox(20,20){}}
%\put(480,200){\framebox(20,20){}}
\put(495,205){\line(1,3){27}} 
\put(550,205){\line(-1,4){20}} 

\put(530,180){\framebox(20,20){}}
\put(550,180){\framebox(20,20){}}
%\put(530,170){\framebox(20,20){}}
%\put(530,190){\framebox(20,20){}}
%\put(550,190){\framebox(20,20){}}
\put(550,205){\line(1,2){40}} 
\put(592,285){\line(1,-2){40}} 

\put(600,180){\framebox(20,20){}}
\put(620,180){\framebox(20,20){}}
\put(640,180){\framebox(20,20){}}

\put(690,180){\framebox(20,20){}}
\put(615,285){\line(1,-1){80}} 
\put(660,285){\line(1,-2){40}} 
\put(545,285){\line(5,-3){140}} 

\put(525,315){\line(1,1){50}} 
\put(590,315){\line(0,1){50}} 
\put(605,365){\line(3,-4){38}} 

%\put(560,60){\framebox(20,20){}}
%\put(560,40){\framebox(20,20){}}
%\put(580,60){\framebox(20,20){}}
%\put(600,60){\framebox(20,20){}}
\put(560,60){\framebox(20,20){}}
\put(580,60){\framebox(20,20){}}
\put(600,60){\framebox(20,20){}}
\put(550,175){\line(1,-3){30}} 
\put(590,85){\line(1,3){30}} 
\put(505,175){\line(5,-2){230}} 

\put(640,60){\framebox(20,20){}}
\put(660,60){\framebox(20,20){}}
\put(680,60){\framebox(20,20){}}
\put(700,60){\framebox(20,20){}}
\put(635,175){\line(1,-2){45}} 

%\put(740,60){\framebox(20,20){}}
%\put(740,40){\framebox(20,20){}}
\put(740,60){\framebox(20,20){}}
\put(560,175){\line(2,-1){180}}
 
\put(705,175){\line(1,-1){90}} 
\put(700,175){\line(1,-2){45}} 
\put(710,175){\line(3,-2){135}} 

\put(780,60){\framebox(20,20){}}
\put(800,60){\framebox(20,20){}}
\put(575,180){\line(2,-1){200}} 
\put(640,175){\line(5,-3){150}} 

\put(850,70){\circle{20}} 

\end{picture}
\caption{The isomorphism between the Bratteli diagram for $\widetilde{BW}_{t}(-q^{2},q)$ 
up to the $4^{th}$ level (on the left) with the Bratteli diagram for $V^{\otimes t}$
up to the $4^{th}$ level (on the right), where
$V$ is the $3$-dimensional irreducible $U_{q}(osp(1|2))$-module}  
\end{center}
\end{figure}

The diagram on the left in Figure \ref{fig:isomorphismsofbrattelidiagrams} is the Bratteli diagram for
$\widetilde{BW}_{t}(-q^{2},q)$ up to the $4^{th}$ level with some of the 
boxes of the vertices filled in. 
If $\lambda$ is a vertex in the Bratteli diagram for $\widetilde{BW}_{t}(-q^{2},q)$, then
$\widetilde{\lambda}$, as given in Lemma \ref{lem:reflectedyoungdiagram}, 
is precisely the Young diagram $\lambda$ with the filled in boxes removed.

We obtain the intertwiners in ${\cal{C}}_{t}$ between the isomorphic irreducible
$U_{q}(\mathfrak{g})$-submodules of $V^{\otimes t}$ defined by 
the projectors $E_{RR} \in {\cal{C}}_{t}$ where
$R$ is a path of length $t$ in the Bratteli diagram for $V^{\otimes t}$,
by applying the homomorphism $\Upsilon$ to the intertwiner matrix units in $\widetilde{BW}_{t}(-q^{2n},q)$.
All we need to do is check that the images of the intertwiners 
under $\Upsilon$ are all well-defined and non-zero.

We construct the intertwiners in ${\cal{C}}_{t}$ recursively.  
To do this, assume that all the matrix units in ${\cal{C}}_{t-1}$ have already been defined 
and that they are non-zero. 
This condition is satisfied for $t=3$ as
the decomposition of $V \otimes V$ into irreducible $U_{q}(\mathfrak{g})$-submodules 
is multiplicity free, thus no intertwiners exist in ${\cal{C}}_{2}$.

In the remainder of this subsection, let $\widetilde{M}$ and $\widetilde{P}$ 
be a pair of paths of length $t$ in the Bratteli diagram for $V^{\otimes t}$
where $shp(\widetilde{M}) = shp(\widetilde{P})$ and $\widetilde{M} \neq \widetilde{P}$. 
Let $M$ and $P$ be the corresponding paths in the Bratteli diagram for
$\widetilde{BW}_{t}(-q^{2n},q)$.
The intertwiner $E_{\widetilde{M}\widetilde{P}} \in {\cal{C}}_{t}$ is defined by
$E_{\widetilde{M}\widetilde{P}} = \Upsilon(e_{MP})$.

Let us firstly deal with the situation that $shp(M)=shp(P)=\lambda$ where
$\lambda$ contains strictly fewer than $t$ boxes. 
Referring back to Subsection \ref{subsect:matrixbirmanwenzlstuff} we see that
\begin{eqnarray*}
E_{\widetilde{M}\widetilde{P}} & = & \Upsilon(e_{MP}) 
= \frac{Q_{\lambda}(-q^{2n},q)}{\sqrt{Q_{\mu}(-q^{2n},q)Q_{\widetilde{\mu}}(-q^{2n},q)}}
                 E_{\widetilde{M}'\widetilde{S}} \Upsilon(e_{t-1}) E_{\widetilde{T}\widetilde{P}'}, \\
& &  \mbox{where } e_{t-1} \in BW_{t}(-q^{2n},q),
\end{eqnarray*}
where $S$ and $T$ are paths of length $t-1$ in the Bratteli diagram for 
$\widetilde{BW}_{t}(-q^{2n},q)$ satisfying the conditions: 
\begin{itemize}
\item[(i)]   $shp(S)=shp(M') = \mu$, and
\item[(ii)]  $shp(T)=shp(P') = \widetilde{\mu}$, and
\item[(iii)] $S'=T'$, and
\item[(iv)]  $shp(S')= \lambda = shp(T')$.
\end{itemize}
Recall that these paths exist. Note that $Q_{\mu}(-q^{2n},q) \neq 0$ and
$Q_{\widetilde{\mu}}(-q^{2n},q) \neq 0$; if it were true that
$Q_{\mu}(-q^{2n},q) = 0$ then $\mu$ would not be a vertex
in the Bratteli diagram for $\widetilde{BW}_{t}(-q^{2n},q)$.
Similar remarks hold for $\widetilde{\mu}$ and $\lambda$. 
Thus $E_{\widetilde{M}\widetilde{P}}$ is well-defined.
(We will later show that $E_{\widetilde{M}\widetilde{P}}$ is non-zero.)

Now let us deal with the situation that $shp(M)=shp(P)=\lambda$ where
$\lambda$ contains exactly $t$ boxes and $shp(M') = shp(P')$.
Referring back to Subsection \ref{subsect:matrixbirmanwenzlstuff}, we see that
$$E_{\widetilde{M}\widetilde{P}} =\Upsilon(e_{MP})
 = \Upsilon\big((1-z_{t}) o_{MP}\big),$$
 where $o_{MP} = o_{M'P'}o_{PP}$ and
 $z_{t} = \sum_{S}e_{SS}$
 with the summation going over all paths $S$ of length $t$ such that $shp(S)$ contains fewer than $t$
 boxes. 

Now let us deal with the situation that $shp(M)=shp(P)=\lambda$ where
$\lambda$ contains exactly $t$ boxes and $shp(M') \neq shp(P')$.
Choose paths $\overline{M}$ and $\overline{P}$ of length $t$ such that 
$shp(\overline{M}) = shp(M)$ and $shp(\overline{P}) = shp(P)$ and
\begin{itemize}
\item[(i)] $\overline{M}'' = \overline{P}''$, and
\item[(ii)] $shp(\overline{M}') = shp(M')$, and
\item[(iii)] $shp(\overline{P}') = shp(P')$.
\end{itemize}
Such paths can always be chosen. Then
$$E_{\widetilde{M}\widetilde{P}} =\Upsilon(e_{MP})
 = \Upsilon\big((1-z_{t}) o_{MP}\big), \hspace{5mm} \mbox{where}$$
 \begin{equation}
 \label{eq:thaiequation}
 o_{MP} = \frac{1-q^{2d}}{\sqrt{(1-q^{2(d+1)})(1-q^{2(d-1)})}} 
           o_{M'\overline{M}'} g_{t-1} o_{\overline{P}'P'}o_{PP},
	   \hspace{2mm} g_{t-1} \in BW_{t}(-q^{2n},q),
\end{equation}
where $d = d(\overline{M},t-1)$ is the integer defined by (\ref{eq:brontethewhalecatslothdogmouse}).
The integer $|d(\overline{M},i)|+1$ is the number of boxes in the hook going through the boxes 
containing the numbers $i$ and $(i+1)$ \cite{tw}.

We now prove that the coefficient on the right hand side of (\ref{eq:thaiequation}) 
is well-defined and non-zero.
It is not difficult to see that the coefficient is well-defined if $|d| \neq 1$, 
and we now show that this is always true.  
As $|d| + 1$ is the length of the hook going through the boxes 
containing the numbers $(t-1)$ and $t$, it is always true that $|d| + 1 \geq 2$ as 
each such hook contains at least two boxes.
Now the only situation in which it could possibly be true that $|d| = 1$ is if the boxes 
containing the numbers $(t-1)$ and $t$ are immediately horizontally or vertically adjacent.
However this cannot occur for the following reason: from the above construction, 
the number $t$ is in the same box in $\overline{M}$ as the number $(t-1)$ is in $\overline{P}$,
and the number $t$ is in the same box in $\overline{P}$ that the number $(t-1)$ is in $\overline{M}$.
It follows  that if the numbers $(t-1)$ and $t$ are 
immediately horizontally or vertically adjacent in $\overline{M}$, 
each must be in the corresponding `swapped' box in $\overline{P}$, 
and then at least one of $\overline{M}$ or $\overline{P}$
{\emph{cannot}} be a standard tableau.  
This contradicts the assumption that both $\overline{M}$ and $\overline{P}$
are standard tableaux, thus $|d| \neq 1$  and the coefficient in
(\ref{eq:thaiequation}) is well-defined.

It remains for us to show that the coefficient in (\ref{eq:thaiequation}) is non-zero.
This follows immediately from the fact that $|d| \neq 0$.
Note that we have not yet proved that
the matrix units are all non-zero.

Let us write $E_{MP}$ to denote $E_{\widetilde{M}\widetilde{P}}$.
We note that the matrix unit $E_{MP} \in {\cal{C}}_{t}$, for $M \neq P$, is 
an intertwiner between the isomorphic irreducible $U_{q}(\mathfrak{g})$-modules
 $E_{PP} (V^{\otimes t})$ and $E_{MM} (V^{\otimes t})$:
 $$E_{MP}: E_{PP} (V^{\otimes t}) \rightarrow E_{MM} (V^{\otimes t}),$$
and that the whole collection of matrix units satisfy
$$E_{QR}E_{ST} = \delta_{RS} E_{QT}.$$
To show that each intertwiner $E_{MP}$ is non-zero, it suffices to note that
each projector $E_{PP}$ is non-zero and that $E_{PP} = E_{PM}E_{MP}$.

We then have the complete set of projectors and intertwiners in ${\cal{C}}_{t}$.
This means that  
${\cal{L}}_{t} = {\cal{C}}_{t}$.
To see this, 
note that the matrix units 
$\{E_{ST} = \Upsilon(e_{ST}) | \ (S,T) \in \Omega_{t}(-q^{2n},q)\}$ are a basis for ${\cal{C}}_{t}$,
as the fact that the superdimension of each finite dimensional irreducible
$U_{q}(\mathfrak{g})$-module $V_{\lambda}$ with integral dominant highest weight
$\lambda$ is non-zero means that Schur's lemma has the same form as it does for ungraded quantum algebras, thus the centraliser ${\cal{L}}_{t}$ is generated by the complete set 
$\{E_{ST} = \Upsilon(e_{ST}) | \ (S,T) \in \Omega_{t}(-q^{2n},q)\}$ of projectors and intertwiners.

Note that ${\cal{J}}_{t} = 0$; to see this,
let $X$ be an arbitrary element of ${\cal{C}}_{t}$, then
$$X = \sum_{(S,T) \in \Omega_{t}(-q^{2n},q)} x_{ST} E_{ST}, \hspace{10mm}
x_{ST} \in \mathbb{C},$$ 
where $x_{ST} \neq 0$ for at least one pair $(S,T)$.  Let $(A,B)$ be such a pair, then
$$\mathrm{str}_{q}(E_{BA}X) = \mathrm{str}_{q}(x_{AB}E_{BA}E_{AB}) = 
x_{AB} \mathrm{str}_{q}(E_{BB}) \neq 0,$$
thus $X  \notin {\cal{J}}_{t}$.  As $X$ was arbitrary, ${\cal{J}}_{t} = 0$.

Note that we obtained  ${\cal{C}}_{t} = {\cal{L}}_{t}$ by using the fact that
$\lambda$ and $\widetilde{\lambda}$ do not appear on the same level of the Bratteli diagram for
$\widetilde{BW}_{t}(-q^{2n},q)$.  
If $\lambda$ and $\widetilde{\lambda}$ did appear on the same level, we could only conclude 
that there is an inclusion of ${\cal{C}}_{t}$  in ${\cal{L}}_{t}$ rather than an equality. 
Of course, in that event, there  may actually be an equality, but
a different method would have to be used to obtain all the intertwiners.

We now present the technical lemmas used above.
\begin{lemma}
\label{lemma:equalityoftraces}
Let $\psi: {\cal{C}}_{t} \rightarrow \mathbb{C}$ be a map defined by
$$\psi(X) = \frac{\mathrm{str}_{q}(X)}{\big(sdim_{q}(V)\big)^{t}},$$
and let $\mathrm{tr}$ be the trace functional on $BW_{t}(-q^{2n},q)$ 
mentioned in (\ref{eq:tracefunctional}).
Then 
$$\psi \big( \Upsilon (a) \big) = \mathrm{tr}(a), \hspace{10mm} \forall a \in 
BW_{t}(-q^{2n},q).$$
\end{lemma}
\begin{proof}
Any functional $\phi$ on $BW_{\infty}(-q^{2n},q)$ satisfying
Eq. (\ref{eq:tracefunctional}) for all $t \in \mathbb{N}$ is identical to 
$\mathrm{tr}$ \cite[Lem. 3.4 (d)]{w2}, and we will show that
$\psi \circ \Upsilon$ is such a functional.

To show that $\psi \circ \Upsilon$ satisfies Eq. (\ref{eq:tracefunctional}), it suffices
to show that for each $t \in \mathbb{N}$,
\begin{equation}
 \label{eq:julietoolie(1)}
-\psi\big(\Upsilon(a) \check{R}_{t-1} \Upsilon(b)\big) = 
-\psi\big( \check{R}_{t-1}\big) \psi\big(\Upsilon(ab)\big), \hspace{5mm}
 \forall a, b \in BW_{t-1}(-q^{2n},q),
 \end{equation}
as the element
$e_{t-1} \in BW_{t}(-q^{2n},q)$ can be written as a function of the $g_{t-1}$'s.  
We will show that Eq. (\ref{eq:julietoolie(1)}) is true 
using Lemmas \ref{lem;kilo} and \ref{ACER1}.

The left hand side of Eq. (\ref{eq:julietoolie(1)}) is
\begin{equation}
\label{eq:julietoolie(2)}
-\mathrm{str}_{q}^{\otimes t} \big(A \check{R}_{t-1} B \big) / \big(sdim_{q}(V) \big)^{t},
\end{equation}
where $\mathrm{str}_{q}^{\otimes t}$ indicates  
that we take the quantum supertrace over all $t$ tensor factors,
and we also write $A = \Upsilon(a)$ and $B = \Upsilon(b)$. 
Now we can regard each $X \in {\cal{C}}_{t-1}$ as an element of ${\cal{C}}_{t}$ under the mapping
$X \mapsto X \otimes \mathrm{id}$, 
then by applying the identity to the first $t-1$ tensor powers of (\ref{eq:julietoolie(2)}) and
taking the quantum supertrace over
the $t^{th}$ tensor power of (\ref{eq:julietoolie(2)}), we obtain, using Lemmas \ref{lem;kilo} and \ref{ACER1},
\begin{equation}
\label{eq:julietoolie(3)}
-\mathrm{str}_{q}^{\otimes t} \big(A \check{R}_{t-1} B \big)/ \big(sdim_{q}(V) \big)^{t}
  = \frac{-\chi_{V}(v^{-1})}{sdim_{q}(V)} \ 
  \frac{\mathrm{str}_{q}^{\otimes (t-1)} \big(A B \big)}{ \big(sdim_{q}(V) \big)^{t-1}}.
\end{equation}
 Now 
 $$\psi\big( \check{R}_{t-1}\big) = \chi_{V}(v^{-1})/sdim_{q}(V),$$
 and the right hand side of (\ref{eq:julietoolie(3)}) equals
 the right hand side of (\ref{eq:julietoolie(1)}).
Now (\ref{eq:julietoolie(1)}) is true for all $a$ and $b$ belonging to $BW_{t-1}(-q^{2n},q)$,
and it remains to show that $\psi \circ \Upsilon$ is a functional on 
$BW_{\infty}(-q^{2n},q)$ satisfying (\ref{eq:tracefunctional}) for 
all $t=1, 2, \ldots$.
This follows from the fact that
$\psi(A \otimes \mathrm{id}) = \psi(A)$ for all $A \in {\cal{C}}_{t}$, thus we can regard
$\psi$ as well-defined in the inductive limit  
${\cal{C}}_{2} \subset {\cal{C}}_{3} \subset {\cal{C}}_{4} \subset \cdots$, which
completes the proof.

\end{proof}

The following lemma, which we used in the proof of Lemma \ref{lemma:equalityoftraces}, appears in
\cite[Lem. 2]{lg} and is proved in \cite[Lem. 3.1]{z1}.
\begin{lemma}
\label{lem;kilo}
Let $V$ be the fundamental $(2n+1)$-dimensional irreducible
$U_{q}(osp(1|2n))$-module with highest weight $\epsilon_{1}$ and
let $\pi$ be the representation of $U_{q}(osp(1|2n))$ afforded by $V$.  
Let 
$\check{\cal{R}}_{V,V} \in \mathrm{End}_{U_{q}(osp(1|2n))}(V \otimes V)$ be as given in
Eq. (\ref{eq:bigjimmyboy(2)}). Then
$$(\mathrm{id} \otimes \mathrm{str})
\left[(\mathrm{id} \otimes \pi(K_{2\rho}))  \check{\cal{R}}_{V,V}^{\pm 1} \right] 
= q^{\pm(\epsilon_{1}, \epsilon_{1} + 2\rho)}\mathrm{id}_{V}=\chi_{V}(v^{\mp 1})\mathrm{id}_{V}.$$
\end{lemma}

\begin{lemma}
\label{ACER1}
For all $x, y \in {\cal{C}}_{t}$, 
\begin{equation}
\label{bigbrother22006}
(\mathrm{id}^{\otimes t} \otimes \mathrm{str}_{q}) x \check{R}_{t} y = 
   \left[ (\mathrm{id}^{\otimes t} \otimes \mathrm{str}_{q}) \check{R}_{t} \right] xy,
\end{equation}
where each element $z \in {\cal{C}}_{t}$ is embedded in ${\cal{C}}_{t+1}$ via the mapping 
$z \mapsto z \otimes \mathrm{id}$.
\end{lemma}
\begin{proof}
Let $W$ be a finite dimensional integrable $U_{q}(osp(1|2n))$-module and $V$ be the $(2n+1)$-dimensional
irreducible $U_{q}(osp(1|2n))$-module.
Let $\{w^{i} | \ i \in I\}$ be a homogeneous basis of $W$ and
$\{v^{j} | \ j \in J\}$ be the homogeneous basis of $V$ used throughout this paper 
where $v^{j}=v_{j}$.

We write the action of $B \in \mathrm{End}_{\mathbb{C}}(W)$ as
$B w^{i} = \sum_{k \in I} B^{k}_{i} w^{k}$, $B^{k}_{i} \in \mathbb{C}$, and the action of 
$A \in \mathrm{End}_{\mathbb{C}}(W \otimes V)$ as
$$A (w^{i} \otimes v^{j}) = \sum_{k \in I, \ l \in J} A^{kl}_{ij} w^{k} \otimes v^{l}, 
\hspace{5mm} A^{kl}_{ij} \in \mathbb{C},$$
with the obvious generalisations to further tensor powers.
For such an $A$, $\left[(\mathrm{id} \otimes \mathrm{tr})A \right] \in 
\mathrm{End}_{\mathbb{C}}(W)$ 
with the action
$$\left[ (\mathrm{id} \otimes \mathrm{tr})A \right]^{i}_{k} = \sum_{j \in J} A^{ij}_{kj}.$$
Similarly 
$\left[ (\mathrm{id} \otimes \mathrm{str})A \right]^{i}_{k} = \sum_{j \in J}(-1)^{[v^{j}]} A^{ij}_{kj}$
and
$$\left[ (\mathrm{id} \otimes \mathrm{str}_{q})A \right]^{i}_{k} = 
       \sum_{j \in J}(-1)^{[v^{j}]} \langle \left. v^{j} \right.^{*} , K_{2 \rho} v^{j} \rangle A^{ij}_{kj},$$
where $\langle \cdot, \cdot \rangle$ is the dual pairing:
$\langle \left. v^{j} \right.^{*}, v^{k} \rangle = \delta_{jk}$, 
and $[\left. v^{j} \right.^{*} ] = [v^{j}]$. 

Clearly, 
\begin{equation}
\label{bigbrotherdavid1}
(\mathrm{id}^{\otimes t} \otimes \mathrm{str}_{q})\check{R}_{t} = \chi_{V}(v^{-1}) (\mathrm{id}_{V})^{\otimes t},
\end{equation}
from Lemma \ref{lem;kilo}.
Now let $$\left\{ v^{\overline{i}} = v^{i_{1}} \otimes v^{i_{2}} \otimes \cdots \otimes v^{i_{t-1}} \left| \right. \ 
\overline{i} = (i_{1}, i_{2} \ldots, i_{t-1}) \in \overline{J} = J^{\times (t-1)} \right\}$$ 
be a homogeneous basis of $V^{\otimes (t-1)}$, then we can write $x \in {\cal{C}}_{t}$ as
$$x = \left(X^{\overline{i} k}_{\overline{j} l}\right)_{\overline{i}, \overline{j} \in \overline{J}, k,l \in J} \hspace{5mm}
\mathrm{and } \hspace{5mm}
x \otimes \mathrm{id}_{V} = \left(X^{\overline{i} k r}_{\overline{j} l s}\right)_{\overline{i}, \overline{j} \in \overline{J}, k, l, r, s \in J}$$
where 
\begin{equation}
\label{eq:GEM}
X^{\overline{i} k r}_{\overline{j} l s} = \delta_{rs} X^{\overline{i} k}_{\overline{j} l}.
\end{equation}
In this notation, we can rewrite (\ref{bigbrotherdavid1}): let
$$\check{R}_{t} = \left(R^{\overline{i} k r}_{\overline{j} l s} \right)_{\overline{i}, \overline{j} \in \overline{J}, k, l, r, s \in J}$$
then we have
\begin{equation}
\label{eq:GEM3}
\left[ (\mathrm{id}^{\otimes t} \otimes \mathrm{str}_{q}) \check{R}_{t} \right]^{\overline{i} k}_{\overline{j} l} = 
\sum_{p} (-1)^{[v^{p}]} \langle \left. v^{p} \right.^{*}, K_{2 \rho} v^{p} \rangle R^{\overline{i} k p}_{\overline{j} l p} = 
\delta_{\overline{i} \overline{j}} \delta_{kl} \chi_{V}(v^{-1}).
\end{equation}

Now let $x$ and $y$ be arbitrary elements of ${\cal{C}}_{t}$ embedded in  
${\cal{C}}_{t+1}$ via the map $z \mapsto z \otimes \mathrm{id}_{V}$.
To determine the left hand side of (\ref{bigbrother22006}), we note that
\begin{eqnarray*}
\left[ (x \otimes \mathrm{id}_{V}) \check{R}_{t} (y \otimes \mathrm{id}_{V})   
   \right]^{\overline{i} k r}_{\overline{j} l s} & = & 
    \sum_{ \overline{a}, \overline{d} \in \overline{J}, b, c, e, f \in J } 
            X^{\overline{i} k r}_{ \overline{a} b c  } 
            R^{ \overline{a} b c  }_{ \overline{d} e f } 
            Y^{ \overline{d} e f }_{\overline{j} l s} \\ 
 & = & 
    \sum_{ \overline{a}, \overline{d} \in \overline{J}, b, e \in J } 
            X^{\overline{i} k r}_{ \overline{a} b r  } 
            R^{ \overline{a} b r  }_{ \overline{d} e s } 
            Y^{ \overline{d} e s }_{\overline{j} l s},
\end{eqnarray*}
from (\ref{eq:GEM}).
For arbitrary $\overline{i}, \overline{j} \in \overline{J}$ and $k, l \in J$, we have
\begin{eqnarray}
 \left[ (\mathrm{id}^{\otimes t} \otimes \mathrm{str}_{q}) x \check{R}_{t} y \right]^{\overline{i} k}_{\overline{j} l}   
    & = & \sum_{p \in J} (-1)^{[v^{p}]} \langle \left. v^{p} \right.^{*}, K_{2 \rho} v^{p} \rangle
          \left[ (x \otimes \mathrm{id}_{V}) \check{R}_{t} (y \otimes \mathrm{id}_{V}) 
          \right]^{\overline{i} k p}_{\overline{j} l p}    \nonumber   \\
    & = & \sum_{p \in J} (-1)^{[v^{p}]} \langle \left. v^{p} \right.^{*}, K_{2 \rho} v^{p} \rangle
          \sum_{\overline{a}, \overline{d} \in \overline{J}, b, c \in J}
          X^{\overline{i} k p}_{\overline{a} b p} 
          R^{\overline{a} b p}_{\overline{d} c p}
          Y^{\overline{d} c p}_{\overline{j} l p}  \nonumber  \\
    & = & \sum_{\overline{a}, \overline{d} \in \overline{J}, b, c \in J}
          X^{\overline{i} k }_{\overline{a} b}
          \left(\sum_{p} (-1)^{[v^{p}]} \langle \left. v^{p} \right.^{*}, K_{2 \rho} v^{p} \rangle 
             R^{\overline{a} b p}_{\overline{d} c p}  \right)
          Y^{\overline{d} c}_{\overline{j} l}. \nonumber \\
 & &  \label{eq:GEM4}
\end{eqnarray}
The sum over $p$ inside the brackets in (\ref{eq:GEM4}) equals 
$\delta_{ \overline{a} \overline{d} } \delta_{bc} \chi_{V}(v^{-1})$ 
from (\ref{eq:GEM3}),
thus (\ref{eq:GEM4}) simplifies to
$$\chi_{V}(v^{-1}) \sum_{\overline{a} \in \overline{J}, b \in J}
          X^{\overline{i} k }_{\overline{a} b}
          Y^{\overline{a} b}_{\overline{j} l}
  = \chi_{V}(v^{-1}) \left[ xy \right]^{\overline{i} k }_{\overline{j} l},$$
proving the lemma.

\end{proof}

\end{subsection}

\end{section}

\section{A representation of the Hecke algebra $H_{t}(q)$ of type $A_{t-1}$}
\label{sec:repoftheHeckealgebra}

The structures of the centraliser algebras of tensor products of finite dimensional irreducible
spinorial representations of $U_{q}(\mathfrak{g})$ is not known.
Musson and Zou presented branching rules for tensoring a finite dimensional irreducible
$U_{q}(\mathfrak{g})$-module with a spinorial representation in \cite[sec. 5]{mussonzou}, but a complete
understanding is lacking.

In this section, we make an observation that has not, to the author's knowledge, appeared in the 
literature. 
Let $V_{1/2}^{+}$ be the two-dimensional irreducible representation of $U_{q}(osp(1|2))$ with even highest
weight vector, then there exists
a representation of $H_{t}(-q)$ in 
$\mathrm{End}_{U_{q}(osp(1|2))}\left[\big(V^{+}_{1/2}\big)^{\otimes t} \right]$.

\begin{subsection}{$U_{q}(osp(1|2))$ and its finite-dimensional irreducible representations}
We fix $\mathfrak{g}=osp(1|2)$ in this section.
Recall the definition of $U_{q}(\mathfrak{g})$.
The generators $e$ and $f$ are odd and $K^{\pm 1}$ are even.

Let $v_{+}$ be a non-zero vector satisfying
$$e v_{+} = 0, \hspace{10mm} K v_{+} = \omega v_{+}, \hspace{10mm}
\omega \in \mathbb{C}.$$
We assume in the rest of this section that $q$ is generic.
It is not difficult to prove that
the highest weight vector $v_{+}$ generates a finite-dimensional irreducible highest weight 
$U_{q}(\mathfrak{g})$-module $V(\omega)$ if and only if
$$\omega = \left\{ \begin{array}{ll}
\pm q^{\lambda}, & \hspace{5mm} \lambda \in \mathbb{Z}_{+}, \\
\pm i q^{\lambda}, & \hspace{5mm} \lambda \in \mathbb{Z}_{+} + 1/2.
\end{array} \right. $$
In both of these cases, $V(\omega)$ is $(2\lambda + 1)$-dimensional and we
label $V(\omega)$ by $V_{\lambda}$. To be more precise, we write 
$V_{\lambda}^{+}$ (resp. $V_{\lambda}^{-}$) where $K v_{+} = q^{\lambda}$ or 
$K v_{+} = i q^{\lambda}$ 
(resp. $K v_{+} = -q^{\lambda}$ or $K v_{+} = -i q^{\lambda}$).

In contrast with the rest of this paper, we
 take the grading of the highest weight vectors of $V_{\lambda}^{\pm}$ to be even in this section,
and denote by $\Pi V^{\pm 1}_{\lambda}$
 the corresponding irreducible $(2\lambda + 1)$-dimensional
$U_{q}(\mathfrak{g})$-modules with odd highest weight vectors.

The tensorial (odd-dim.) irreducible representations of 
$U_{q}(\mathfrak{g})$ are deformations of irreducible representations of 
$U(\mathfrak{g})$. 
The spinorial representations are not deformations of representations of $U(\mathfrak{g})$ as
$U(\mathfrak{g})$ admits no 
irreducible even-dimensional representations.

It is easy to prove the following branching rules:
\begin{lemma}
\begin{equation}
\label{eq:branchingrulesaussiesrules}
V_{\lambda}^{+} \otimes V_{1/2}^{+} \cong \left\{
\begin{array}{ll} 
V^{+}_{\lambda + 1/2} \oplus \Pi V^{+}_{\lambda - 1/2}, & 
\mbox{if } \lambda \in \{1,2,\ldots\}, \\
V^{-}_{\lambda + 1/2} \oplus \Pi V^{-}_{\lambda - 1/2}, & 
\mbox{if } \lambda \in \{1/2,3/2,5/2,\ldots\}.
\end{array} \right.
\end{equation}
\end{lemma}
We often employ the two-dimensional irreducible 
$U_{q}(\mathfrak{g})$-module $V^{+}_{1/2}$ in this section.
From (\ref{eq:branchingrulesaussiesrules}) it is convenient to define a Bratteli diagram for
$\big(V^{+}_{1/2}\big)^{\otimes t}$ similar to the Bratteli diagram for tensor products
of the two-dimensional irreducible representation of $U_{q}(sl(2))$:
let the $t^{th}$ level of the Bratteli diagram for $\big(V^{+}_{1/2}\big)^{\otimes t}$ consist
of all Young diagrams $\mu$ with $t$ boxes arranged in at most two rows: 
$\mu = [\mu_{1}, \mu_{2}]$, then $\mu$ labels
the highest weight of an
irreducible $U_{q}(\mathfrak{g})$-module $V_{(\mu_{1}-\mu_{2})/2}^{\pm}$.
A vertex $\mu$ on the $m^{th}$ row is connected 
to a vertex $\nu$ on the $(m+1)^{st}$ row if and only if $\mu$ and $\nu$ differ by exactly one box.

It is interesting to observe from (\ref{eq:branchingrulesaussiesrules}) that for each fixed $t$,
all the irreducible $U_{q}(\mathfrak{g})$-summands of 
$\big(V^{+}_{1/2}\big)^{\otimes t}$ are labelled with the same superscript (`$+$' or `$-$'); and that the superscript is `$+$' if $t \equiv 0, 1 \pmod{4}$ and `$-$' otherwise.
Furthermore, if $V_{\lambda}^{\pm}$ (resp. $\Pi V_{\lambda}^{\pm}$)
is a summand in $\big(V^{+}_{1/2}\big)^{\otimes t}$, then
$\Pi V_{\lambda}^{\pm}$ (resp. $V_{\lambda}^{\pm}$) is not a summand
in $\big(V^{+}_{1/2}\big)^{\otimes t}$.

\end{subsection}

\subsection{$R$-matrix for $V^{+}_{1/2} \otimes V^{+}_{1/2}$}
\label{sec:R-matrices}

We now write down the $R$-matrix for 
$\mathrm{End}_{\mathbb{C}}(V^{+}_{1/2} \otimes V^{+}_{1/2})$ and thereby
obtain a representation of the Braid group $B_{t}$ on $t$ strings in
$\mathrm{End}_{U_{q}(\mathfrak{g})}\left[\big(V_{1/2}^{+}\big)^{\otimes t}\right]$, which also
gives a representation of the Hecke algebra $H_{t}(-q)$ of type $A_{t-1}$.

To write down an $R$-matrix for 
$\mathrm{End}_{\mathbb{C}}(V^{+}_{1/2} \otimes V^{+}_{1/2})$ we need to slightly change 
the method we used to write down $R$-matrices for representations of 
$U_{q}(osp(1|2n))$ previously. 
In particular, we take the highest weight of $V^{+}_{1/2}$ to be complex:
write $q = e^{i \phi}$, $\phi \in \mathbb{C}$ and let $v_{+}$ be 
the highest weight vector of $V^{+}_{1/2}$, then
$$K v_{+} = iq^{1/2} v_{+} = q^{1/2 + \pi/2\phi} v_{+}.$$
Let $\{ v_{1/2}, v_{-1/2}  \}$ be a homogeneous basis of $V^{+}_{1/2}$ where $v_{1/2} = v_{+}$ and
$v_{-1/2} = f v_{+}$.
Noting that $V^{+}_{1/2}$ is two-dimensional,
let $$\widetilde{R}^{V} = 1 \otimes 1 + (q^{-1}-q)(e \otimes f)$$
and define ${\cal{E}} \in \mathrm{End}_{\mathbb{C}}(V^{+}_{1/2} \otimes V^{+}_{1/2})$ by
$${\cal{E}}(v_{j} \otimes v_{k}) = q^{(j + \pi/2\phi)(k + \pi/2\phi)} (v_{j} \otimes v_{k}), 
\hspace{5mm} j, k \in \{ 1/2, -1/2\},$$
then upon writing $\pi$ to denote the representation of $U_{q}(\mathfrak{g})$ afforded by 
$V^{+}_{1/2}$ and writing 
${\cal{R}}_{\frac{1}{2} \frac{1}{2}} = {\cal{E}} \cdot (\pi \otimes \pi)\widetilde{R}^{V}$, it can be shown, in a similar way as we did previously for tensorial representations, that
$${\cal{R}}_{\frac{1}{2} \frac{1}{2}} \cdot \left[(\pi \otimes \pi) \Delta(x)\right]
 = \left[(\pi \otimes \pi) \Delta'(x)\right] \cdot {\cal{R}}_{\frac{1}{2} \frac{1}{2}},
 \hspace{5mm} \forall x \in U_{q}(\mathfrak{g}).$$
A similar process can be followed to create $R$-matrices for tensor products of other irreducible spinorial representations of $U_{q}(\mathfrak{g})$.

A standard argument shows that the even map 
$\check{\cal{R}}_{\frac{1}{2} \frac{1}{2}} \in \mathrm{End}_{\mathbb{C}}
(V^{+}_{1/2} \otimes V^{+}_{1/2})$ 
defined by
$$\check{\cal{R}}_{\frac{1}{2} \frac{1}{2}}(v_{j} \otimes v_{k}) = 
 P \circ \left( {\cal{R}}_{\frac{1}{2} \frac{1}{2}}(v_{j} \otimes v_{k}) \right), 
\hspace{5mm} j, k \in \{ 1/2, -1/2\},$$
where $P$ is the graded permutation operator,
is $U_{q}(\mathfrak{g})$-linear:
$$\check{\cal{R}}_{\frac{1}{2}\frac{1}{2}} \cdot \left[( \pi \otimes \pi) \Delta(x) \right]
  = \left[( \pi \otimes \pi) \Delta(x) \right] \cdot \check{\cal{R}}_{\frac{1}{2}\frac{1}{2}}, 
\hspace{10mm} \forall x \in U_{q}(\mathfrak{g}).$$
With respect to the ordered basis 
$\{v_{1/2} \otimes v_{1/2}, v_{1/2} \otimes v_{-1/2}, v_{-1/2} \otimes v_{1/2}, 
v_{-1/2} \otimes v_{-1/2} \}$, 
$\check{\cal{R}}_{\frac{1}{2} \frac{1}{2}}$ is explicitly
$$\check{\cal{R}}_{\frac{1}{2} \frac{1}{2}} = i^{\pi/2\phi}\left( \begin{array}{cccc}
i q^{1/4} & 0 & 0 & 0 \\
0 & 0 & q^{-1/4} & 0 \\
0 & q^{-1/4} & i (q^{1/4} + q^{-3/4}) & 0 \\
0 & 0 & 0 & i q^{1/4}
\end{array}
\right).$$

Recall from (\ref{eq:branchingrulesaussiesrules}) that
$V^{+}_{1/2} \otimes V^{+}_{1/2} \cong V^{-}_{1} \oplus \Pi V^{-}_{0}$, then a direct 
calculation gives
$$\check{\cal{R}}_{\frac{1}{2} \frac{1}{2}} \cdot w_{1} = i^{1 + \pi/2\phi}q^{1/4} w_{1},
\hspace{10mm} \forall w_{1} \in V^{-}_{1} \subset V^{+}_{1/2} \otimes V^{+}_{1/2},$$
$$\check{\cal{R}}_{\frac{1}{2} \frac{1}{2}} \cdot w_{0} = i^{1 + \pi/2\phi}q^{-3/4} w_{0},
\hspace{10mm} \forall w_{0} \in \Pi V^{-}_{0} \subset V^{+}_{1/2} \otimes V^{+}_{1/2},$$
thus
$$\left(\check{\cal{R}}_{\frac{1}{2} \frac{1}{2}} - i^{1 + \pi/2\phi}q^{1/4}\right)
  \left(\check{\cal{R}}_{\frac{1}{2} \frac{1}{2}} - i^{1 + \pi/2\phi}q^{-3/4}\right) = 0.$$
By inspection,
\begin{equation}
\label{partialqsupertrace}
(\mathrm{id} \otimes \mathrm{str})\left[(\mathrm{id} \otimes \pi(K)) 
\check{\cal{R}}^{\pm 1}_{\frac{1}{2} \frac{1}{2}} \right] = \chi_{V^{+}_{1/2}}(v^{\mp 1}) \pi(\mathrm{id}),
\end{equation}
where $\chi_{V^{+}_{1/2}}(v^{-1}) = -i^{\pi/2\phi} q^{3/4}$
and  $\chi_{V^{+}_{1/2}}(v) = -i^{-\pi/2\phi} q^{-3/4}$.

Defining 
$$\check{R}_{j}^{\pm 1} = \mathrm{id}^{\otimes (j-1)} \otimes 
\check{\cal{R}}_{\frac{1}{2} \frac{1}{2}}^{\pm 1} \otimes
\mathrm{id}^{\otimes (t-j-1)} \in \mathrm{End}_{U_{q}(\mathfrak{g})}
\left[\big(V^{+}_{1/2}\big)^{\otimes t}\right],$$
and fixing ${\cal{C}}_{t}^{\frac{1}{2}}$ to be the complex algebra generated by
$\left\{\check{R}_{j}^{\pm 1} | \ j=1, 2, \ldots, t-1 \right\}$,
we obtain a representation of the Braid group $B_{t}$
via the homomorphism
$\overline{\rho}: B_{t} \rightarrow {\cal{C}}_{t}^{\frac{1}{2}}$ defined by
$\overline{\rho}: \sigma_{j}^{\pm 1} \mapsto \check{R}_{j}^{\pm 1}$.
The map 
$$\psi(X) = \frac{\mathrm{str}_{q}(X)}{(sdim_{q}(V^{+}_{1/2}))^{t}}, \hspace{5mm}
X \in {\cal{C}}^{\frac{1}{2}}_{t}$$ is a Markov trace and
it should be possible to construct link invariants following \cite{gtb} 
 using this Markov trace, Eq. (\ref{partialqsupertrace}) and 
Lemma \ref{ACER1} adapted to this section.

\subsection{A representation of the Hecke algebra $H_{t}(q)$}
\label{sec:HeckeAlgebra}

Following \cite{w88} we define the Hecke algebra $H_{t}(q)$ of type $A_{t-1}$
to be the complex associative unital algebra generated by the elements
$\{g_{1}, g_{2}, \ldots, g_{t-1}\}$ subject to the relations
\begin{eqnarray*}
g_{j} g_{j+1} g_{j} & = & g_{j+1} g_{j} g_{j+1}, \hspace{10mm}
                          j = 1, 2, \ldots, t-1,   \\
g_{j}g_{k} & = & g_{k}g_{j}, \hspace{20mm} |j-k| > 1,   \\
(g_{j})^{2} & = & (q-1)g_{j} + q, \hspace{7mm}
                          j = 1, 2, \ldots, t-1,  
\end{eqnarray*}
where $q$ is a complex number. Then, we immediately have the following result:
\begin{lemma}
\label{lem:repofheckealgebra}
The algebra homomorphism  
$\rho: H_{t}(-q) \rightarrow {\cal{C}}_{t}^{\frac{1}{2}}$ defined by
$$\rho: g_{j} \mapsto  i^{1-\pi/2\phi} q^{3/4} \check{R}_{j},$$
yields a representation of $H_{t}(-q)$. 
\end{lemma}
It may be possible to construct projectors and intertwiners in
$\mathrm{End}_{U_{q}(\mathfrak{g})}\left[ \big(V^{+}_{1/2}\big)^{\otimes t} \right]$ by applying
$\rho$ to matrix units (eg those given in \cite{tw}) in $H_{t}(-q)$, but we have not gone through the details. 
As with the representations 
of the Birman-Wenzl-Murakami algebra we defined 
earlier in this paper, we expect that the representation
of $H_{t}(-q)$ given in Lemma \ref{lem:repofheckealgebra} is not faithful.

Furthermore, 
even if all the projectors onto and intertwiners between the irreducible $U_{q}(\mathfrak{g})$-summands of
$\big(V^{+}_{1/2}\big)^{\otimes t}$ can be obtained by applying
$\rho$ to matrix units in $H_{t}(-q)$, 
the fact that $V^{+}_{1/2}$ has vanishing superdimension may result in
${\cal{C}}_{t}^{\frac{1}{2}}$ being a proper subalgebra of
$\mathrm{End}_{U_{q}(\mathfrak{g})}\left[ \big(V^{+}_{1/2} \big)^{\otimes t} \right]$ from Schur's lemma. 
However, we have not explored this issue and leave it for further study.

\noindent
{\bf{Acknowledgements}}

I would like to thanks to my Ph.D supervisor, Assoc. Prof. R. B. Zhang, for many useful discussions on this work.
This work was supported in part by an Australian Postgraduate Award, which was gratefully received.

\end{document}